\documentclass[final, 3p, times]{elsarticle}

\usepackage{amssymb}
\usepackage{amsmath}
\usepackage{mathtools}
\usepackage{diagrams}
\usepackage{multicol}
\usepackage{color}
\usepackage{subfigure}
\diagramstyle[labelstyle=\scriptstyle]
\newarrow{Dto}{<}{-}{-}{-}{>}
\usepackage[all]{xy}
\diagramstyle[labelstyle=\scriptstyle]
\input xy \xyoption{all}
\newcommand{\vect}[2]{\ensuremath{\boldsymbol{\mathrm{#1}}_{#2}}}

\newcommand{\Dform}[2]{\ensuremath{#1^{#2}}}

\newcommand{\chain}[2]{\ensuremath{#1_{(#2)}}}
\newcommand{\cochain}[2]{\ensuremath{#1^{(#2)}}}

\newcommand{\boundary}{\ensuremath{\partial}}
\newcommand{\coboundary}{\ensuremath{\delta}}

\newcommand{\x}[1]{\ensuremath{x^{#1}}}
\newcommand{\y}[1]{\ensuremath{y^{#1}}}
\newcommand{\dx}[1]{\ensuremath{\mathrm{d}x^{#1}}}
\newcommand{\dy}[1]{\ensuremath{\mathrm{d}y^{#1}}}

\newcommand{\diff}{\ensuremath{\mathrm{d}}}
\newcommand{\grad}{\ensuremath{\mathrm{grad}}}
\renewcommand{\div}{\ensuremath{\mathrm{div}}}
\newcommand{\curl}{\ensuremath{\mathrm{curl}}}

\newcommand{\Int}[1]{\ensuremath{\mathcal{I}_{#1}}}
\newcommand{\Red}[1]{\ensuremath{\mathcal{R}_{#1}}}

\newcommand{\Mat}[2]{\ensuremath{\mathrm{#1}_{#2}}}

\newcommand{\KnotVector}[3]{\ensuremath{\vect{#1}{#2} = \left\{ \right. #3 \left.\right\}}}

\newcommand{\Span}[1]{\ensuremath{\mathrm{span} \left\{#1\right\}}}

\newtheorem{example}{Example}[section]

\journal{}

\begin{document}

\begin{frontmatter}

\title{High order gradient, curl and divergence conforming spaces, with an application to compatible NURBS-based IsoGeometric Analysis}

\author[3ME]{R.R. Hiemstra}
\ead{R.R.Hiemstra@tudelft.nl}
\author[3ME]{R.H.M. Huijsmans}
\ead{R.H.M.Huijsmans@tudelft.nl}
\author[Aero]{M.I.Gerritsma}
\ead{M.I.Gerritsma@tudelft.nl}

\address[3ME]{Department of Marine Technology, Mekelweg 2, 2628CD Delft}
\address[Aero]{Department of Aerospace Technology, Kluyverweg 2, 2629HT Delft}

\begin{abstract}
Conservation laws, in for example, electromagnetism, solid and fluid mechanics, allow an exact discrete representation in terms of line, surface and volume integrals. We develop high order interpolants, from any basis that is a partition of unity, that satisfy these integral relations exactly, at cell level. The resulting gradient, curl and divergence conforming spaces have the property that the conservation laws become completely independent of the basis functions. This means that the conservation laws are exactly satisfied even on curved meshes. As an example, we develop high order gradient, curl and divergence conforming spaces from NURBS - non uniform rational B-splines - and thereby generalize the compatible spaces of B-splines developed in \cite{Buffa:2011b}. We give several examples of 2D Stokes flow calculations which result, amongst others, in a point wise divergence free velocity field.
\end{abstract}


\end{frontmatter}
\begin{quotation}
	\textit{Be careful of the naive view that a physical law is a mathematical relation between previously defined quantities. The situation is, rather, that a certain mathematical structure represents a given physical structure. \citet{Burke:1985}}
\end{quotation}

\section{Introduction \label{sec:Introduction}}
\begin{figure}
		\centering
		\subfigure[$ \int_L \grad \; T \; dl = T(b)-T(a)$]{\includegraphics[trim = 0cm 0cm 12cm 0cm, clip,width =0.3\textwidth]{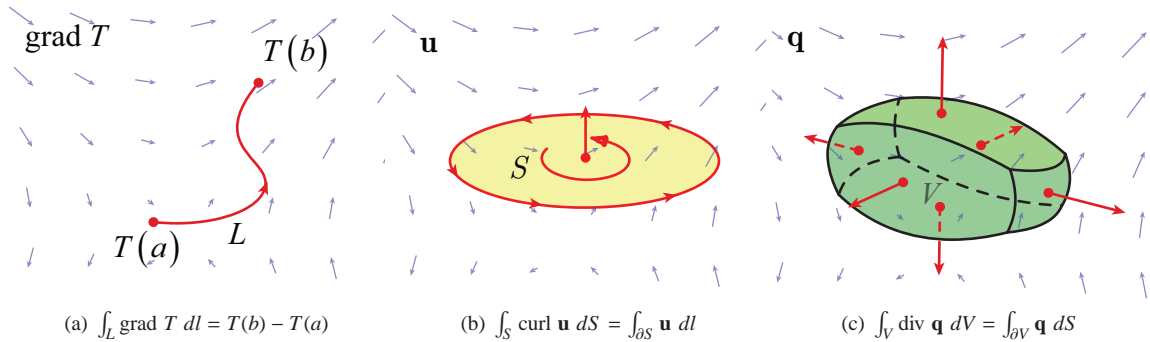}}
		\subfigure[$ \int_S \curl \;\vect{u}{} \; dS = \int_{\partial S} \vect{u}{} \; dl$]{\includegraphics[trim = 6cm 0cm 6cm 0cm, clip,width =0.3\textwidth]{figures/Stokes}}
		\subfigure[$ \int_V \div \; \vect{q}{} \;dV = \int_{\partial V} \vect{q}{} \; dS$]{\includegraphics[trim = 12cm 0cm 0cm 0cm, clip,width =0.3\textwidth]{figures/Stokes}}
		\caption{The fundamental theorem of calculus, Stokes circulation theorem and Gauss divergence theorem.}
		\label{fig:Stokes}
\end{figure}

In deriving mathematical models for physical theories, we frequently start with analysis on finite dimensional geometric objects, like a control volume and its bounding surfaces. We assign global, 'measurable', quantities to these different geometric objects and set up balance statements. Take for example the global balance in (\ref{eq:LocalVsGlobal}), where the total mass/momentum/energy $E$ inside a control volume $V$ is only conserved (no change in time) if the in- and outgoing mass/momentum/energy fluxes $Q$ over the bounding surfaces $\partial V$ cancel.
\begin{align}
	\begin{diagram}
  & \quad \frac{\partial}{\partial t} E \left(V\right) = Q \left(\partial V \right) \quad && \rTo^{\scriptstyle{\lim V \rightarrow P}} && \quad \frac{\partial}{\partial t} e = \div \; \vect{q}{}. \quad & \\
  & \text{(discrete or global)}  																							&&																					 && \text{(differential or local)}  										&  	
	\end{diagram}
	\label{eq:LocalVsGlobal}
\end{align}
This is exactly Gauss divergence theorem, depicted in Figure \ref{fig:Stokes}c. Other balance equations in $\mathbb{R}^3$ involve the fundamental theorem of calculus, relating a global quantity associated with a curve $L$, to the values of a quantity at the boundary points $\boundary L$, and Stokes circulation theorem, which relates the amount of rotation in a surface $S$ to the amount of circulation around the bounding curve $\boundary S$.

While the association of physical quantities with geometry is clear in the global sense, it remains obscured when the mathematical model is written in local form, in (\ref{eq:LocalVsGlobal}), as a differential equation. The local variables, i.e. the source field $\vect{q}{}$ and density $e$, obtained from a limiting process by shrinking the integration domain $V$ up to a point $P$, although mathematically well defined, seem to have lost their geometric significance.

Classical numerical methods, in particular finite difference and nodal finite element methods, take the differential statement as a starting point for discretization. These methods expand their unknowns in terms of nodal interpolations only, and thereby disregard the underlying geometry of the physics. This can lead to instabilities, and perhaps more dangerously, to internal inconsistencies such as the violation of fundamental conservation principles. Where instabilities lead to outright failure of a numerical method, inconsistencies can lead to unphysical solutions that can go unnoticed by the human eye \cite{Perot:2011}. This becomes more and more pronounced, as the trend is to simulate increasingly larger and complex non-linear phenomena, like multi-phase flows, fluid structure interaction and magnetohydrodynamics. 

To capture the behavior of a physical phenomena well, a discretization method should't only approximate the spaces of the infinite dimensional system, but should also follow the structure induced by the relations between them; in particular the structure induced by the fundamental balance equations depicted in Figure \ref{fig:Stokes}. By choosing degrees of freedom, not solely associated with nodes, but also with edges, faces and volumes in the mesh, we are able to exactly satisfy these relations in the discrete setting.

We follow the pioneering work of \citet{Tonti:1975,Tonti:1976}, \citet{Mattiussi:1997}, \citet{Bossavit:1998,Bossavit:2000} and the recent advances in \textit{Discrete Exterior Calculus} \cite{Desbrun:2005,Hirani:2003}, \textit{Finite Element Exterior Calculus} \cite{Arnold:2006,Arnold:2010}, \textit{Compatible} \cite{Buffa:2011b,Buffa:2010,Evans:2012,Back:2011,Ratnani:2012} and \textit{Mimetic Methods} \cite{Bochev:2010,Bochev:2006,Bouman:2011,Gerritsma:2011,Kreeft:2012a,Kreeft:2012b,Kreeft:2011,Palha:2011}. These methods, do not focus on one particular physical problem, but identify and discretize the underlying structure that constitutes a wide variety of physical field theories. They are said to be '\textit{compatible}' with the geometric structure of the underlying physics, i.e. they '\textit{mimic}' important properties of the physical system. This leads, amongst others, to naturally stable and consistent numerical schemes that have discrete conservation properties by construction and are applicable to a wide variety of physical theories. Furthermore, they offer insight into the properties of existing numerical schemes.

We develop arbitrary order interpolants, starting from any basis that is a partition of unity, which satisfy the fundamental integral theorems exactly. The gradient, curl and divergence conforming spaces have the property that the conservation laws become completely independent of the basis functions. This means that the conservation laws are exactly satisfied at the coarsest level of discretization and on arbitrarily curved meshes. It is remarkable that inf-sup stability is automatically guaranteed when this physical structure is encoded in the discretization \cite{Kreeft:2012b,Kreeft:2011}.

As an example we apply our approach to NURBS (non-uniform rational B-splines)  and thereby generalize the compatible spaces of B-splines introduced in \cite{Buffa:2011b}. NURBS, the standard in Computer Aided Design (CAD), have only recently become popular in the finite element community due to the very succesfull IsoGeometric Analysis (IGA) paradigm \cite{Hughes:2005}. IGA employs CAD technologies, such as B-splines and NURBS, directly within the finite element spaces, and thereby integrates computer aided design with finite element analysis (FEA). IsoGeometric Analysis not only bridges the gap between CAD and FEA, it also provides exact geometry representation at the coarsest level of discretization \cite{Hughes:2005}; possesses increased robustness and accuracy \cite{Lipton:2010}; and refining strategies become practically applicable \cite{Hughes:2005}. Development is ongoing and IsoGeometric Analysis has over the years reached a certain level of maturity. For an overview we refer the reader to \cite{Cottrell:2009}.

\subsection{Geometric structure of partial differential equations (PDE's)}
\begin{figure}
	\centering
	\includegraphics[width=0.7\textwidth]{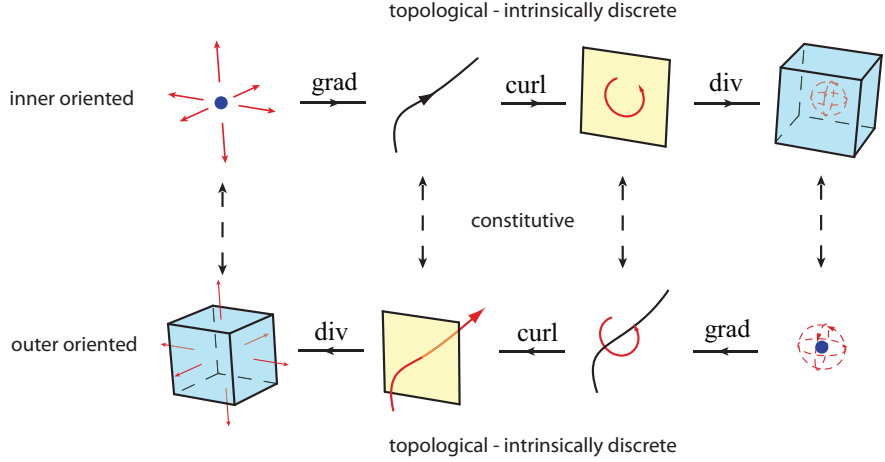}
	\caption{Physical quantities in $\mathbb{R}^3$ can be associated with either inner or outer oriented points, curves, surfaces and volumes. They are horizontally connected by means of the balance equations, i.e. fundamental theorem of calculus, Stokes circulation theorem and Gauss divergence theorem, which are metric free. The metric and material dependent relations, i.e. the constitutive equations represent the relation between physical quantities associated with dual geometric objects.}
	\label{fig:structure}
\end{figure}

Geometry induces physics with a clear geometric structure. Consider Figure \ref{fig:structure}. Physical quantities are naturally related to geometric objects such as points, curves, surfaces and volumes. Furthermore, we can consider two separate types of orientation with respect to the geometric object: inner and outer orientation. Temperature, for example, is measured in a point; strain of a fiber along a curve; magnetic flux through a surface; mass flowing out of a volume; and an amount of rotation in a surface. Note there is room for interpretation since an amount of rotation could similarly be associated with a vortex filament.

Physical quantities are horizontally connected by the fundamental theorem of calculus, Stokes circulation theorem and Gauss divergence theorem (Figure \ref{fig:Stokes}). These balance equations can exactly be described in terms of discrete physical quantities, that could have been obtained from some measurement process, and besides are independent of the shape of the geometry; i.e. they are intrinsically discrete and topological in nature. 

The metric, notions of length, area, volume, angle, along with material properties of the medium,  are described by the constitutive equations, which require a differential formulation. These represent the relations between physical quantities associated with dual geometric objects. Consider for example Hooke's law which relates the strain of a fiber along an inner oriented curve, with the stress through an outer oriented surface. Another example is the proportional relationship between point-wise defined temperature and kinetic energy contained in a volume for a perfect gas. 

The schematic in Figure \ref{fig:structure} can serve as a template to design compatible numerical methods for partial differential equations (PDE's). By reformulating a PDE in terms of a set of first order equations one can distinguish between the topological nature of the balance laws and the metric nature of the constitutive equations. This idea is by no means new. In fact, Figure \ref{fig:structure} is known as a Tonti diagram, after Enzo Tonti. For more details, we refer the reader to the work of \citet{Tonti:1975} and \citet{Mattiussi:1997}.

We consider Stokes flow in a domain $\Omega$, filled with an incompressible fluid with constant viscosity $\nu$. Under these assumptions, Stokes flow can be described by the following equations,
\begin{align}
	& \nu \Delta \vect{u}{}  -  \grad \; p = \vect{f}{} & &\text{and}& &\div \; \vect{u}{} = 0,&
\end{align}
where $\vect{u}{}$ is the velocity, $p$ the pressure and $\vect{f}{}$ the right hand side forcing. Using the the operator splitting, $-\Delta \vect{u}{} = \curl \; \curl  \; \vect{u}{} + \grad \; \div \; \vect{u}{}$, and using the incompressibility constraint, $\div \; \vect{u}{} = 0$, we can decouple Stokes flow into the following first order equations, 
\begin{align}
	& \curl \; \vect{u}{}  = \vect{\omega}{} ,    & &\curl \; \vect{\omega}{} + \grad \; p = - \vect{f}{},& &\div \; \vect{u}{}= 0,&
	\label{eq:Stokes0}
\end{align}
where $\vect{\omega}{}$ is the vorticity. Observe that Stokes circulation theorem (Figure \ref{fig:Stokes}b) means that $\vect{u}{}$ in $\curl \; \vect{u}{} $ is associated with inner oriented curves, while Gauss divergence theorem (Figure \ref{fig:Stokes}c) implies that $\vect{u}{}$ in $\div \; \vect{u}{} $ is associated with outer oriented surfaces. The variable $\vect{u}{}$ thus has a different geometric interpretation depending on the type of balance law and needs to be treated accordingly.

It is possible to discretize variables on staggered grids of opposite orientation, as in many finite volume methods, and make an explicit relation between quantities associated with dual geometric objects. We will however follow a more finite element type of method and circumvent the use of a dual (staggered) grid by using integration by parts. While $\grad$, $\curl$ and $\div$ are the differential operators that naturally occur in the fundamental theorem of calculus, Stokes circulation theorem, and Gauss divergence theorem, we can define $\grad^{\star}$, $\curl^{\star}$ and $\div^\star$ as their formal Hilbert adjoint up to a possible boundary term \cite{Kreeft:2012a},
\begin{align}
  &	(\vect{\alpha}{},-\grad^{\star} \; \beta ) 			 = (\div \; \vect{\alpha}{}, \beta ), & 		
  & (\vect{\gamma}{}, \curl^{\star} \;\vect{\delta}{} ) = (\curl \; \vect{\gamma}{}, \vect{\delta}{} ),&
  &	(\epsilon,-\div^{\star} \; \vect{\zeta}{} ) 			 = (\grad \; \epsilon, \vect{\zeta}{} ) & 	
  \label{eq:AdjointDerivative}	
\end{align}
While the $\grad$, $\curl$ and $\div$ operators are purely topological (metric free) and allow an exact discrete representation (Section \ref{sec:DiscreteModeling}), the $\grad^{\star}$, $\curl^{\star}$ and $\div^{\star}$ are metric dependent and can only be approximated in the discrete setting.

\begin{figure}
		\centering
		\includegraphics[trim = 2.5cm 0cm 0cm 0cm, clip,scale = 0.65]{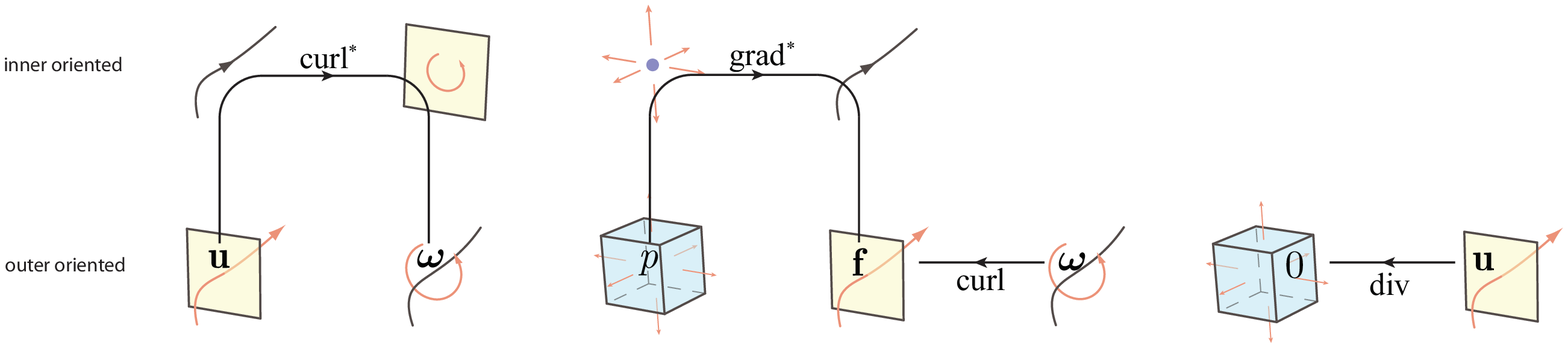}
		\caption{Tonti diagram illustrating mixed formulation for Stokes flow with outer oriented variables. Velocity $\vect{u}{}$ is associated with outer oriented surfaces, pressure $p$ is defined in outer oriented volumes and the vorticity $\vect{\omega}{}$ is related to outer oriented lines.}
		\label{diag:OuterOrientation}
\end{figure}

We can use these operators in a mixed Galerkin setting, where the resulting mixed formulation depends solely on physical considerations. Take for example the mixed formulation of the Stokes problem depicted in Figure \ref{diag:OuterOrientation}, where the variables are associated with outer oriented geometric elements. In this case we apply the topological equation $\div \; \vect{u}{} = 0$ for conservation of mass, which will eventually give a point-wise divergence free velocity field in the discrete setting. Furthermore, inflow boundary conditions can be enforced strongly, while tangential velocity is prescribed in a weak sense.

\begin{figure}
		\centering
		\includegraphics[trim = 2.5cm 0cm 0cm 0cm, clip,scale = 0.65]{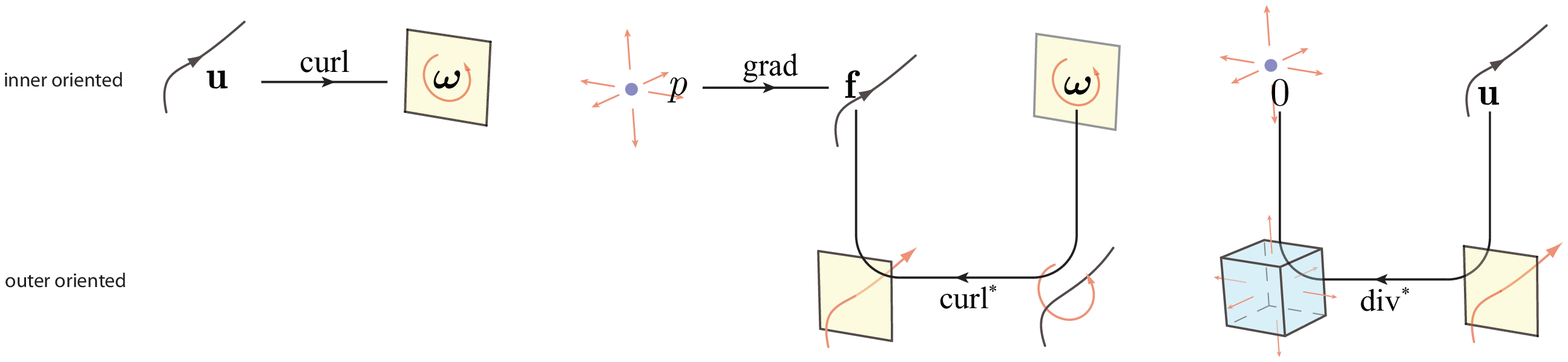}
		\caption{Tonti diagram illustrating mixed formulation for Stokes flow with inner oriented variables. Velocity $\vect{u}{}$ is associated with inner oriented lines, pressure $p$ is associated with inner oriented  points and the vorticity $\vect{\omega}{}$ is related to inner oriented surfaces.}
		\label{diag:InnerOrientation}
\end{figure}

We could also choose to associate all variables with inner oriented geometric elements. This mixed formulation is shown in Figure \ref{diag:InnerOrientation}. In this case the equation $\curl \; \vect{u}{} = \vect{\omega}{}$ is topological and allows an exact discrete representation. Furthermore, tangential velocity can be enforced strongly, while normal velocity is prescribed in a weak sense. 

\vspace{0.25cm}

The paper is structured as follows. In Section \ref{sec:DifferentialModeling} we introduce differential forms as the natural mathematical objects representing physical fields. Differential forms allow a continuous description of physics that maintains its geometric content and allows us to separate the topological from the metric dependent structure. Topological structure is intrinsically discrete, and can exactly be encoded in terms of discrete objects known as chains (discrete geometric objects) and cochains (discrete quantities), see Section \ref{sec:DiscreteModeling}. Metric structure, on the other hand, requires a continuous description. In Section \ref{sec:CommutingProjection}, we derive a new procedure to develop high order cochain interpolants, from any basis that is a partition of unity, that by construction satisfy the balance laws. As an example, we develop high order gradient, curl and divergence conforming basis functions from NURBS and apply these newly developed spaces in a mixed Galerkin setting, Section \ref{sec:MixedGalerkin}, to Stokes flow and perform some numerical calculations to asses the numerical approach (section \ref{sec:NumericalResults}); finally, conclusions are drawn in section \ref{sec:Conclusion}.

%These kind of equilibrium equations can exactly be described in terms of global (discrete) physical quantities, that could have been obtained from some measurement process, and besides are independent of the shape of the control volume; i.e. they are intrinsically discrete and topological in nature.

%The mathematical concepts of stability and consistency guarantee that the truncation error or residual of a numerical scheme tends to zero with mesh refinement. This brute force approach, that has been the cornerstone of computational science in the last 40 years, however, says nothing about the physical significance of a solution on coarse to medium-coarse grids. And this is what is desired in practice: physically significant results on coarse grids.
\section{Differential modeling \label{sec:DifferentialModeling}}
The mathematical description of physical phenomena such as electromagnetism, solid and fluid mechanics, relies heavily on the use of line, surface and volume integrals. The differential objects that appear as their integrands are called differential forms (see Example \ref{ex:DifferentialForms2}), and are studied in the mathematical field of differential geometry \cite{Burke:1985,Flanders:1963,Frankel:2011}. 

A description of physics in terms of differential forms offers significant benefits. Most notable, differential forms maintain a clear relation with the underlying geometry and therefore allow a separation of the metric dependent content, of a physical theory, from its topological (metric-free) part. This relation with geometry also makes it possible to transfer operations on the geometry to operations on the variables associated with that geometry. Furthermore, representation of variables in terms of differential forms offers a generalized concept of derivative, and thereby reduces the fundamental theorem of calculus, Stokes circulation theorem and Gauss divergence theorem to one equation, known as the generalized Stokes theorem. Most material presented in this section can be found in the references cited above.

\subsection{Differential forms}
Consider a sufficiently smooth $n$-dimensional domain $\Omega$ with boundary $\partial \Omega$ and local coordinates $\vect{x}{} = \left( \x{1},\x{2},...,\x{n} \right)$. A $k$-form $\Dform{a}{k}$ is a mathematical expression of the following form,
\begin{align}
   & \Dform{a}{k} = \sum_{I} a_I (\vect{x}{}) \; \dx{I} & & \text{with} & &\dx{I} = \dx{i_1} \wedge \dx{i_2} \wedge ... \wedge \dx{i_k},&
   \label{eq:Dform}
\end{align}
where the indices satisfy, $1 \leq i_1 \leq i_2 \leq ... \leq i_k \leq n$, the $a_I (\vect{x}{})$ are smooth functions that denote the spatial distribution of some physical quantity, and the basis elements $\dx{I}$ refer to the associated geometry. The collection of all $k$-forms is a vector space $\Lambda^{k}(\Omega)$ of dimension equal to the binomial coefficient $\binom{n}{k}$. Some examples of differential forms in $\mathbb{R}^3$ are depicted below.

\begin{example}[Differential forms in $\mathbb{R}^3$]
	Temperature, $T^{0}$, is a 0-form, a scalar function that assigns to every point $\vect{x}{}$ in domain $\Omega$ a real valued temperature; force, $	F^{1}$, is a 1-form, since it can naturally be integrated along a curve to obtain work; flux density, $Q^{2}$, is a 2-form and can be integrated over a surface to obtain the global flux; mass density $\rho^{3}$, is a 3-form, since it is readily integrated over a volume to yield mass.
\begin{align*}
		&\text{\textbf{0-form}}& & \text{Temperature:}  &	  T^{0} &= T(\vect{x}{}), \\
		&\text{\textbf{1-form}}& & \text{Force density} &    F^{1} &= f_1(\vect{x}{}) \dx{1} + f_2(\vect{x}{}) \dx{2} + f_3(\vect{x}{}) \dx{3}, \\
		&\text{\textbf{2-form}}& & \text{Flux density}  &   	Q^{2} &= q_1(\vect{x}{}) \dx{2} \wedge \dx{3} + q_2(\vect{x}{}) \dx{3} \wedge \dx{1} + q_3(\vect{x}{}) \dx{1} \wedge \dx{2}, \\
		&\text{\textbf{3-form}}& & \text{Mass density}	 & \rho^{3} &= \rho(\vect{x}{}) \dx{1} \wedge \dx{2} \wedge \dx{3}. 		
\end{align*}
\label{ex:DifferentialForms2}
\end{example}

The $\wedge$, is called the wedge product, which is a map, $ \wedge \; : \; \Lambda^{k}(\Omega) \times \Lambda^{l}(\Omega) \mapsto \Lambda^{k+l}(\Omega)$, $k+l \leq n$. The wedge product is skew symmetric, $\Dform{a}{k} \wedge \Dform{b}{l} = (-1)^{k+l} \; \Dform{b}{l} \wedge \Dform{a}{k}$. So, $\dx{i} \wedge \dx{j} = - \dx{j} \wedge \dx{i}$, which implies that the orientation changes under a permutation of the elements and linear independence, since $\dx{i} \wedge \dx{i} = -\dx{i} \wedge \dx{i} = 0$.

$k$-forms are the natural integrands over $k$-dimensional geometric objects,
\begin{align}
	\left\langle \Dform{a}{k},\Omega_k \right\rangle :=	\int_{\Omega_{k}} \Dform{a}{k} \in \mathbb{R} 
	\label{eq:DualityPairing}
\end{align}
Here $\left\langle .,.\right\rangle$ denotes duality pairing between the $k$-form and a $k$-dimensional geometric object. This duality, between geometry and physics has important consequences. It induces physics with a clear geometric structure. This allows us to separate the metric and material dependent content in a physical theory from the topological, metric independent part. Furthermore, it makes it possible to transfer operations on the geometry to operations on the variables associated with that geometry.

\subsection{Topological structure}
Integration itself is a metric free operation. If $\Phi \; : \; \Omega' \mapsto \Omega$, is a smooth map from an $n$-dimensional reference domain $\Omega'=[0,1]^n$ to the $n$-dimensional physical domain $\Omega$, the pullback,  $\Phi^{\star} \; : \; \Lambda^k(\Omega) \mapsto \Lambda^{k}(\Omega')$, maps $k$-forms on $\Omega$ to $k$-forms on $\Omega'$. This means we can perform the integration of a $k$-form in the reference domain $\Omega'$,
\begin{align}
	\int_{\Phi \left(\Omega'_{k} \right)} \Dform{a}{k} = \int_{\Omega'_{k}}  \Phi^{\star} \Dform{a}{k} \quad \Longleftrightarrow  \quad \left\langle \Dform{a}{k},\Phi\left(\Omega'_k\right) \right\rangle = \left\langle \Phi^{\star}\Dform{a}{k},\Omega'_k \right\rangle.
	\label{eq:PullBack}
\end{align}
This means that the pullback $\Phi^{\star}$ is the formal adjoint of the map $\Phi$ in the duality pairing defined in (\ref{eq:DualityPairing}). Important properties of the pull back are linearity and commutation with the wedge product, $\Phi^{\star} \left(\Dform{a}{k}\wedge \Dform{b}{l} \right) =  \left(\Phi^{\star}\Dform{a}{k}  \right) \wedge  \left(\Phi^{\star} \Dform{b}{l} \right)$.

Topological structure is induced by balance laws, which relate a quantity associated with a geometric object to another quantity which is associated with its boundary. In the theory of differential forms, this balance is represented by the generalized Stokes theorem, the mother of all equations,
\begin{align}
	\int_{\Omega_{k}} \diff \Dform{\omega}{k-1} = \int_{\partial \Omega_{k}} \Dform{\omega}{k-1} \quad \Longleftrightarrow \quad
	\left\langle \diff \Dform{\omega}{k-1},\Omega_{k}\right\rangle = \left\langle \Dform{\omega}{k-1},\boundary \Omega_{k}\right\rangle.
	\label{eq:Stokes}
\end{align}
Again we can observe the duality between geometry and physics: the boundary operator, $\boundary$, is the adjoint of the exterior derivative $\diff$. While the boundary operator is a map $\boundary \; : \; \Omega_{k} \mapsto \Omega_{k-1}$, the exterior derivative is a map, $\diff \; : \; \Lambda^{k-1}(\Omega) \mapsto \Lambda^{k}(\Omega)$. 

The exterior derivative $\diff$ is a coordinate independent generalization of the well known vector calculus identities, the gradient, curl and divergence operators in $\mathbb{R}^3$. The Stokes theorem (\ref{eq:Stokes}) thus generalizes the fundamental theorem of calculus ($k=1$), Stokes circulation theorem ($k=2$) and Gauss divergence theorem ($k=3$), depicted in Figure \ref{fig:Stokes}. Furthermore, the exterior derivative is no more difficult to compute in a curved coordinate system than it is in a Cartesian frame. One simply takes the differential of the components and follows the properties of the wedge product.
\begin{example}[Action of the exterior derivative in $\mathbb{R}^3$] Consider local coordinates $\left(\x{1},\x{2},\x{3} \right)$. The action of the exterior derivative on the 0-form $\Dform{T}{0}$, 1-form $\Dform{F}{1}$ and 2-form $\Dform{Q}{2}$ from Example \ref{ex:DifferentialForms2}, is given by
\begin{align*}
	\diff \Dform{T}{0} &= \frac{\partial T}{\partial \x{1}} \dx{1} + \frac{\partial T}{\partial \x{2}} \dx{2} + \frac{\partial T}{\partial \x{3}} \dx{3} \\
		\diff 	\Dform{F}{1} &= \left( \frac{\partial f_3}{\partial \x{2}}  - \frac{\partial f_2}{\partial \x{3}} \right) \dx{2} \wedge \dx{3} +
													  \left( \frac{\partial f_1}{\partial \x{3}}  - \frac{\partial f_3}{\partial \x{1}} \right) \dx{3} \wedge \dx{1} + 
													  \left( \frac{\partial f_2}{\partial \x{1}}  - \frac{\partial f_1}{\partial \x{2}} \right) \dx{1} \wedge \dx{2}  \\
		\diff \Dform{Q}{2} &=   \left( \frac{\partial q_1}{\partial \x{1}} + \frac{\partial q_2}{\partial \x{2}} + \frac{\partial q_2}{\partial \x{2}}   \right) \dx{1} \wedge \dx{2} \wedge \dx{3}.
\end{align*}
Observe that the components are given by those of the familiar $\grad$, $\curl$ and $\div$ operators, respectively.
\end{example}

The exterior derivative has a number of important properties: it is a linear operator (\ref{eq:diff1}); satisfies a Leibniz rule for differentiation (\ref{eq:diff2}); its metric free nature is reflected by commutation with the pull back  (\ref{eq:diff3}); and is exact  (\ref{eq:diff4});
\begin{align}
	& \diff \left( \Dform{a}{k} + \Dform{b}{k} \right) = \diff \Dform{a}{k} + \diff \Dform{b}{k}  \label{eq:diff1} \\
	& \diff \left(\Dform{a}{k} \wedge \Dform{b}{l}  \right) 
																= \diff \Dform{a}{k} \wedge \Dform{b}{l} - (-1)^{k} \Dform{a}{k} \wedge \diff \Dform{b}{l} \label{eq:diff2}\\
	& \diff \Phi^{\star} = \Phi^{\star} \diff \label{eq:diff3} \\
	& \diff \circ \diff = 0 \label{eq:diff4}
\end{align}
This last property is analogous to the vector calculus identities, $\curl \; \grad = 0$ and $\div \; \curl=0$.

The $n+1$-spaces of differential forms in an $n$-dimensional domain $\Omega$, satisfy the following sequence, known as the \textit{de Rahm} complex,
\begin{align}
\begin{diagram}
	\mathbb{R}&  \rTo  & \Lambda^{0}(\Omega) & \rTo^{\diff} & \Lambda^{1}(\Omega) & \rTo^{\diff} & ... & \rTo^{\diff} & \Lambda^{n}(\Omega) & \rTo^{\diff} & 0.	\\
\end{diagram}
\label{Diag:TopologicalLaws}
\end{align}
which is exact on contractible domains. The main focus of this paper, see Section \ref{sec:CommutingProjection}, is to construct discrete spaces of differential forms, that follow the same structure, such that conservation and balance laws can be strongly enforced.

\subsection{Metric structure}
Constitutive equations depend on the local metric, notions of length, angle, area, volume etc., and material properties of the medium under consideration. In order to measure the local metric we require a point wise inner product of $k$-forms, $\left(.,.\right) \; : \; \Lambda^{k} \left(\Omega \right) \times \Lambda^{k} \left(\Omega \right) \mapsto \mathbb{R}$. In particular, a Riemanian metric gives rise to the Hilbert space $L^2$ inner product on $\Lambda^{k}(\Omega)$,
\begin{align}
	\left( \Dform{\alpha}{k}, \Dform{\beta}{k} \right)_{\Omega} := \int_{\Omega} \left( \Dform{\alpha}{k}, \Dform{\beta}{k} \right) \diff \Omega 
																															 = \int_{\Omega} \Dform{\alpha}{k} \wedge \star \Dform{\beta}{k}.
\label{eq:hodge}
\end{align}
Here $\star$ denotes the Hodge star operator, which is a map $\star \: : \; \Lambda^{k} \left(\Omega \right) \mapsto \Lambda^{n-k} \left(\Omega \right)$.
\begin{example}[Action of the Hodge star in $\mathbb{R}^3$]
In $\mathbb{R}^3$ with orthonormal coordinates $\vect{y}{} = \left(\y{1},\y{2},\y{3}\right)$ we have
\begin{align*}
  & \star 1 = \dy{1} \wedge \dy{2} \wedge \dy{3},&  		& \star \dy{1} \wedge \dy{2} \wedge \dy{3} = 1&
\end{align*}
\begin{align*}
  & \star \dy{1} = \dy{2} \wedge \dy{3},& 	& \star \dy{2} = \dy{3} \wedge \dy{1},& & \star \dy{3} = \dy{1} \wedge \dy{2}.&
\end{align*}
\end{example}

The Hodge not only maps $k$-forms into $(n-k)$-forms, but also changes its type of orientation from inner to outer and vice versa. We can therefore apply it to connect two copies of an exact sequence into the following structure,
\begin{align}
\begin{diagram}
  \mathbb{R} & \rTo	& \Lambda^{0}(\Omega)	& \rTo^{\diff} & \Lambda^{1}(\Omega) & \rTo^{\diff} & \Lambda^{2}(\Omega) & \rTo^{\diff} & \Lambda^{3}(\Omega) & \rTo	& 0	 \\
             &      & \dDto_{\star}       &              &    \dDto_{\star}    &              &   \dDto_{\star}     &              &    \dDto_{\star}    &      &    \\
           0 & \lTo & \Lambda^{3}(\Omega) & \lTo^{\diff} & \Lambda^{2}(\Omega) & \lTo^{\diff} & \Lambda^{1}(\Omega) & \lTo^{\diff} & \Lambda^{0}(\Omega) & \lTo & \mathbb{R}
\end{diagram}
\label{eq:DoubleDeRahm}
\end{align}
Note that this is exactly the structure given in Figure \ref{fig:structure}! While the exterior derivative $\diff$ models the topological structure of the balance laws, the Hodge $\star$ models the metric structure of the constitutive equations.

Using Leibniz's rule for differentiation (\ref{eq:diff2}), and noting that a double application of the Hodge star changes nothing up to a possible minus sign, $\star \star = (-1)^{k(n-k)}$, we can derive an adjoint operator for the exterior derivative,
\begin{align*}
	\diff \left(\Dform{\alpha}{k-1} \wedge \star \Dform{\beta}{k}  \right) 
			&= \diff \Dform{\alpha}{k-1} \wedge \star \Dform{\beta}{k} - (-1)^{k} \Dform{\alpha}{k-1} \wedge \diff \star \Dform{\beta}{k} \\
			&= \diff \Dform{\alpha}{k-1} \wedge \star \Dform{\beta}{k} - (-1)^{n(k+1)+1} \Dform{\alpha}{k-1} \wedge \star \left(\star \diff \star \right) \Dform{\beta}{k},
\end{align*}
where $\diff^{\star} = (-1)^{n(k+1)+1} \star \diff \star$ is the codifferential. Using the inner product property of the Hodge (\ref{eq:hodge}) and using integration by parts, we find that the codifferential is the Hilbert adjoint of the exterior derivative with an additional boundary term,
\begin{align}
	\left( \diff \Dform{\alpha}{k-1} , \Dform{\beta}{k} \right)_{\Omega} - 	\left( \Dform{\alpha}{k-1} , \diff^{\star} \Dform{\beta}{k} \right)_{\Omega}
	= \int_{\partial \Omega} \Dform{\alpha}{k-1} \wedge \star \Dform{\beta}{k}.
	\label{eq:coderivative}
\end{align}
The coderivative $\diff^{\star}$ is thus a generalization of the $\grad^{\star}$, $\curl^{\star}$ and $\div^{\star}$ introduced in (\ref{eq:AdjointDerivative}). Equation (\ref{eq:coderivative}) will prove to be very important in the mixed Galerkin setting explained in Section \ref{sec:MixedGalerkin}. The exterior derivative and codifferential can be combined to construct the Laplacian $\Delta$. The Hodge Laplacian is a map, $\Delta \; : \; \Lambda^{k} \left( \Omega \right) \mapsto \Lambda^{k} \left( \Omega \right)$,
\begin{align}
	\Delta = \diff \diff^{\star} + \diff^{\star} \diff.
\end{align}
\section{Discrete modeling \label{sec:DiscreteModeling}}
Classical numerical methods, in particular finite difference and nodal finite element methods, expand their unknowns in terms of nodal interpolations and run into trouble when it comes to conservation. Conservation, by definition (Generalized Stokes theorem (\ref{eq:Stokes})) is a relation between a global 'measurable' quantity associated with a geometric object and another global 'measurable' quantity associated with its boundary. By choosing degrees of freedom associated with nodes, as well as edges, faces and volumes in the mesh, we are able to exactly satisfy these relations in the discrete case.

\begin{figure}
	\centering
	\includegraphics[width=0.75\textwidth]{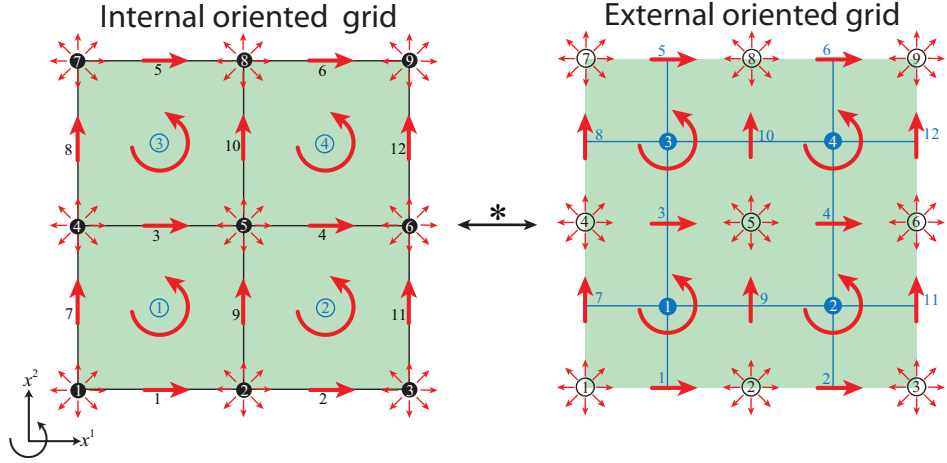}%
	\caption{The discrete Hodge invokes a connection between variables associated with dual geometric objects.}
	\label{fig:DiscreteHodge}
\end{figure}

The constitutive equations describe a relation between physical quantities associated with dual geometric objects. In discrete space, this is modeled by a discrete Hodge $\ast$ operator which invokes a relation between staggered grids of opposite orientation, see Figure \ref{fig:DiscreteHodge}. 
Staggered finite volume methods explicitly build such a dual (staggered) grid and thus explicitly construct a discrete Hodge to connect variables on both grids. In this paper we use a Galerkin finite element type of approach and circumvent the use of a dual grid by applying integration by parts, using (\ref{eq:coderivative}).

In this section we introduce geometry in terms of algebraic topology  \cite{Hatcher:2002,Singer:1976}. Topology describes the relations between oriented geometric objects (chains) and what 'lives' on those objects (cochains), however, without the notion of distance or measure. Algebraic topology thus allows us to encode the purely discrete and topological content of the physics into the discrete model, without any approximation.

Algebraic topology can in many ways be regarded as the discrete analogous of differential geometry. The duality pairing between chains and cochains is analogue to that of integration of differential forms. By defining a formal adjoint of the boundary of a chain, we obtain a discrete derivative acting on cochains, known as the coboundary operator. This discrete derivative is constructed such that it exactly satisfies Stokes theorem in terms of chains and cochains. Because all objects and operators that will be introduced in this section are topological and thus completely metric-free, they have the same form and value on topologically equivalent grids. The discretization of these relations is equivalent on a nice uniform Cartesian mesh as it is on a highly curved grid. Furthermore, they do not change on moving meshes as long as the topology stays the same.
%
%
%\textbf{We first show how to calculate with discrete geometric objects, called chains, examples of which are sets of oriented points, edges, faces and volumes. By duality pairing between geometry and physical variables, we then derive how to perform calculations on cochains, the discrete quantities associated with chains.}

\subsection{Cell complexes, chains and the boundary operator}
The pairing between physical quantities and oriented geometric objects, such as points, edges, faces and volumes, leads to a straightforward discretization of the $n$-dimensional computational domain $\Omega$. The topological description of $\Omega$ is given by a so called cell complex $D$, a formalized concept of the discretization of space, which is a partitioning in terms of $k$-dimensional sub-domains ($k=0,...,n$) called $k$-cells, which we denote by $\sigma_{(k),i}$.

\begin{figure}
	\centering
\subfigure[Invertible map]{\includegraphics[trim = 0cm 0cm 0cm 0cm, clip,scale=0.6]{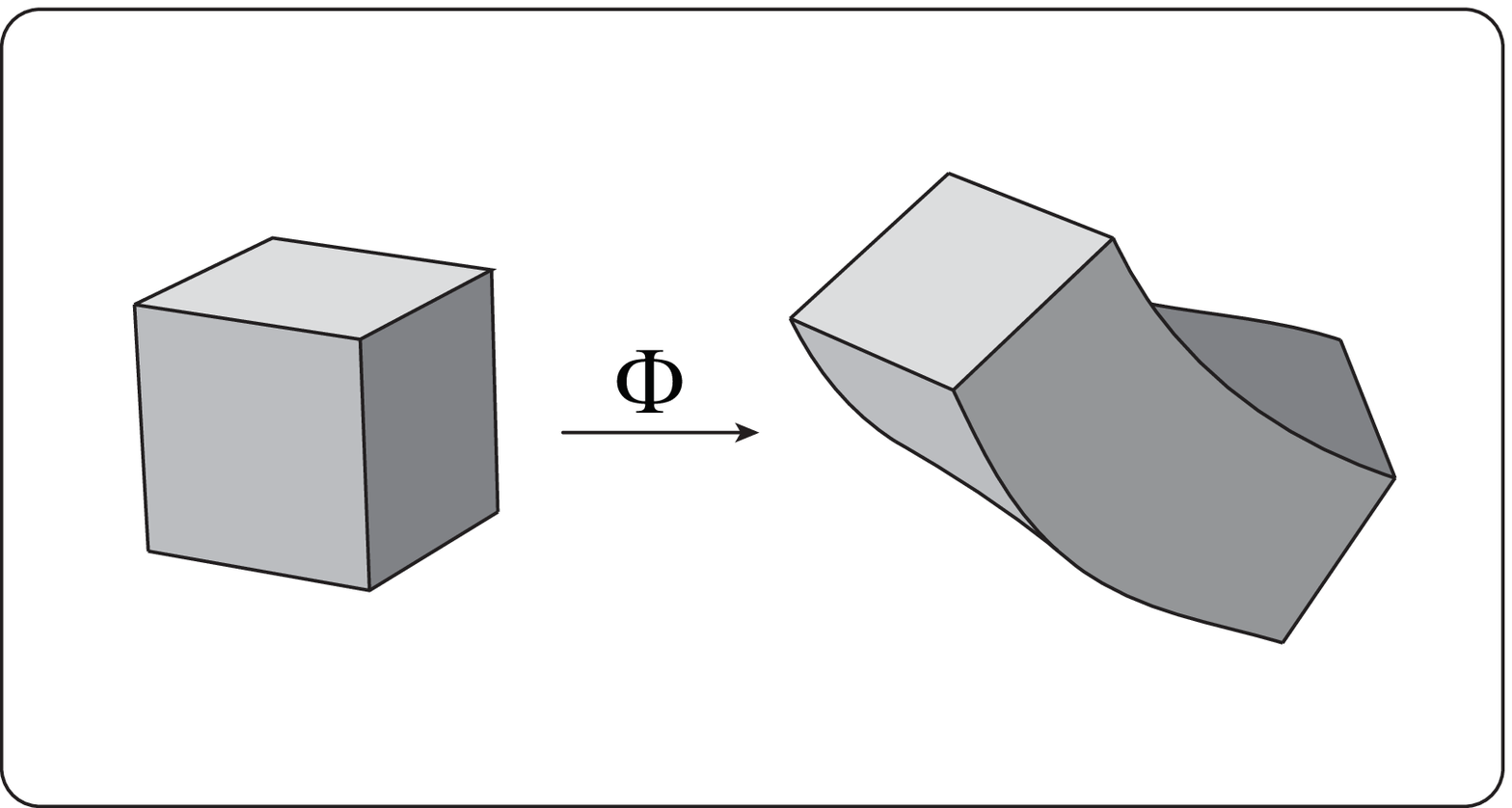}} \hfill
\subfigure[Numbering and orientation]{\includegraphics[trim = 0cm 0cm 0cm 0cm, clip,scale=0.6]{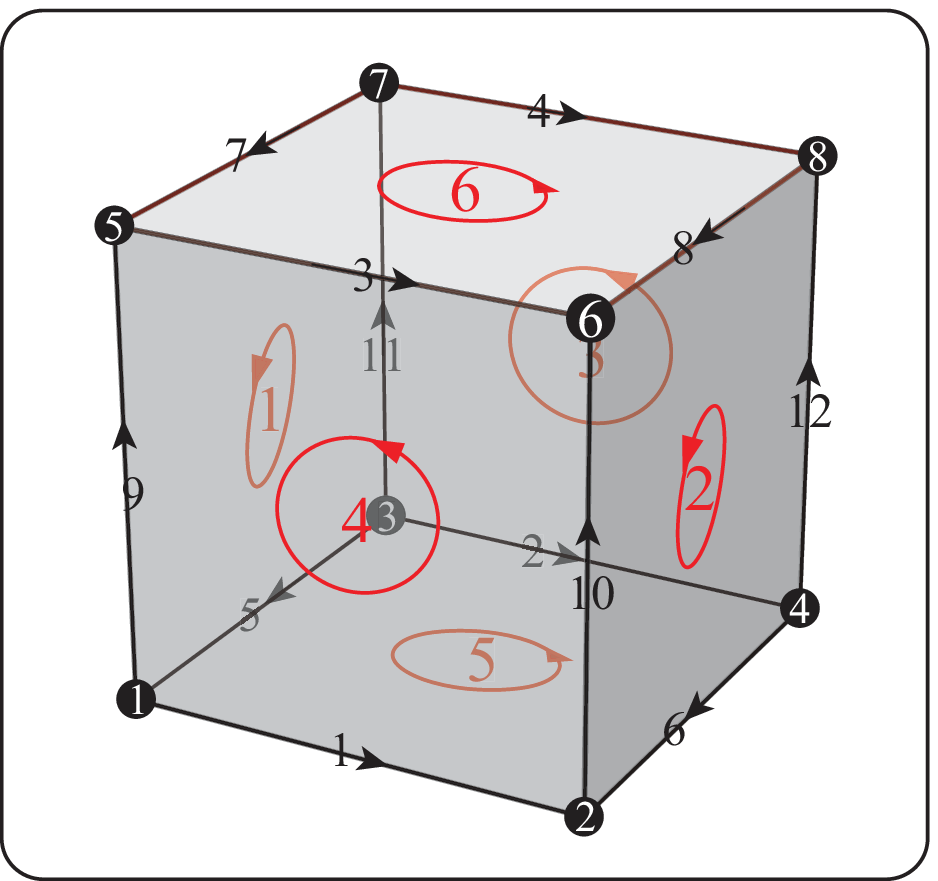}}
\subfigure[Partitioning of reference domain in terms of of 0-cells (points), 1-cells (edges), 2-cells (faces) and 3-cells (volumes)]{\includegraphics[trim = 0cm 0cm 0cm 0cm, clip,scale=0.595]{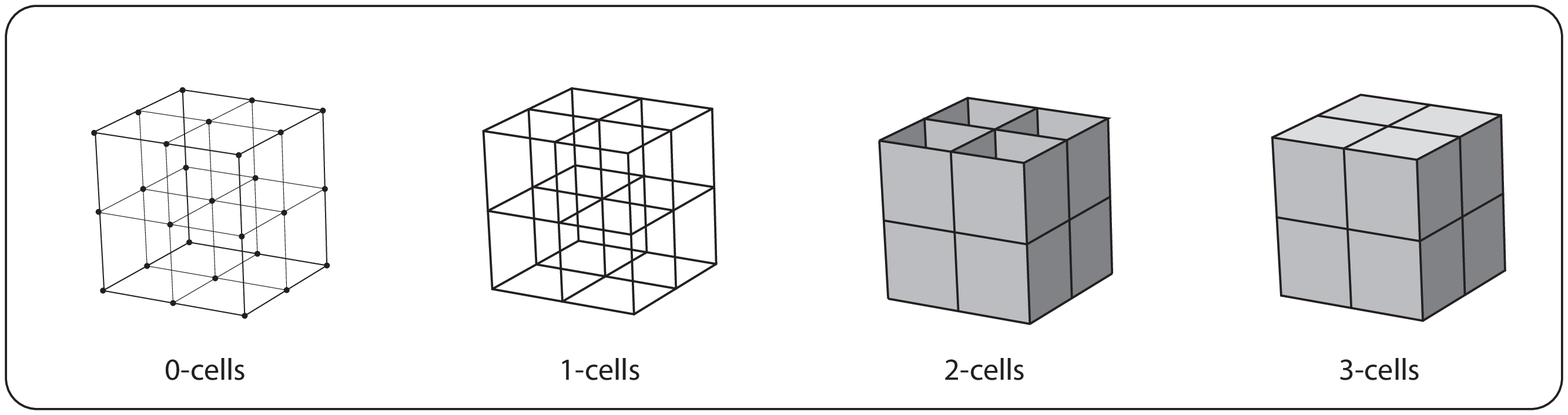}}
	\caption{(a) The curved domain $\Phi(\Omega') =\Omega$ has the same topology as the reference domain $\Omega'$. (c) Once the reference domain $\Omega'$ is partitioned into 0-cells (points), 1-cells (edges), 2-cells (faces) and 3-cells (volumes), every cell can be given an arbitrary numbering and orientation (b).}
	\label{fig:CellComplex}
\end{figure}

Figure \ref{fig:CellComplex} depicts a cell complex in $\mathbb{R}^3$, where we can distinguish between 0-cells (points), 1-cells (edges), 2-cells (faces) and 3-cells (volumes). For the present work it is most convenient to think of a cell complex as a union of $k$-cubes. An alternative would be to partition the cell complex in $k$-simplices as is done in e.g. \cite{Desbrun:2005,Hirani:2003,Arnold:2006,Arnold:2010}. From the topological viewpoint, both descriptions are equivalent \cite{Dieudonne:1989}.

Once all $k$-cells in the cell complex $D$ have been numbered and given an orientation (see Figure \ref{fig:CellComplex}b), they can be collected to form a $k$-chain, $\chain{c}{k} \in \chain{C}{k} (D)$,
\begin{align}
	\chain{c}{k} = \left\{c^j \cdot \sigma_{(k),j} \right\}_{j=0}^s \quad \text{with} \; c^{j} \in \left\{-1,0,1 \right\}.
\end{align}
In the representation of a chain we use superscript for the coefficients $c^j$ and subscript for the basis $k$-cells $\sigma_{(k),j}$. Chains can be used as a discrete representation of the geometry. In the description of geometry, however, we will only be concerned with chains with coefficients, $c^j$ of $\left\{-1,0,1 \right\}$. These correspond to either a cell with orientation opposite to the chosen default one, a cell not part of the chain, or a cell with default orientation. In actual matrix calculations we directly use the coefficients as the column vector $\vect{c}{} = \begin{pmatrix} c^0 & c^1 & \cdots & c^s \end{pmatrix}^T$.

\begin{figure}[h]
	\centering
	\includegraphics[width=0.75\textwidth]{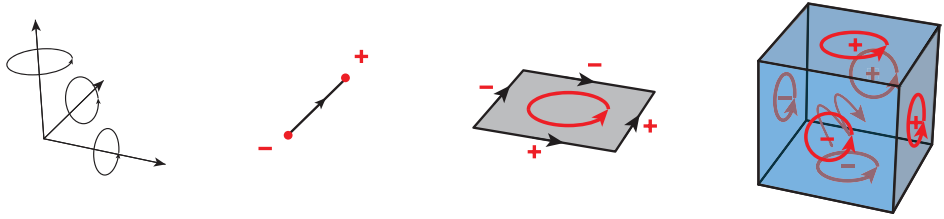}
	\caption{Compatible orientation (+) and incompatible orientation (-) between volumes and boundary faces, faces and boundary edges, edges and boundary points.}
	\label{fig:Orientation}
\end{figure}

Connectivity between volumes, faces, edges and points, is encoded in the boundary operator, $\boundary$. The boundary operator is a linear map $\boundary \; : \; \chain{C}{k} \mapsto \chain{C}{k-1}$ defined as,
\begin{align*}
	\boundary \chain{c}{k} = \boundary \left\{c^j \cdot \sigma_{(k),j} \right\}_{j=0}^s  = \left\{c^j \cdot \boundary \sigma_{(k),j} \right\}_{j=0}^s 
\end{align*}
The boundary $\boundary$ of a $k$-chain returns a unique $(k-1)$-chain and can be calculated as a linear combination of its $k$-cells boundaries. More precisely, $ \boundary \sigma_{(k),j} = \left\{ e^i_j \; \sigma_{(k-1),i} \right\}_{i=0}^r$, where $e^i_j$ has the value,
\begin{enumerate}
	\item $e^i_j = 0$, if $\sigma_{(k-1),i}$ is not part of the boundary of $\sigma_{(k),j}$.
	\item $e^i_j = 1$, if $\sigma_{(k-1),i}$ is part of the boundary of  $\sigma_{(k),j}$ with compatible orientation.
	\item $e^i_j =-1$, if $\sigma_{(k-1),i}$ is part of the boundary of  $\sigma_{(k),j}$ with incompatible orientation
\end{enumerate}
Figure \ref{fig:Orientation} illustrated what we mean with compatible and incompatible orientation.

Hence,
\begin{align*}
		\boundary \chain{c}{k} = \left\{c^j \cdot \boundary \sigma_{(k),j} \right\}_{j=0}^s = \left\{  \sum_{j=0}^s c^j e^i_j \cdot \sigma_{(k-1),i} \right\}_{i=0}^r = 
		\left\{ d^i \cdot \sigma_{(k-1),i} \right\}_{i=0}^r = \chain{d}{k-1}.
\end{align*}
So the boundary of $\chain{c}{k}$ is the unique $(k-1)$-chain, $ \chain{d}{k-1} $, of which the coefficients are given by $d^i = \sum_{j=0}^s c^j e^i_j$. The action of the boundary operator therefore allows the matrix vector product $\vect{d}{} = \mathrm{E}_{k-1,k} \vect{c}{}$. Here $\mathrm{E}_{k-1,k}$ is a  $\mathrm{rank}\left(C_{(k-1)} \right) \times \mathrm{rank}\left(C_{(k)} \right)$ incidence matrix with coefficients $\left(\mathrm{E}_{k-1,k}\right)_{ij} = e^i_j$. Note that in matrix calculations we only use the coefficients, not the basis chains. Example \ref{ex:IncidenceMatrix} shows the incidence matrices associated with the numbered and oriented cell complex of Figure \ref{fig:CellComplex}b.
\begin{example}[Incidence matrices in $\mathbb{R}^3$]
Consider the numbered and oriented cell complex in Figure \ref{fig:CellComplex}b. The connectivity between points, edges, faces and volumes  is encoded in the following incidence matrices. $\Mat{E}{0,1}$ maps from 1-chains to 0-chains, $\Mat{E}{1,2}$ maps from 2-chains to 1-chains and $\Mat{E}{2,3}$ maps from 3-chains to 2-chains.   
\setcounter{MaxMatrixCols}{12}
%\begin{small}
%	\begin{align*}
%		\Mat{E}{0,1} = 
%		\begin{bmatrix*}[r]
%		 % 1  &  2  &  3  &  4  &  5  &  6  &  7  &  8   	  
%			-1  &  1  &  0  &  0  &  0  &  0  &  0  &  0  \\
%			 0  &  0  & -1  &  1  &  0  &  0  &  0  &  0  \\
%			 0  &  0  &  0  &  0  & -1  &  1  &  0  &  0  \\
%			 0  &  0  &  0  &  0  &  0  &  0  & -1  &  1  \\
%			 1  &  0  & -1  &  0  &  0  &  0  &  0  &  0  \\
%			 0  &  1  &  0  & -1  &  0  &  0  &  0  &  0  \\
%			 0  &  0  &  0  &  0  &  1  &  0  & -1  &  0  \\
%			 0  &  0  &  0  &  0  &  0  &  1  &  0  & -1  \\
%			-1  &  0  &  0  &  0  &  1  &  0  &  0  &  0  \\
%			 0  & -1  &  0  &  0  &  0  &  1  &  0  &  0  \\
%			 0  &  0  & -1  &  0  &  0  &  0  &  1  &  0  \\
%			 0  &  0  &  0  & -1  &  0  &  0  &  0  &  1  \\
%		\end{bmatrix*}
%		, \quad
%		\Mat{E}{1,2} = 
%		\begin{bmatrix*}[r]
%		 % 1  &  2  &  3  &  4  &  5  &  6  &  7  &  8  &  9  & 10  & 11  & 12 
%			 0  &  0  &  0  &  0  & -1  &  0  &  1  &  0  & -1  &  0  &  1  &  0  \\
%			 0  &  0  &  0  &  0  &  0  & -1  &  0  &  1  &  0  & -1  &  0  &  1  \\
%			 0  &  1  &  0  & -1  &  0  &  0  &  0  &  0  &  0  &  0  & -1  &  1  \\
%			 1  &  0  & -1  &  0  &  0  &  0  &  0  &  0  & -1  &  1  &  0  &  0  \\
%			 1  & -1  &  0  &  0  &  1  & -1  &  0  &  0  &  0  &  0  &  0  &  0  \\
%			 0  &  0  &  1  & -1  &  0  &  0  &  1  & -1  &  0  &  0  &  0  &  0  \\
%		\end{bmatrix*}
%		, \quad
%		\Mat{E}{2,3} = 
%		\begin{bmatrix*}[r]
%			 -1  &  1  &  -1  &  1  &  -1  &  1
%		\end{bmatrix*}
%	\end{align*}
%\end{small}
\begin{small}
	\begin{align*}
		\Mat{E}{0,1} = 
		\begin{bmatrix*}[r]
		 % 1  &  2  &  3  &  4  &  5  &  6  &  7  &  8  &  9  & 10  & 11  & 12 	  
			-1  &  0  &  0  &  0  &  1  &  0  &  0  &  0  & -1  &  0  &	 0  &  0 \\
			 1  &  0  &  0  &  0  &  0  &  1  &  0  &  0  &  0  & -1  &	 0  &  0 \\
			 0  & -1  &  0  &  0  & -1  &  0  &  0  &  0  &  0  &  0  &	-1  &  0 \\
			 0  &  1  &  0  &  0  &  0  & -1  &  0  &  0  &  0  &  0  &	 0  & -1 \\ 
			 0  &  0  & -1  &  0  &  0  &  0  &  1  &  0  &  1  &  0  &	 0  &  0 \\
			 0  &  0  &  1  &  0  &  0  &  0  &  0  &  1  &  0  & -1  &	 0  &  0 \\
			 0  &  0  &  0  & -1  &  0  &  0  & -1  &  0  &  0  &  0  &	 1  &  0 \\
			 0  &  0  &  0  &  1  &  0  &  0  &  0  & -1  &  0  &  0  &	 0  &  1
		\end{bmatrix*}
		, \quad
		\Mat{E}{1,2} = 
		\begin{bmatrix*}[r]
		 % 1  &  2  &  3  &  4  &  5  &  6
			 0  &  0  &  0  &  1  &  1  &  0   \\
			 0  &  0  &  1  &  0  & -1  &  0   \\
			 0  &  0  &  0  & -1  &  0  &  1   \\	
			 0  &  0  & -1  &  0  &  0  & -1   \\	
			-1  &  0  &  0  &  0  &  1  &  0   \\	
			 0  & -1  &  0  &  0  & -1  &  0   \\	
			 1  &  0  &  0  &  0  &  0  &  1   \\	
			 0  &  1  &  0  &  0  &  0  & -1   \\	
			-1  &  0  &  0  & -1  &  0  &  0   \\	
			 0  & -1  &  0  &  1  &  0  &  0   \\	
			 1  &  0  & -1  &  0  &  0  &  0   \\	
			 0  &  1  &  1  &  0  &  0  &  0  
		\end{bmatrix*}
		, \quad
		\Mat{E}{2,3} = 
		\begin{bmatrix*}[r]
			 -1  \\  1  \\  -1  \\  1  \\  -1  \\  1
		\end{bmatrix*}
	\end{align*}
\end{small}
One can readily check that the boundary of the boundary is empty,
\begin{align*}
	\Mat{E}{0,1} \Mat{E}{1,2} = \begin{pmatrix}0 & ... & 0 \end{pmatrix}^T \quad \text{and} \quad \Mat{E}{1,2} \Mat{E}{2,3} = \begin{pmatrix}0 & ... & 0 \end{pmatrix}^T
\end{align*}
\label{ex:IncidenceMatrix}
\end{example}

\begin{figure}
	\centering
	\includegraphics[width=0.75\textwidth]{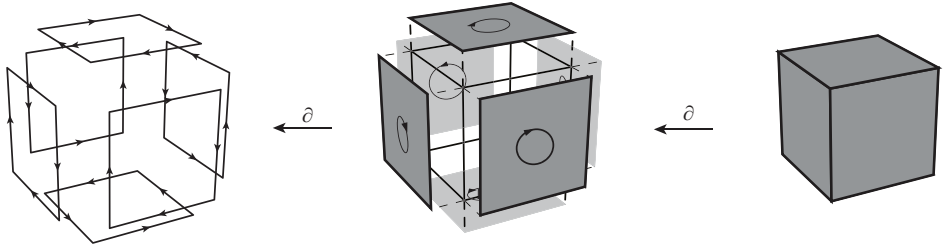}%
	\caption{A double application of the boundary operator leads to an empty chain. Note that all edges are oriented in the opposite way, so any value associated with them cancels.}
	\label{fig:BoundaryEmpty}
\end{figure}

We now introduce the most important property of the boundary operator, which leads to some profound consequences in the description of geometry and physics. Picture for example a finite dimensional volume $V$ in $\mathbb{R}^3$. It is clear that its bounding surface $\boundary V$ encloses all of $V$ (this is in fact the definition of the boundary). This means that the surface $\boundary V$ has no boundary itself; it is boundaryless. Figure \ref{fig:BoundaryEmpty} illustrates that taking the boundary twice of a volume leads to an empty 1-chain. 

In general, taking the boundary twice of a $k$-chain leads to an empty $(k-2)$-chain,
\begin{align}
	\boundary \boundary \chain{c}{k} = \chain{0}{k-2} \quad \text{for all} \;  \chain{c}{k} \in \chain{C}{k}(D).
	\label{eq:BE}
\end{align}
The set of $k$-chains and boundary operators thus gives rise to an exact sequence on contractible domains, the chain complex $\left( \chain{C}{k}(D), \boundary  \right)$
\begin{align}
\begin{diagram}
	\cdots &   \lTo^{\boundary}  & \chain{C}{k-1}(D) & \lTo^{\boundary} & \chain{C}{k}(D) &  \lTo^{\boundary} & \chain{C}{k+1}(D) & \lTo^{\boundary} & \cdots 
\end{diagram}
\label{eq:emptyboundary}
\end{align}

\subsection{Cochains and the coboundary operator}
So far we have learned how to calculate with discrete geometric objects, such as oriented points, edges, faces and volumes. By duality pairing with discrete geometry, we introduce cochains as discrete analogues of differential forms. Cochains can describe for instance  integral values along a chain. Integration itself is however nothing but the 'measurement technique' that assigns global values to cells. For our purpose, it is most convenient to think of cochains in the most general way as any set of global values associated with oriented points, edges, faces and volumes in the mesh (cell complex). More clearly, they need not be integral values. This fact will allow us maximum freedom in the projection of differential forms - to be discussed in the next section - where we use cochains as degrees of freedom in combination with suitable basis functions based on NURBS.

A $k$-cochain is an expression which looks like,
\begin{align}
	\cochain{a}{k} = \left\{a_i \cdot \sigma^{(k),i} \right\}_{i=0}^s \quad \text{with} \; a_i \in \mathbb{R} \quad \text{and} \quad \sigma^{(k),i} \in \cochain{C}{k}.
\end{align}
In contrast to chains, we use subscript for the coefficients $a_i$ and superscript to denote the basis  $\sigma^{(k),i}$. While we considered chains only with coefficients $\left\{-1,0,1 \right\}$, we allow the coefficients of cochains to be any real number. As with chains, we shall use the coefficients directly as the column vector $\vect{a}{} = \begin{pmatrix} a_0 & a_1 & \cdots & a_s  \end{pmatrix}^T$ in matrix calculations.

Cochains are linear functionals on chains, and by choosing the basis cochains $\sigma^{(k),i}$ dual to that of the space of chains, $\left\langle \sigma^{(k),i},\sigma_{(k),j}\right\rangle = \delta^i_j$, we can define the duality pairing between chains and cochains as,
\begin{align}
	\left\langle \cochain{a}{k}, \chain{c}{k} \right\rangle := \sum_{i=0}^s  a_i \cdot c^i  = \vect{a}{}^T \vect{c}{} \in \mathbb{R}
	\label{eq:DualityPairing2}
\end{align}
Note the similarity between duality pairing of differential forms and geometry by means of integration (\ref{eq:DualityPairing}). In the discrete setting, integration is replaced by summation (\ref{eq:DualityPairing2}).

We can now define a discrete Stokes theorem, the mother of all equations, in terms of chains and cochains. As in the continuous setting where the exterior derivative is the formal adjoint of the boundary operator (\ref{eq:Stokes}), we can define the coboundary operator on cochains as the formal adjoint of the boundary operator on chains,
\begin{align}
	\left\langle 	\coboundary \cochain{b}{k-1} , \chain{c}{k} \right\rangle  :=  \left\langle \cochain{b}{k-1} ,  \boundary \chain{c}{k} \right\rangle
	\label{eq:DiscreteStokes}
\end{align}
for all $ \cochain{b}{k-1} \in \chain{C}{k-1}(D)$ and $\chain{c}{k} \in \chain{C}{k}(D)$. While the boundary is a map, $\boundary \; : \; \chain{C}{k} \mapsto \chain{C}{k-1}$, the coboundary is a map,  $\coboundary \; : \; \cochain{C}{k-1} \mapsto \cochain{C}{k}$. Analogous to the exterior derivative $\diff$ acting on forms in the continuous setting, the coboundary operator acts as a discrete derivative on cochains. The coboundary is thus a discrete version of the gradient $(k=1)$, curl $(k=2)$ and divergence $(k=3)$ operators from vector calculus.

\begin{figure}
	\centering
	\includegraphics[trim = 0cm 2cm 0cm 0cm, clip,width=0.75\textwidth]{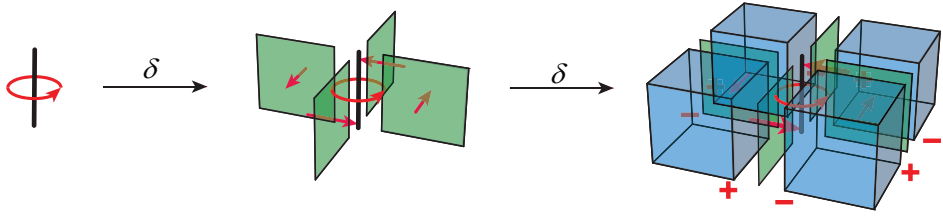}%
	\caption{A double application of the coboundary operator leads to an empty cochain. Note that the two additions of face values per volume cancel each other out.}
	\label{fig:CoboundaryEmpty}
\end{figure}

Earlier we remarked the fact that $\boundary \circ \boundary = 0$ , has important consequences in physics. Amongst others, it implies that $\coboundary \circ \coboundary = 0$, see Figure \ref{fig:CoboundaryEmpty}. This is easily proven using (\ref{eq:DiscreteStokes}) twice and (\ref{eq:BE}), 
\begin{align*}
	  \left\langle \coboundary \coboundary \cochain{b}{k-1}, \chain{e}{k+1} \right\rangle 
	 \stackrel{(\ref{eq:DiscreteStokes})}{=} \left\langle \coboundary \cochain{b}{k-1} , \boundary \chain{e}{k+1}  \right\rangle 
	 \stackrel{(\ref{eq:DiscreteStokes})}{=} \left\langle \cochain{b}{k-1} , \boundary \boundary \chain{e}{k+1} \right\rangle 
	 \stackrel{(\ref{eq:BE})}{=} 0
\end{align*}
This property is analogous to (\ref{eq:diff4}) and the important vector calculus identities $\curl \; \grad = 0$ and $\div \; \curl = 0$ in continuous space.

We can now set up the following exact sequence, known as a cochain complex $\left( \cochain{C}{k}(D), \coboundary  \right)$,
\begin{align}
\begin{diagram}
	\cdots &   \rTo^{\coboundary}  & \cochain{C}{k-1}(D) & \rTo^{\coboundary} & \cochain{C}{k}(D) &  \rTo^{\coboundary} & \cochain{C}{k+1}(D) & \rTo^{\coboundary} & \cdots 
\end{diagram}
\end{align}

Since we have a matrix representation of the boundary operator, we can derive a matrix representation of the coboundary operator, using the adjoint property of the boundary and coboundary operator (\ref{eq:DiscreteStokes}),
\begin{align}
	\left\langle \cochain{b}{k-1} , \boundary \chain{c}{k} \right\rangle \stackrel{(\ref{eq:DualityPairing2})}{=} \vect{b}{}^T \left( \Mat{E}{k-1,k} \vect{c}{} \right) = \left( \Mat{E}{k-1,k}^T \vect{b}{} \right)^T \vect{c}{} \stackrel{(\ref{eq:DualityPairing2})}{=} \left\langle \coboundary \cochain{b}{k-1},\chain{c}{k} \right\rangle.
\end{align}
We can conclude that the matrix representation of the coboundary operator is given by $\Mat{D}{k,k-1} = \Mat{E}{k-1,k}^T$. We can therefore perform differentiation in the discrete setting using the matrix representation $\vect{a}{} = \Mat{D}{k,k-1} \vect{b}{}$. Furthermore, $\Mat{D}{k+1,k} \Mat{D}{k,k-1} = \begin{pmatrix} 0 & \cdots & 0 \end{pmatrix}^T$, so the matrix representation of the gradient, curl and divergence exactly preserves the null space of the differential operators.
\section{Commuting projection \label{sec:CommutingProjection}}
Balance equations, represented by the Generalized Stokes theorem (\ref{eq:Stokes}), allow an exact discrete representation in terms of chains and cochains using (\ref{eq:DiscreteStokes}). Since the structure of these equations is inherently discrete and metric free, we do not expect that basis functions play any role here. The constitutive equations - the material and metric dependent relations -  on the other hand, require a continuous formulation of both the geometry and the field variables, and here is where the basis functions come in to play. 

A continuous representation of field variables involves approximation by a suitable projection, $	\pi_h \; : \; \Lambda^{k} (\Omega) \mapsto \Lambda^{k}_h (\Omega)$,
from an infinite dimensional space $ \Lambda^{k} (\Omega) $ of differential forms, to a finite dimensional conforming subspace  $ \Lambda^{k}_h (\Omega) \subset \Lambda^{k} (\Omega)$. Approximation always involves a loss in information. It is however important that the error is bounded by a certain constant $C < \infty$, 
\begin{align}
	& \text{approximation property:} \quad \left\| \Dform{a}{k} - \pi_h \Dform{a}{k}  \right\|_{\Lambda^k} \leq C \cdot h^p \quad \text{for all } \Dform{a}{k} \in \Lambda^{k}(\Omega)
	\label{eq:approximation}
\end{align}
for a norm defined on $\Lambda^k(\Omega)$, where $p$ denotes the approximation order and $h$ is a measure of the the maximum partitioning size.  In order to be consistent, the projection should reproduce all of $ \Lambda^{k}_h(\Omega)$,  which means that,
\begin{align}
	& \text{consistency property:} \quad \pi_h \circ \pi_h \; \Dform{a}{k} = \pi_h \Dform{a}{k} \quad \text{for all } \Dform{a}{k} \in \Lambda^{k}(\Omega)
	\label{eq:consistency}
\end{align}
These two requirements are in general met by polynomial approximations. For interpolation estimates of for example NURBS and splines, see \cite{Bazilevs:2006a,Bramble:1970} respectively. 

To capture the behavior of a physical phenomena well, a discretization method should not only approximate the spaces of the infinite dimensional system, but should also follow the structure induced by the relations between them, i.e. the structure induced by the balance laws and constitutive equations (Figure \ref{fig:structure}). Only then will a discrete representation of the physics have any physical significance. The focus of this paper is to derive a projection of differential forms which commutes with differentiation,
\begin{align}
\label{CD:1}
\centering
\xymatrix{
	\Lambda^{k} \left(\Omega \right) \ar[d]_{\pi_h} \ar[r]^{\diff} & 	\Lambda^{k+1} \left(\Omega \right) \ar[d]^{\pi_h}  \\
	\Lambda^{k}_h \left(\Omega\right)  \ar[r]^{\diff} &	\Lambda^{k+1}_h \left(\Omega \right)
}
\end{align}
Projection and subsequent differentiation should give the same result as first applying the derivative and then do the projection. This commuting diagram property will guarantee that conservation and balance laws remain exactly satisfied in the discrete setting. Furthermore, this property is the key to naturally stable and consistent numerical algorithms. It can be shown that (\ref{CD:1}) implies a conforming Hodge decomposition, which together with the Poincar\'e inequality proves inf-sup stability. For more details we refer the reader to the papers by Kreeft et. al. \cite{Kreeft:2012b,Kreeft:2011}.

To assure that the balance laws are exactly satisfied, even on curved meshes, we require that the projection is independent of geometric transformations. This means that the projection should be constructed such that it commutes with the pull back as well,
\begin{align}
\label{CD:PullBack}
\centering
\xymatrix{
	\Lambda^{k} \left(\Phi(\Omega') \right) \ar[d]_{\pi_h} \ar[r]^{\Phi^{\star}} & 	\Lambda^{k} \left(\Omega' \right) \ar[d]^{\pi_h}  \\
	\Lambda^{k}_h \left(\Phi(\Omega')\right)  \ar[r]^{\Phi^{\star}} &	\Lambda^{k}_h \left(\Omega' \right)
}
\end{align}
This can be achieved by constructing the projection such that it preserves the integrals of a differential $k$-form, $\Dform{a}{k}$, along a chain $\chain{c}{k}$. Using the fact that the pull back $\Phi^{\star}$ is the adjoint of the map $\Phi$ (\ref{eq:PullBack}), we obtain the desired result,
\begin{align}
	\int_{\chain{c}{k}} \pi_h  \Phi^{\star} \Dform{a}{k}  =
	\int_{\chain{c}{k}} \Phi^{\star} \Dform{a}{k} = 
	\int_{\Phi(\chain{c}{k})} \Dform{a}{k} = 
	\int_{\Phi(\chain{c}{k})} \pi_h \Dform{a}{k} = 
	\int_{\chain{c}{k}} \Phi^{\star} \pi_h \Dform{a}{k}.
\end{align}

\subsection{Reduction, change of bases and reconstruction}
We extend the commuting diagram in (\ref{CD:1}) to define operations between the continuous (Section \ref{sec:DifferentialModeling}) and discrete formalism (Section \ref{sec:DiscreteModeling}). Following \citet{Bochev:2006}, we introduce two separate operators, the reduction $\Red{}$ and  reconstruction $\Int{}$ - which map differential forms to cochains and cochains to a finite dimensional representation of differential forms, respectively. Additionally we introduce a change of basis from cochains to a different type of cochains. As we shall see, this change of basis will allow maximum freedom in the choice of reconstruction method.

\vspace{0.25cm}
\textbf{Reduction} $\Red{} \; : \; \Lambda^{k}(\Omega) \mapsto \cochain{C}{k}(D) $ is an abstraction of the measurement process and is given by the DeRahm map (\ref{eq:DeRahm}), which is defined by means of integration,
\begin{align}
	\left\langle \Red{} \left( \Dform{a}{k} \right), \chain{c}{k} \right\rangle = \int_{\chain{c}{k}} \Dform{a}{k}.
	\label{eq:DeRahm}
\end{align}

To be able to use the coboundary as a discrete derivative, as outlined in Section \ref{sec:DiscreteModeling}, we require that the reduction commutes with differentiation,
\begin{align}
\label{CD:2}
\centering
\xymatrix{
	\Lambda^{k} \left(\Omega \right) \ar[d]_{\Red{}} \ar[r]^{\diff} & 	\Lambda^{k+1} \left(\Omega \right) \ar[d]^{\Red{}}  \\
	\cochain{C}{k} \left(D\right)  \ar[r]^{\coboundary} &	\cochain{C}{k+1} \left(D\right)
}
\end{align}
which is indeed the case as can be proven by using Stokes theorem and the duality between the boundary and coboundary operator,
\begin{align*}
	\left\langle \Red{} \left( \diff \Dform{a}{k} \right) ,\chain{c}{k+1} \right\rangle  \stackrel{(\ref{eq:DeRahm})}{=}  \int_{\chain{c}{k+1}} \diff \Dform{a}{k} \stackrel{(\ref{eq:Stokes})}{=} \int_{\boundary \chain{c}{k+1}} \Dform{a}{k} \stackrel{(\ref{eq:DeRahm})}{=}   
	\left\langle \Red{} \left( \Dform{a}{k} \right), \boundary \chain{c}{k+1} \right\rangle \stackrel{(\ref{eq:DiscreteStokes})}{=}  \left\langle  \coboundary \Red{} \left( \Dform{a}{k} \right) , \chain{c}{k+1} \right\rangle 
\end{align*}

\vspace{0.25cm}

\textbf{Change of basis} $ \mathcal{F} \; : \; \cochain{\bar{C}}{k}(D) \mapsto \cochain{C}{k}(D)$ is an invertible map from a cochain $\cochain{\bar{a}}{k} = \left\{\bar{a}_j \cdot \bar{\sigma}^{(k),j}  \right\}_{i=0}^s \in \cochain{\bar{C}}{k}(D)$, to a different type of cochain $\cochain{a}{k} = \left\{a_j \cdot \sigma^{(k),j}  \right\}_{i=0}^s \in \cochain{C}{k}(D)$ which are defined by means of integration (\ref{eq:DeRahm}). The change of basis $\mathcal{F}$ makes it possible to use a much wider range of reconstruction methods. It provides the freedom to use degrees of freedom other than nodal values, curve, surface and volume integrals, while maintaining all the advantages of a compatible approach. For the reconstruction of 0-forms, for example, we no longer require nodal interpolants, but can use any basis that possesses partition of unity, for example NURBS. In practice, the change of basis is given by a square invertible matrix equivalent to that in an interpolation / histopolation problem. 

Since $\mathcal{F}$ is a map from cochains to cochains, we require that it commutes with discrete differentiation by means of the coboundary operator,
\begin{align}
\label{CD:3}
\centering
\xymatrix{
	\cochain{C}{k} \left(D\right)  \ar[d]_{\mathcal{F}^{-1}} \ar[r]^{\coboundary} & 	\cochain{C}{k+1} \left(D\right) \ar[d]^{\mathcal{F}^{-1}}  \\
	\cochain{\bar{C}}{k} \left(D\right)  \ar[r]^{\coboundary} &	\cochain{\bar{C}}{k+1} \left(D\right)
}
\end{align}
The commuting diagram in (\ref{CD:3}) guarantees that the mathematical structure remains exactly the same after change of basis. The discrete modeling tools introduced in Section \ref{sec:DiscreteModeling} therefore remain unaltered and we can perform discrete differentiation using the matrix representation $\Mat{D}{k+1,k} = \Mat{E}{k,k+1}^T$, as explained in the end of Section \ref{sec:DiscreteModeling}. 

\vspace{0.25cm}
\textbf{Reconstruction} $	\Int{} \; : \; \cochain{\bar{C}}{k}(D) \mapsto \Lambda_{h}^{k}(\Omega) $ maps cochains $\cochain{\bar{a}}{k} = \left\{\bar{a}_j \cdot \bar{\sigma}^{(k),j}  \right\}_{i=0}^s \in \cochain{\bar{C}}{k}(D)$ back to a finite dimensional representation of a differential form $\Dform{a}{k}_h \in \Lambda_{h}^{k}(\Omega) \subset \Lambda^{k} (\Omega)$. In order to use the discrete calculus introduced in Section \ref{sec:DiscreteModeling} the reconstruction is required to commute with differentiation,
\begin{align}
\label{CD:4}
\centering
\xymatrix{
	\cochain{\bar{C}}{k} \left(D\right)  \ar[d]_{\Int{}} \ar[r]^{\coboundary} & 	\cochain{\bar{C}}{k+1} \left(D\right) \ar[d]^{\Int{}}  \\
	\Lambda^{k}_h \left(\Omega\right)  \ar[r]^{\diff} &	\Lambda^{k+1}_h \left(\Omega\right)
}
\end{align}
This way we can perform the gradient, curl and divergence in the discrete setting, while continuous representations can be reconstructed where and whenever required.

\vspace{0.5cm}
The commuting property of projection with differentiation (\ref{CD:1}) can be seen as the composition of the commuting relations in (\ref{CD:2}), (\ref{CD:3}) and (\ref{CD:4}) and is illustrated in the following diagram,
\begin{align}
\label{CD:5}
\centering
\xymatrix{
	\Lambda^{k} \left(\Omega \right) \ar@<-2.5ex>@/_1pc/[ddd]_{\pi_h} \ar[d]_{\Red{}} 
	\ar[r]^{\diff} & 	\Lambda^{k+1} \left(\Omega \right) \ar[d]^{\Red{}}  \ar@<2.5ex>@/^1pc/[ddd]^{\pi_h}          \\
	\cochain{C}{k}(D) \ar@<-0.5ex>[d]_{\mathcal{F}^{-1}}  \ar[r]^{\coboundary} & \cochain{C}{k+1}(D)  \ar@<-0.5ex>[d]^{\mathcal{F}^{-1}}      \\
	\cochain{\bar{C}}{k}(D) \ar@<-0.5ex>[d]_{\Int{}}  \ar[r]^{\coboundary} & \cochain{\bar{C}}{k+1}(D)  \ar@<-0.5ex>[d]^{\Int{}}    \\
	\Lambda^{k}_h \left(\Omega\right)  \ar[r]^{\diff} &	\Lambda^{k+1}_h \left(\Omega \right)
}
\end{align}

The choice of reconstruction, $\Int{} \; : \; \cochain{\bar{C}}{k}(D) \mapsto \Lambda^{k}_h \left(\Omega\right)$, uniquely defines that of $\Int{} \; : \; \cochain{\bar{C}}{k+1}(D) \mapsto \Lambda^{k+1}_h \left(\Omega\right)$ through commutation with differentiation (\ref{CD:4}). This in practice means that we choose a basis for $\Lambda^{0}_h \left(\Omega\right)$ and derive the remaining finite dimensional spaces of differential forms, $\Lambda^{0}_h \left(\Omega\right)$, for $k=1,...,n$, using the commuting properties.

Furthermore, from (\ref{eq:consistency}) it naturally follows that the change of basis is related to the reduction and reconstruction by, $\mathcal{F} \; := \; \Red{} \circ \Int{}$. Therefore the choice of basis for 0-forms determines not only the basis for all other spaces of differential forms, but also uniquely defines the change of basis $ \mathcal{F} \; : \; \cochain{\bar{C}}{k}(D) \; \mapsto  \cochain{C}{k}(D)$ for $k=0,1,...,n$.

\begin{figure}
	\centering
	\includegraphics[width=0.75\textwidth]{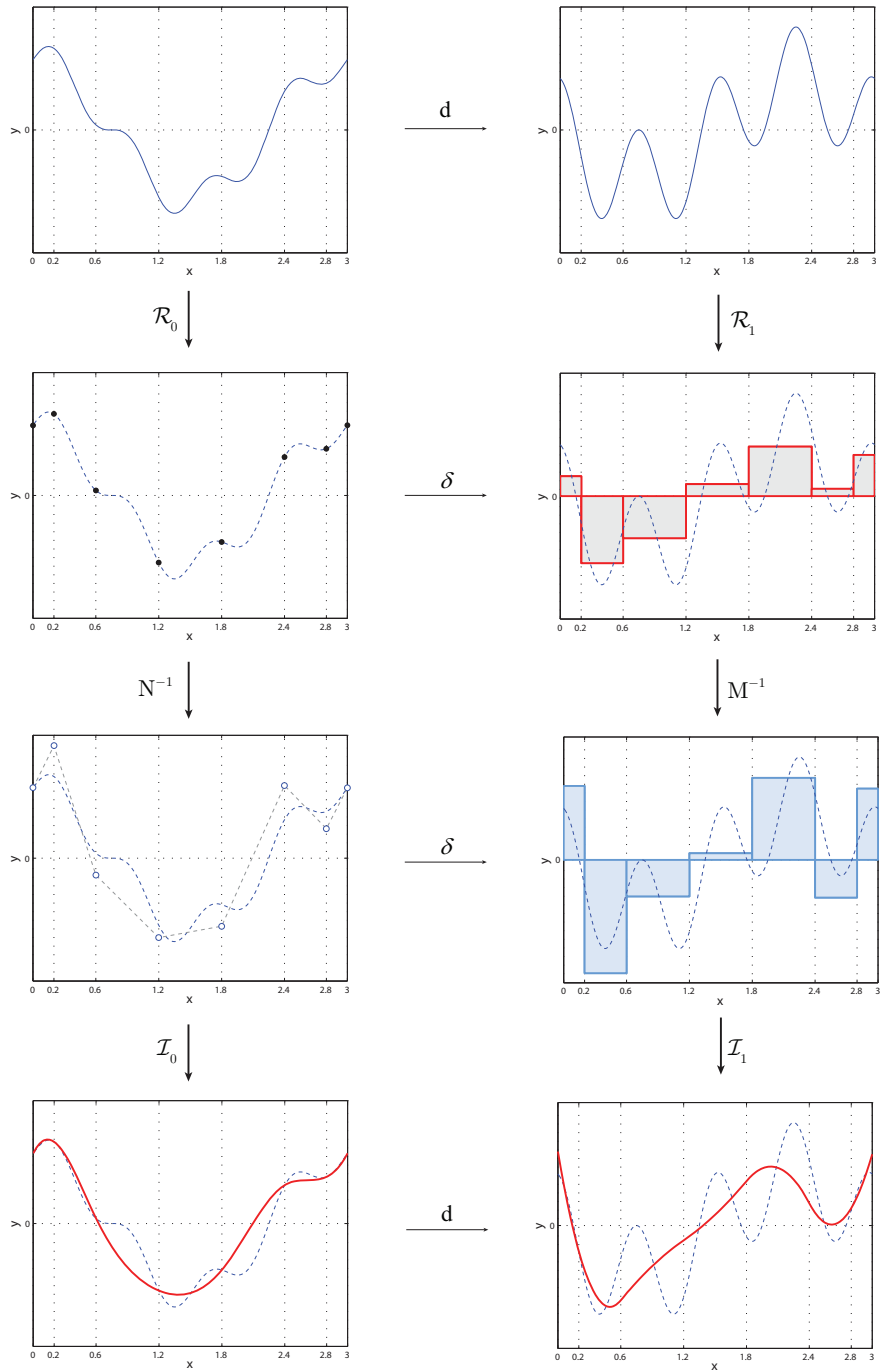}
	\caption{While a zero form is projected by interpolating a set of nodal values (left), its derivative (right) is projected by interpolating a set of integral values; a process called histopolation. This construction guarantees that the projection commutes with differentiation. Furthermore, the projection is broken down in three stages: 1) the reduction process $\mathcal{R}$ reducing a 0-form and its derivative to a set of nodal and integral values, respectively; 2) A change of basis from nodal and integral values to new 'node' and 'edge' type of degrees of freedom. This change of basis is equivalent to solving an interpolation problem $\mathrm{N}^{-1}$ and histopolation problem $\mathrm{M}^{-1}$, respectively. 3) the reconstruction process $\mathcal{I}$ which is simply given by a linear combination of the new 'node' (0-cochain), or 'edge' (1-cochain) degrees of freedom with the appropriate basis functions $\left\{ N_i(x) \right\}_{i=0}^n$ and $\left\{ M_i(x) \right\}_{i=1}^n$ respectively.	The resolution is purposely kept low to illustrate the concepts of reduction, change of basis and reconstruction.}
	\label{fig:commutation}
\end{figure}

We will study the commuting properties in (\ref{CD:5}) in more detail in the univariate case, see Figure \ref{fig:commutation}. Once all spaces and operators are determined in the univariate case, we shall use the tensor product to develop multivariate spaces of discrete differential forms.

\subsection{Interpolation and histopolation}
The commuting relations in (\ref{CD:5}) are graphically illustrated in Figure \ref{fig:commutation} in the case of 1D-space. Projection of a 0-form is equivalent to solving an interpolation problem where a set of data measurements are interpolated, while the projection of a 1-form is equivalent to solving a histopolation problem, where a set of line integrals are preserved in the projection.

\vspace{0.25cm}

Say we wish to approximate the 0-form in Figure \ref{fig:commutation}, which represents for example a temperature field $T(x)$. We seek the finite dimensional approximation $T_h(x) = \pi_h T(x)$  in the following space,
\begin{align}
	 & \Lambda^0_h (\Omega) := \left\{ \Int{} \left\{\bar{T}_j \cdot \bar{\sigma}^{(0),j} \right\}_{j=0}^n = \sum_{j=0}^n \bar{T}_j \; N_j(x), \; \text{for all} \;	\bar{T}_j \in \mathbb{R} \right\},
	 \label{eq:projection0}
\end{align}
where the degrees of freedom $\bar{T}_j$ are the coefficients in the 0-cochain $\left\{\bar{T}_j \cdot \bar{\sigma}^{(0),j} \right\}_{j=0}^n$ and the $\left\{N_j(x)\right\}_{j=0}^{n}$ are linear independent basis functions that posses a partition of unity, $\sum_{j=0}^n N_{j}(x) = 1$. 

The projection of a 0-form $T_h(x) = \pi_h T(x)$ is required to follow $\Red{} \left( T(x) \right) = \Red{} \left( T_h(x) \right)$, which means that $T_h(x)$ interpolates $T(x)$ at a chosen set of nodes $\left\{x_i \right\}_{i=0}^{n}$. Then, using (\ref{eq:projection0}), we can set up the $n+1$ by $n+1$ system of linear equations,
\begin{align}
	\pi_h T(x_i) = \sum_{j=0}^n \bar{T}_j \; N_{j}(x_i) = T(x_i) \quad \text{for} \; i=0,1,...,n
	\label{eq:interpolation}
\end{align}
and find the unique set of coefficients $\left\{\bar{T}_j \right\}_{i=0}^n$. This system is guaranteed to have a solution when the matrix $\mathrm{N}$ with coefficients $\mathrm{N}_{ij} = N_{j}(x_i)$ is invertible, which is the case when $x_i \in \Span{N_i(x)}$, for $i=0,...,n$ \cite{DeBoor:1978}. Note that the change of basis $\mathcal{F} = \Red{} \circ \Int{}$ is thus given by this square interpolation matrix $\mathrm{N}$. The inverse map  $\mathcal{F}^{-1} = \mathrm{N}^{-1}$, depicted in the middle of the left hand side of Figure \ref{fig:commutation}, is thus equivalent to solving the interpolation problem.

\vspace{0.25cm}

If we would like to approximate the derivative of the temperature field, the 1-form $\Dform{u}{1}(x) = \diff T(x)$, we require a basis for 1-forms,
\begin{align}
	 & \Lambda^1_h (\Omega) = \left\{ \Int{} \left\{ \bar{u}_j \cdot \bar{\sigma}^{(1),j} \right\}_{i=1}^n = \sum_{i=1}^n \bar{u}_j \; M_j(x), \; \text{for all} \; 			\bar{u}_j \in \mathbb{R} \right\}
	 \label{eq:projection1}
\end{align}
where the degrees of freedom $\bar{u}_j$ are the coefficients in the 1-cochain $\left\{ \bar{u}_j \cdot \bar{\sigma}^{(1),j} \right\}_{i=1}^n$ and therefore the basis functions $\left\{M_j(x)\right\}_{j=1}^n$ can be associated with edges.

The projection of a 1-form $ \Dform{u}{1}_h = \pi_h \Dform{u}{1}$ is required to follow $\Red{} \left(\Dform{u}{1} \right) = \Red{} \left(\Dform{u}{1}_h \right)$, which means that the approximation should preserve the line integrals of the 1-form $u^1(x) = \diff \; T(x)$ along the edges $e_i = [x_{i-1},x_i]$. We can therefore set up the following $n$ by $n$ system of linear equations,
\begin{align}
	\int_{e_i} \pi_h \Dform{u}{1} = \sum_{j=1}^n \bar{u}_j \; \int_{e_i} M_{j}(x) = \int_{e_i} \Dform{u}{1}(x) = T(x_i)-T(x_{i-1}) \quad \text{for} \; i=1,2,...,n
	\label{eq:histopolation}
\end{align}
This process (right side of Figure \ref{fig:commutation}) is called histopolation, from interpolation of a histogram. In this case, the change of basis $\mathcal{F} = \Red{} \circ \Int{}$, is equal to the invertible histopolation matrix $\mathrm{M}$, of  which the coefficients are given by,
\begin{align*}
	\mathrm{M}_{ij} = \int_{e_i} M_j(x)
\end{align*}
The inverse map  $\mathcal{F}^{-1} = \mathrm{M}^{-1}$, depicted in the middle of the right hand side of Figure \ref{fig:commutation}, is thus equivalent to solving the histopolation problem.

\subsection{Edge functions}
The idea now is to choose a basis $\left\{N_i(x)\right\}_{i=0}^n$ for 0-forms, and to derive the new edge type of basis functions $\left\{M_i(x)\right\}_{i=1}^n$ by using the commuting properties in (\ref{CD:4}). Connectivity of edges $e_i$ and points $x_i$, $\partial e_i = x_{i} - x_{i-1}$, implies that we can apply discrete differentiation using the coboundary (\ref{eq:DiscreteStokes}), $\bar{u}_i = \bar{T}_i - \bar{T}_{i-1}$ , for $i=1,2,...,n$. Hence, every coefficient $\bar{T}_i$ can be written as $\bar{T}_i = \bar{T}_0 + \sum_{j=1}^i \left(\bar{T}_j - \bar{T}_{j-1}) \right) = \bar{T}_0 + \sum_{j=1}^i \bar{u}_j$.

Using this result in the projection of a zero form and applying the exterior derivative,
\begin{align*}
	\diff \pi_h \; T(x) = \diff \sum_{i=0}^n \bar{T}_i \;  N_{i}(x)
									&= \bar{T}_0 \; \diff \sum_{i=0}^n N_i(x) + \sum_{i=0}^n  \left( \sum_{j=1}^i \bar{u}_j \right) \diff N_{i}(x).
\end{align*}
Partition of unity implies that, $\diff \sum_{j=0}^n N_j(x) = 0$; hence the first term cancels. The second term can be rewritten such that we obtain a relation for the edge basis functions $\left\{M_i(x)\right\}_{i=1}^n$,
\begin{align*}
		\diff \pi_h \; T^0(x)
									&= \bar{u}_1\; \diff N_1(x) + \left\{\bar{u}_1+\bar{u}_2\right\} \diff N_2(x) + ... \\ 
								  &= \left\{\diff N_1(x) +...+\diff N_n(x) \right\} \bar{u}_1 + 
								     \left\{\diff N_2(x)+...+\diff N_n(x) \right\} \bar{u}_2 + ...
								   = \sum_{i=1}^n \bar{u}_i \sum_{j=i}^n \diff N_j(x).
\end{align*}
By defining the edge basis functions $\left\{M_i(x) \right\}_{i=1}^n$ as,
\begin{align}
		M_i(x)	=  \sum_{j=i}^n \diff N_j(x) = - \sum_{j=0}^{i-1} \diff N_j(x), \quad \text{for} \; i=1,2,...,n,
\end{align}
we have obtained the commuting projection,
\begin{align}
	\diff \pi_h \; T^0(x) = \diff \sum_{i=0}^n \bar{T}_i \;  N_{i}(x)
											  = \sum_{i=1}^n \left(\bar{T}_i - T_{i-1} \right) \; M_{i}(x)
											  = \sum_{i=1}^n \bar{u}_i \; M_{i}(x) = \pi_h \diff \; T^0(x).
	\label{eq:CommutingProjection}
\end{align}
As constructed, we can perform differentiation entirely independent from the basis functions. In fact, in $\diff \Dform{T}{0}_h = \sum_{i=1}^n \left(\bar{T}_i - T_{i-1} \right) \; M_{i}(x) = \sum_{i=1}^n \bar{u}_i \; M_{i}(x) = \Dform{u}{1}$, the basis functions $M_i(x)$ can be canceled out and we are left with discrete differentiation, $\bar{u}_i = \bar{T}_i - \bar{T}_{i-1}$ using the coboundary process as described in Section \ref{sec:DiscreteModeling}. The conservation laws thus become completely independent of the basis functions. This means that they are exactly satisfied, even on coarse and on arbitrary curved meshes, as already proven by the commutation property of the projection with the pullback (\ref{CD:PullBack}).

While the set $\left\{N_i\right\}_{i=0}^n$ sums up to 1 (partition of unity), the set of edge basis functions, $\left\{M_i\right\}_{i=1}^n$, scale such that they provide unit integral,
\begin{align}
	\bar{M}_i := \int_{x_0}^{x_n} M_i (x) = 1. \quad \text{for} \; i=1,2,...,n.
\end{align}
We proof this in the case that the set $\left\{N_i(x)\right\}_{i=0}^n$ has interpolating end-conditions, i.e.  $ T^0_{h}(x_0) = \sum_{i=0}^{n} \bar{T}_i \; N_i (x_0) = \bar{T}_0$ and $T_h(x_n) = \sum_{i=0}^{n} \bar{T}_i \; N_i (x_n) = \bar{T}_n$. Then we have by Stokes theorem (\ref{eq:Stokes}),
	\begin{align*}
		\int_{x_0}^{x_n} \diff T^0_h(x)  =  \int_{x_0}^{x_n} \sum_{i=1}^n \left(\bar{T}_i - \bar{T}_{i-1} \right) M_i(x)
																	   =  \bar{T}_n \; \bar{M}_n +  \sum_{i=1}^{n-1} \bar{T}_i \left(\bar{M}_i-\bar{M}_{i+1}\right) -  \bar{T}_0 \; \bar{M}_1 = \bar{T}_n - \bar{T}_0.
	\end{align*}
Since this should hold for all $\left\{\bar{T}_i \in \mathbb{R} \right\}_{i=0}^n$, $\bar{M}_1 = \bar{M}_n = 1$ and consequently $\left\{\bar{M}_i = 1\right\}_{i=1}^n$.
%\begin{small}
%\begin{align*}
%	\begin{pmatrix} \bar{T}_0 & \cdots & \bar{T}_i & \cdots & \bar{T}_n  \end{pmatrix}
%	\begin{pmatrix} -1     	&   				& 				&		\oslash	\\
%									\ddots 	&  \ddots  	&      		&						\\
%									       	&    1     	& -1  		&						\\
%									       	&          	& \ddots	&	\ddots 		\\
%								  \oslash &  					&					&		1			
%	\end{pmatrix}
%	\begin{pmatrix} \bar{M}_1 \\ \vdots \\ \bar{M}_i \\ \bar{M}_{i+1} \\ \vdots \\ \bar{M}_n  \end{pmatrix} =
%	\begin{pmatrix} \bar{T}_0 & \cdots &\bar{T}_i & \cdots & \bar{T}_n  \end{pmatrix}
%	\begin{pmatrix} -1 \\ \vdots \\ 0 \\ \vdots \\ 1  \end{pmatrix}
%\end{align*}
%%\end{small}
%for which the result should hold for all $\bar{T}_i \in \mathbb{R}$, $i=0,1,...,n$. Although we are left with an overdetermined system of $n+1$ equations and $n$ unknowns, we can use back substitution to find that $\left\{\bar{M}_i = 1\right\}_{i=1}^n$ is the unique solution.

\begin{figure}
	\centering
	\subfigure[]{\includegraphics[trim = 0cm 1cm 0cm 1cm,clip,width=0.49\textwidth]{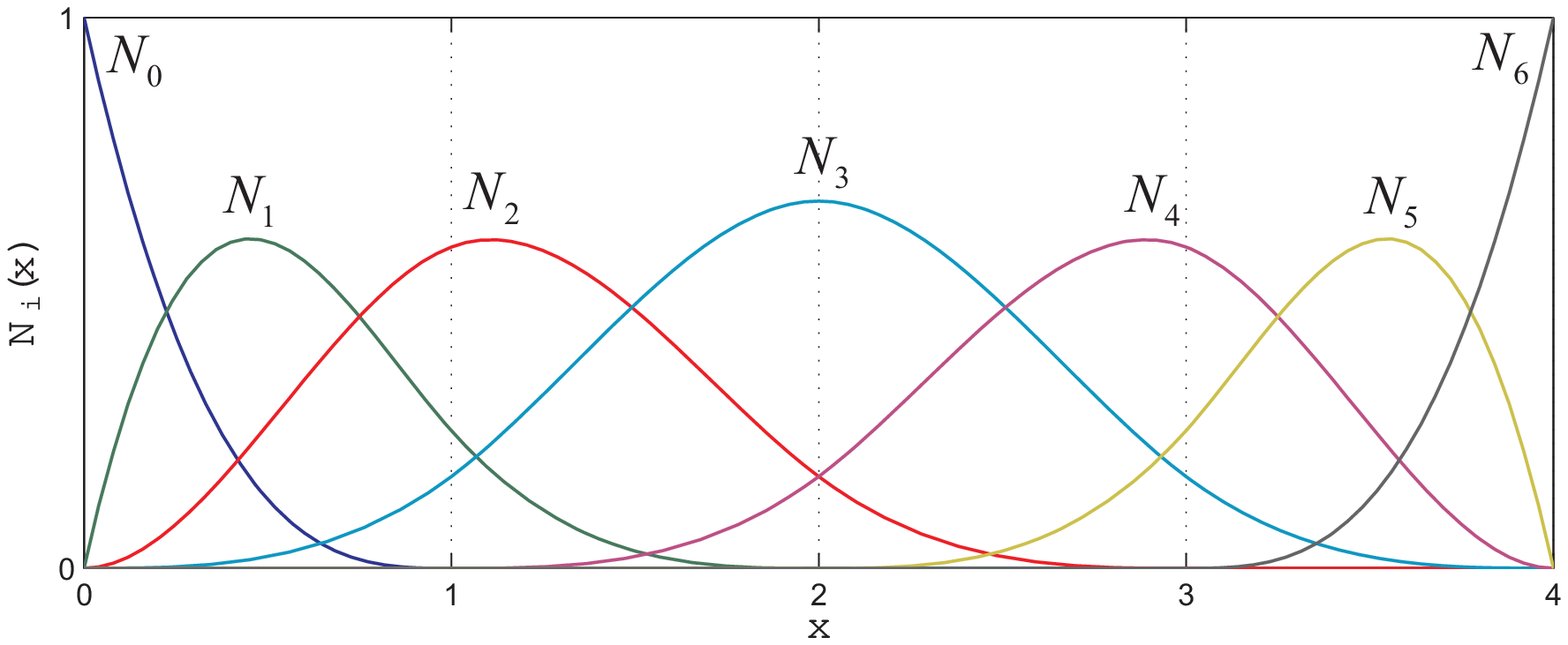}}
	\subfigure[]{\includegraphics[trim = 0cm 1cm 0cm 1cm,clip,width=0.49\textwidth]{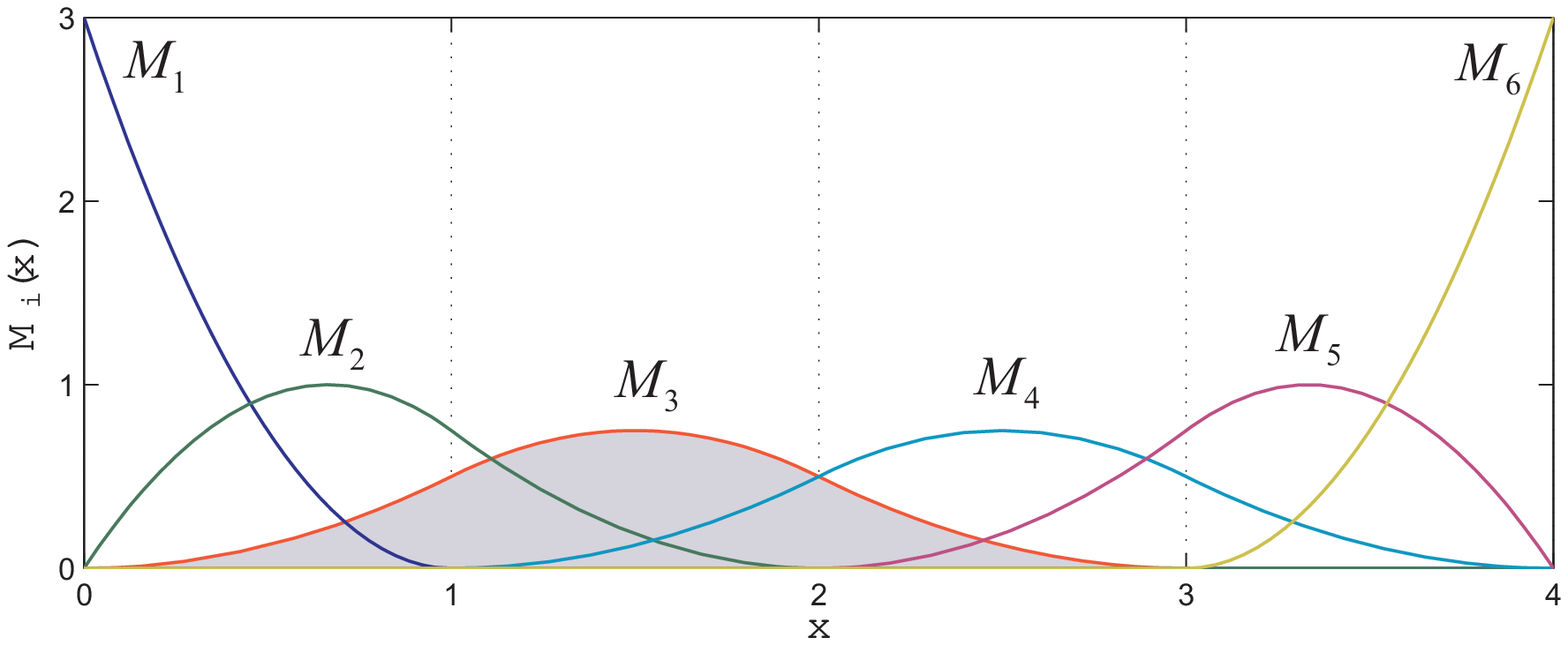}}	\\
	\subfigure[]{\includegraphics[trim = 0cm 1cm 0cm 1cm,clip,width=0.49\textwidth]{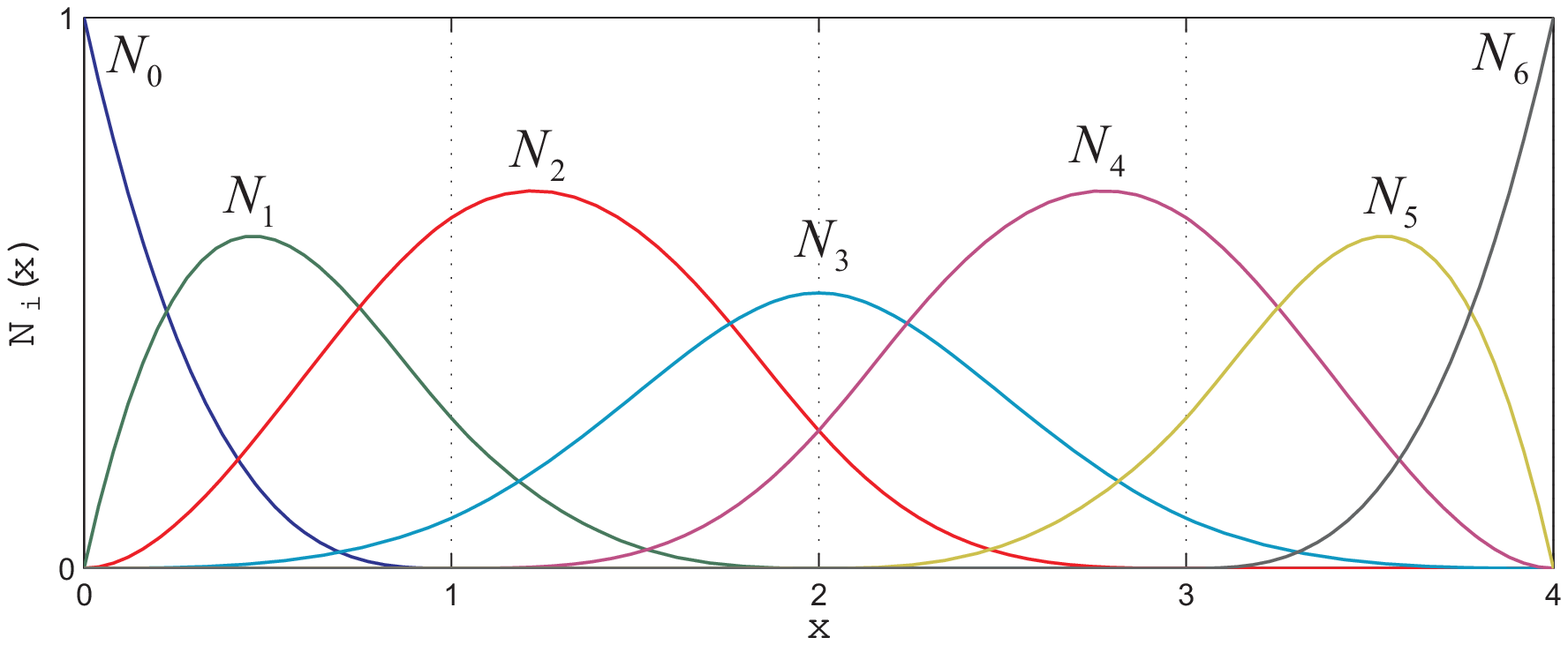}}
	\subfigure[]{\includegraphics[trim = 0cm 1cm 0cm 1cm,clip,width=0.49\textwidth]{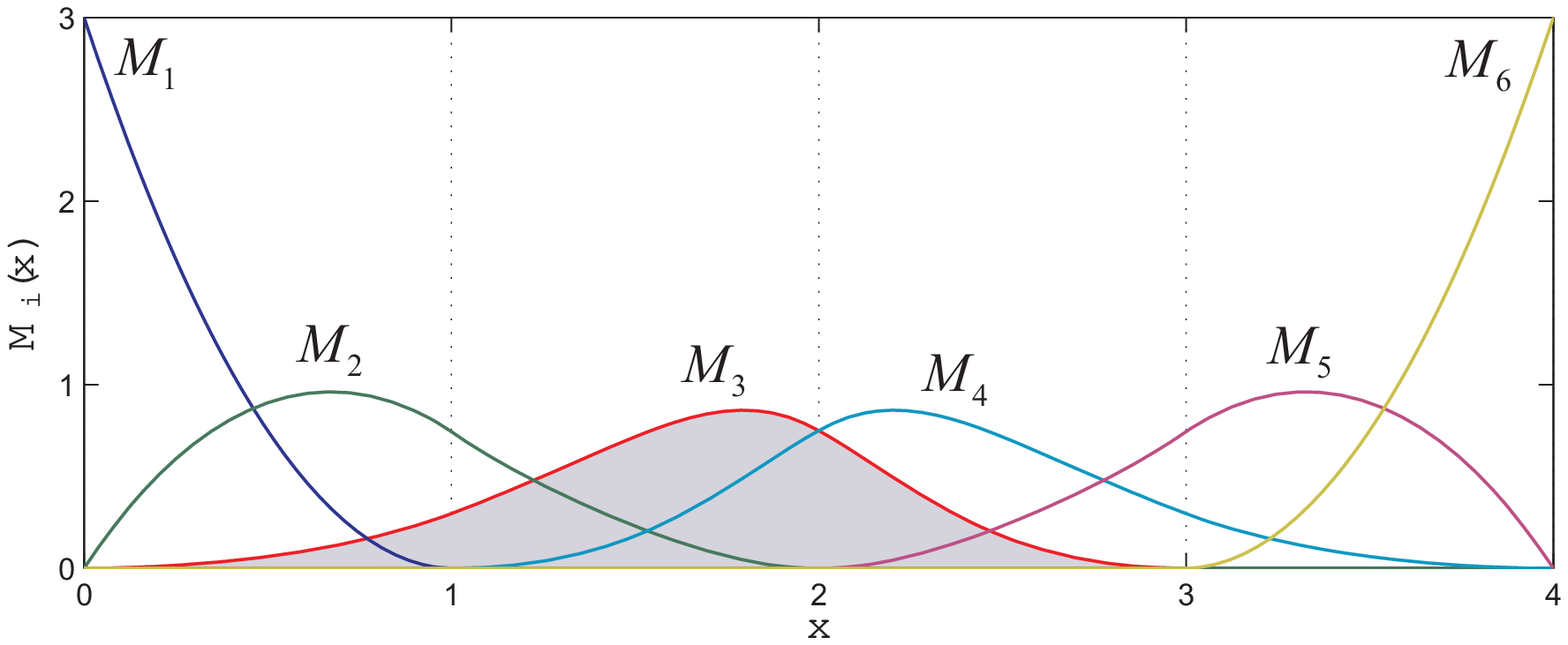}}	\\
%	\subfigure[]{\includegraphics[width=0.49\textwidth]{figures/NurbsNodalBasis4.eps}}
%	\subfigure[]{\includegraphics[width=0.49\textwidth]{figures/NurbsEdgeBasis4.eps}}	\\
%	\subfigure[]{\includegraphics[width=0.49\textwidth]{figures/NurbsNodalBasis5.eps}}
%	\subfigure[]{\includegraphics[width=0.49\textwidth]{figures/NurbsEdgeBasis5.eps}}
	\caption{Cubic NURBS basis functions (left) and associated quadratic edge functions (right) corresponding to a knot vector $\KnotVector{t}{}{0,0,0,0,1,2,3,4,4,4,4}$ under the influence of a different choice in weights. Every weight $w_i$ associated with basis function $N_i(x)$ is set equal to one, except $w_3$: Figure (a) and (b) correspond to $w_3=1$ and thus simplify to B-splines; (c) and (d) correspond to $w_3=1/2$. Note that the 'node' NURBS basis features partition of unity $\sum_{i=0}^n N_{i}(x) = 1$ on $[0,4]$, while the edge functions feature unit integral $\int_0^4 M_{i}(x) = 1$.}
%	\caption{Cubic NURBS basis functions (left) and associated quadratic edge functions (right) corresponding to a knot vector $\KnotVector{t}{}{0,0,0,0,1,2,3,4,4,4,4}$ under the influence of a different choice in weights. All weights $\left\{ w_i\right\}_{i=0}^{n}$ are set all equal to one, except $w_3$: Figure (a) and (b) correspond to $w_3=1$ and thus simplify to B-splines; (c) and (d) correspond to $w_3=1/2$; (e) and (f) correspond to $w_3=1/4$; and (g) and (h) correspond to $w_2=0$. Note that the 'nodal' NURBS basis features partition of unity $\sum_{i=0}^n N_{i}(x) = 1$ on $[0,4]$, while the edge functions feature unit integral $\int_0^4 M_{i}(x) = 1$.}
	\label{fig:NurbsEdgeFunctions}
\end{figure}

\vspace{0.25cm}

Important is to realize that the derivation of the edge basis functions is valid for any basis $\left\{N_i(x)\right\}_{i=0}^n$ (polynomial or non-polynomial) that is linear independent and is a partition of unity. Examples are Lagrange polynomials, Bezier, B-splines and NURBS. If the basis is nodal, which is the case for Lagrange polynomials, then the $\bar{u}_i$ have geometric as well as physical meaning since they represent the line integrals along the mesh edges $e_i$. Edge functions based on Lagrange polynomials are used in the Mimetic framework presented in \cite{Gerritsma:2011,Kreeft:2012a,Kreeft:2012b,Kreeft:2011}. In case of a non-nodal basis $\left\{N_i(x)\right\}_{i=0}^n$, such as Bezier polynomials, B-splines and NURBS, the degrees of freedom $\bar{u}_i$ only have a geometric interpretation, since they are discrete values associated with mesh edges $e_i$. They are however not related to differential forms by means of integration. 

Compatible spaces of B-splines have already been used in \cite{Buffa:2011b,Buffa:2010,Evans:2012,Ratnani:2012}. They make use of the fact that the derivative of a B-spline of order $p$ is a B-spline of order $p-1$ and can be written in the following form \cite{DeBoor:1978},
\begin{align} 
	\diff \sum_{i=0}^n \bar{T}_i \; B_{i,p}(x) = \sum_{i=1}^n \left(\bar{T}_i - \bar{T}_{i-1} \right) \; c_{i} \cdot B_{i,p-1}(x),    
\end{align}
where the $c_i$ are constants that depend on $p$, and the regularity of the mesh. We can directly observe that using B-splines, the 'node' functions, $N_{i}(x)$, are given by $N_{i}(x) = B_{i,p}(x)$, while the 'edge' functions, $M_i(x)$, are given by $M_{i}(x) = c_i \cdot B_{i,p-1}(x)$. The latter are called Curry Schoenberg B-splines, which was in fact the initial B-spline representation derived by \citet{Schoenberg:1946}. It was later that B-splines were scaled to form a partition of unity. 

In this paper we derive compatible spaces of discrete differential forms from NURBS \cite{Piegl:1997}, non-uniform rational B-splines. NURBS are rational functions of B-splines and are given as,
\begin{align}
	N_i(x) = \frac{w_i \cdot B_{i,p}(x)}{\sum_{i=0}^n w_i \cdot B_{i,p}(x)}.
\end{align}
NURBS trivially satisfy the required partition of unity, and can therefore be used as a basis to derive compatible spaces of discrete differential forms. Figure \ref{fig:NurbsEdgeFunctions} shows an example of cubic NURBS basis functions and their corresponding edge functions for different choices of the NURBS weights $w_i$.

\subsection{Multivariate discrete differential forms}
A basis for differential forms in $n$-dimensional space is simply obtained by means of the tensor product of the derived univariate 'node' and 'edge' functions. A basis for 2-forms in 2D is for example obtained by applying the tensor product of edge functions in both directions. We can define the following finite dimensional approximation spaces for 0-, 1-, and 2-forms in $\mathbb{R}^2$,
\begin{align*}
	\Lambda^0_h \left(\Omega \right) &= 																																	  \Span{N_{i}(\x{1}) \otimes N_{j}(\x{2})}_{i=0,j=0}^{n_1,n_2} \\
	\Lambda^1_h \left(\Omega \right) &= 																																		\Span{M_{i}(\x{1}) \otimes N_{j}(\x{2})}_{i=1,j=0}^{n_1,n_2} \times 
						 	\Span{N_{i}(\x{1}) \otimes M_{j}(\x{2})}_{i=0,j=1}^{n_1,n_2} \\
	\Lambda^2_h \left(\Omega \right) &= 																																	 	\Span{M_{i}(\x{1}) \otimes M_{j}(\x{2})}_{i=1,j=1}^{n_1,n_2} 
\end{align*}

\begin{figure}
	\centering
		\subfigure[Inner oriented cell complex]{\includegraphics[trim = 0cm 1cm -1cm 1.5cm,clip,width=0.60\textwidth]{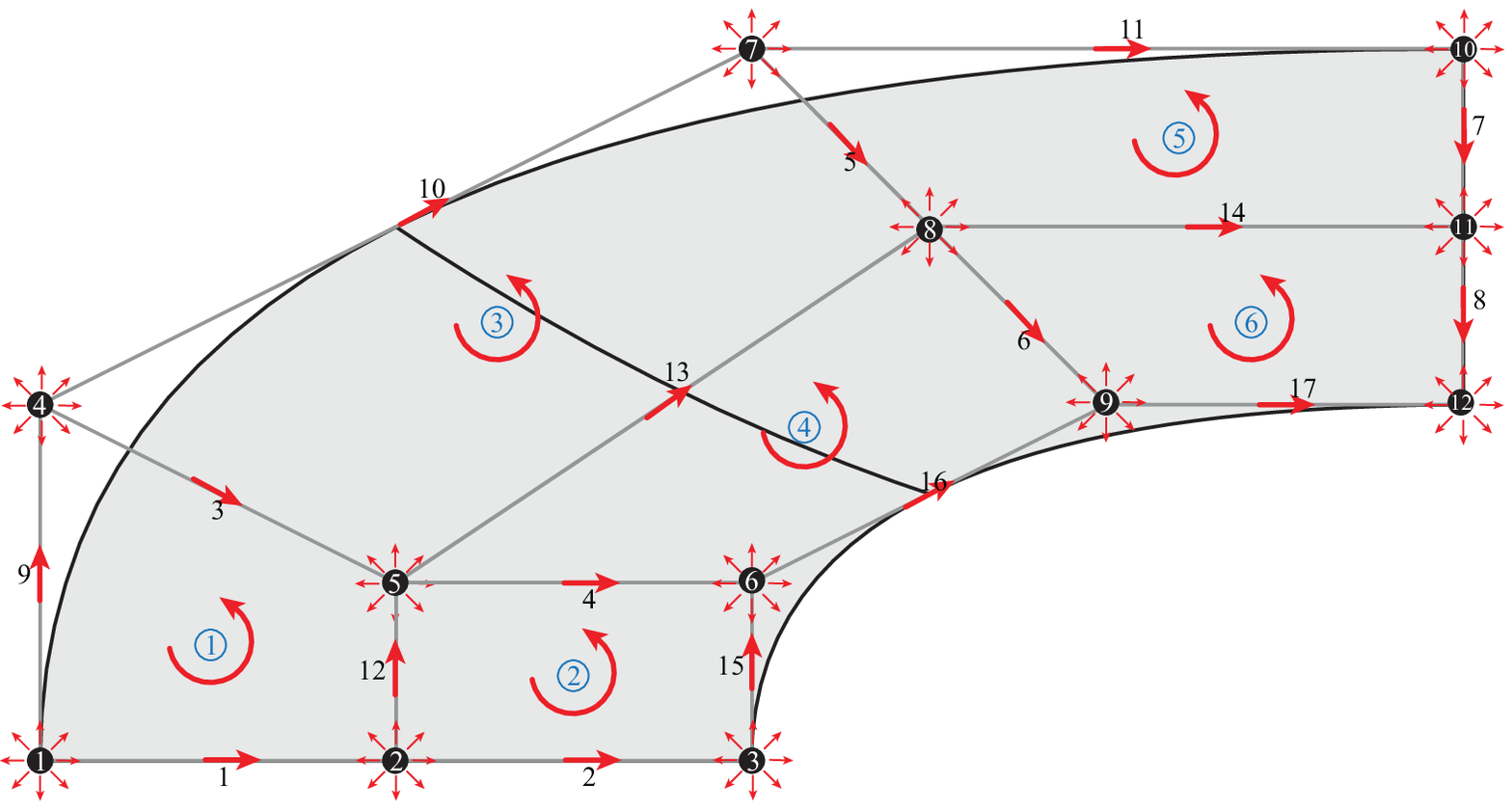}}
		\subfigure[Inner oriented 'node' function]{\includegraphics[trim = 0cm 1cm -2cm 1cm,clip,width=0.45\textwidth]{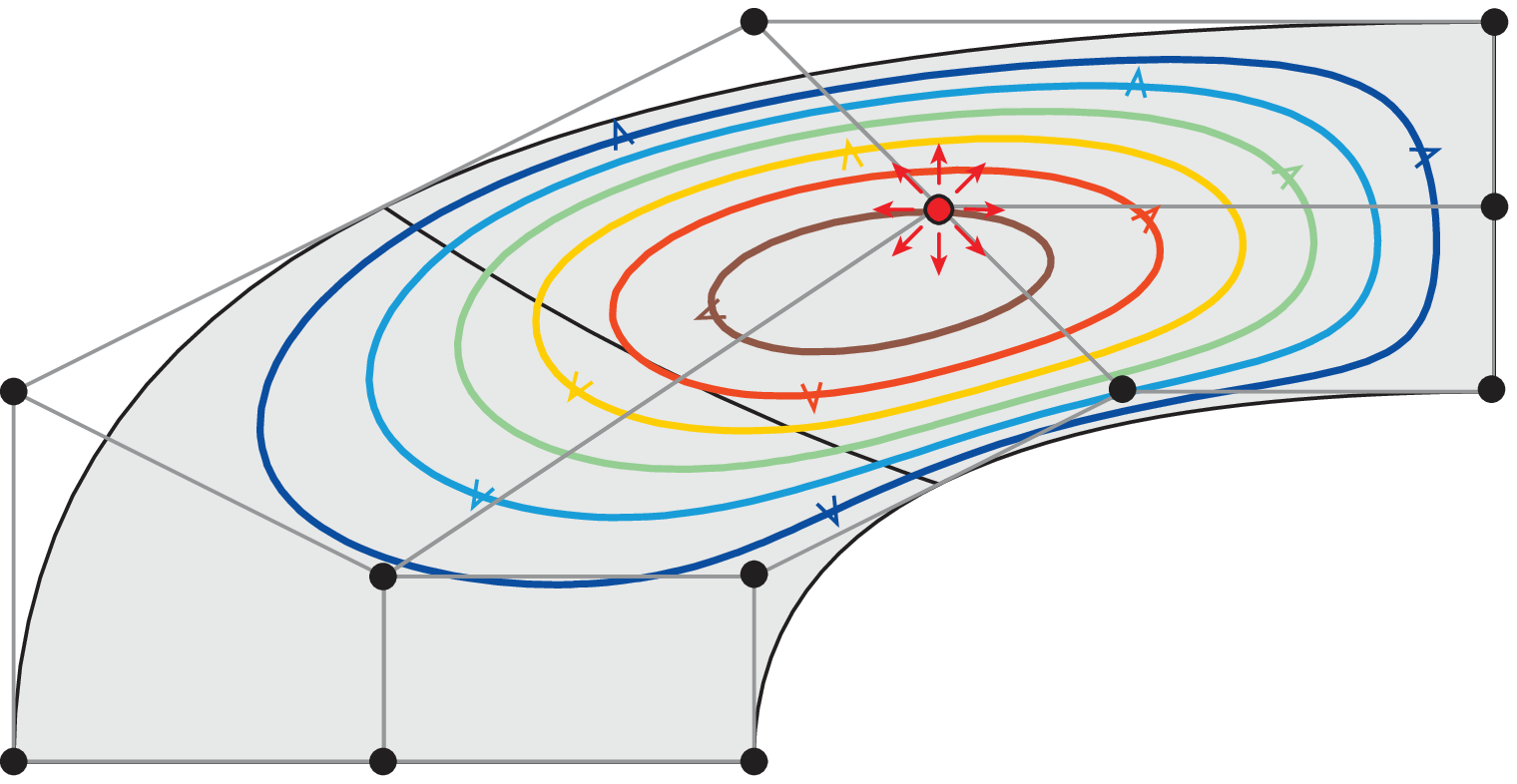}}
		\subfigure[Inner oriented 'node' function]{\includegraphics[trim = -2cm 1cm 0cm 1cm,clip,width=0.45\textwidth]{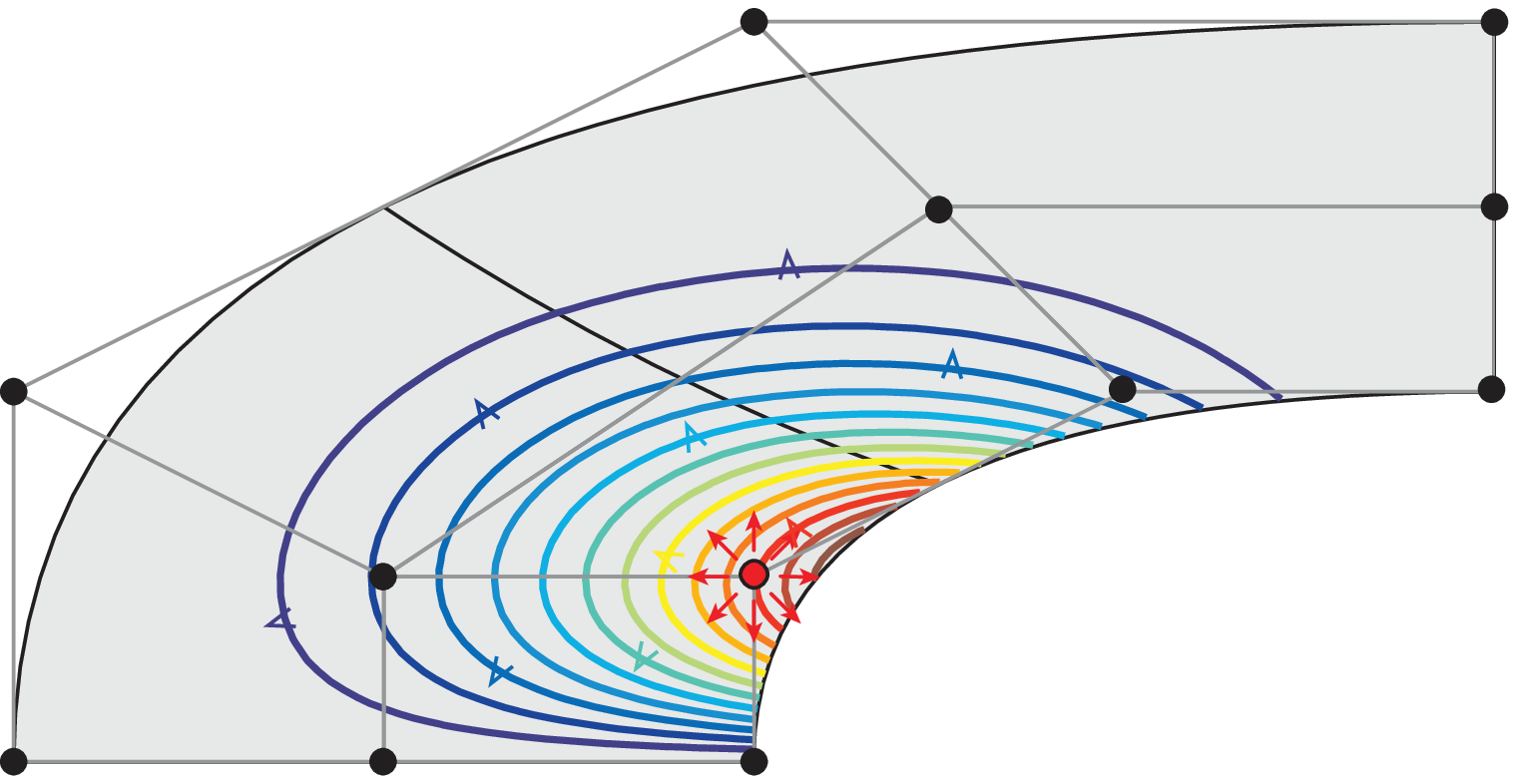}}
		\subfigure[Inner oriented 'edge' function]{\includegraphics[trim = 0cm 1cm -2cm 1cm,clip,width=0.45\textwidth]{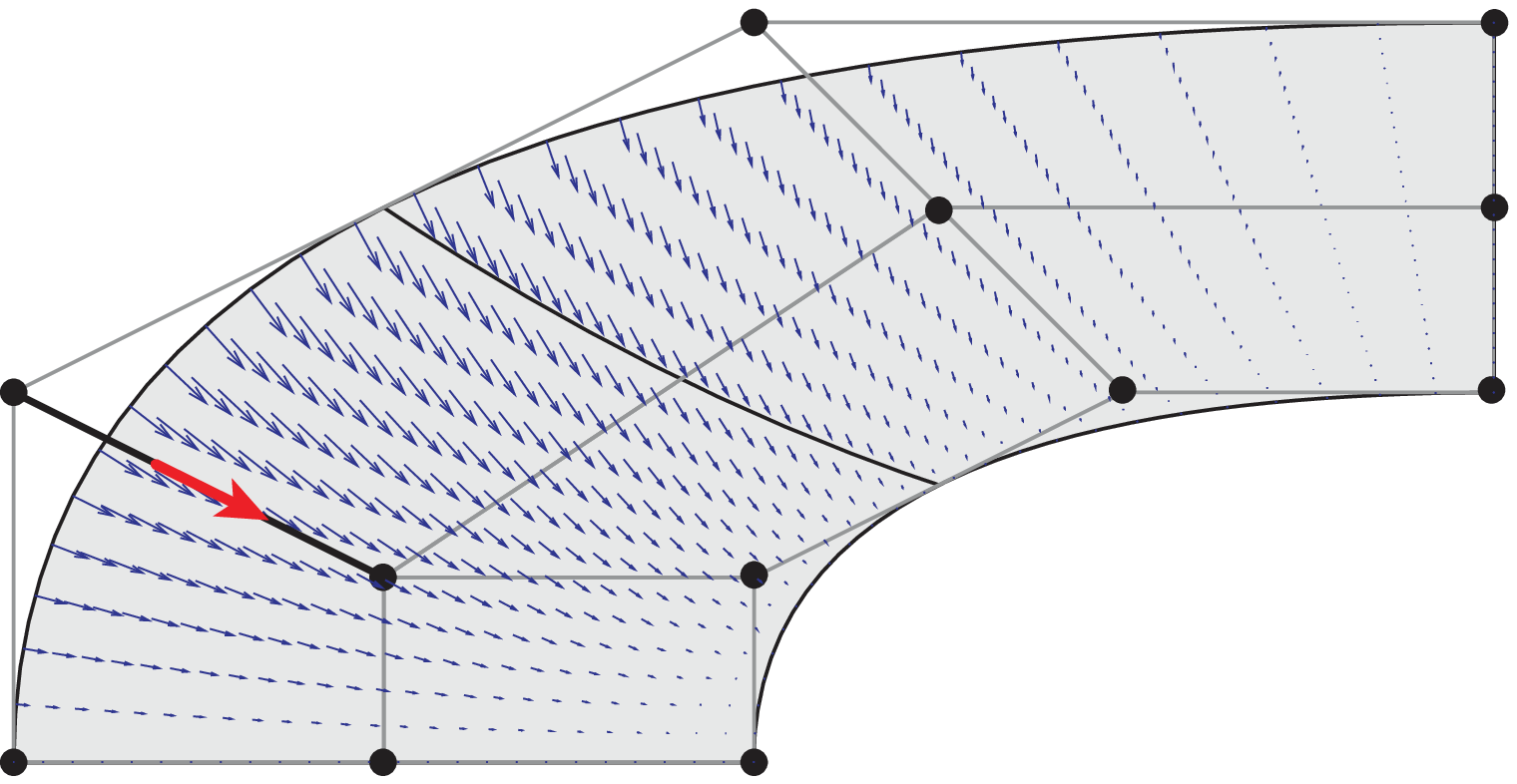}}
		\subfigure[Inner oriented 'edge' function]{\includegraphics[trim = -2cm 1cm 0cm 1cm,clip,width=0.45\textwidth]{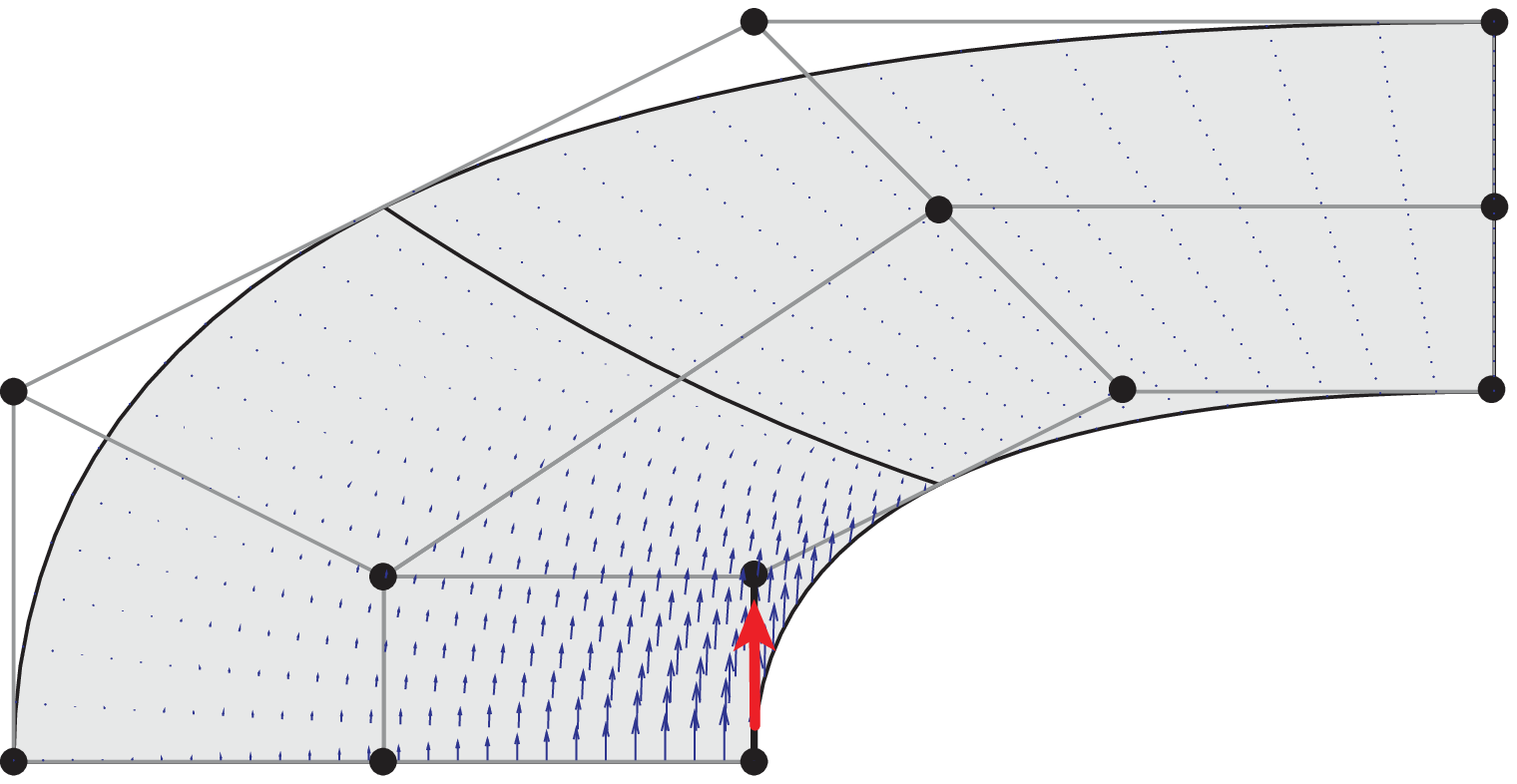}}
		\subfigure[Inner oriented 'face' function]{\includegraphics[trim = 0cm 1cm -2cm 1cm,clip,width=0.45\textwidth]{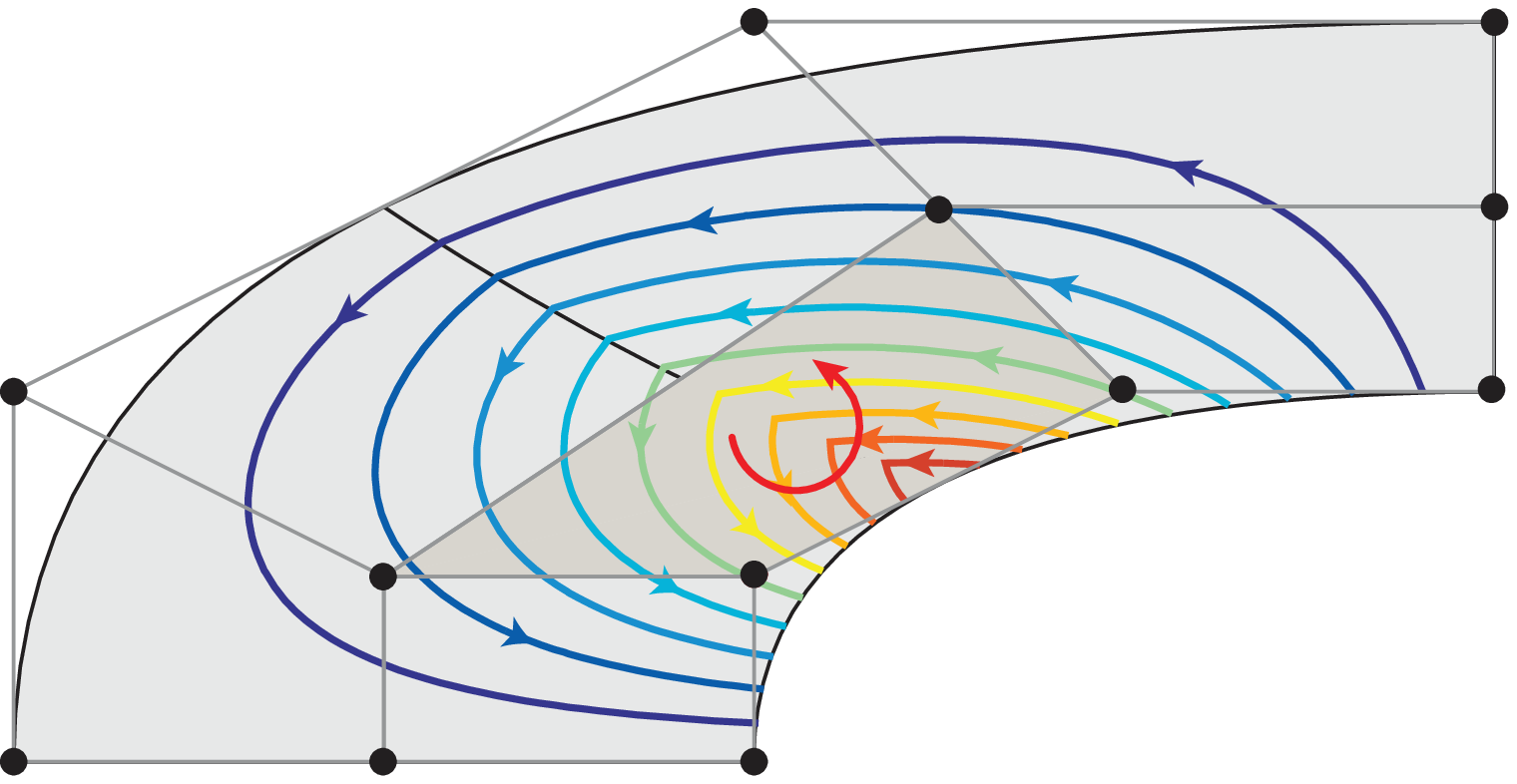}}
		\subfigure[Inner oriented 'face' function]{\includegraphics[trim = -2cm 1cm 0cm 1cm,clip,width=0.45\textwidth]{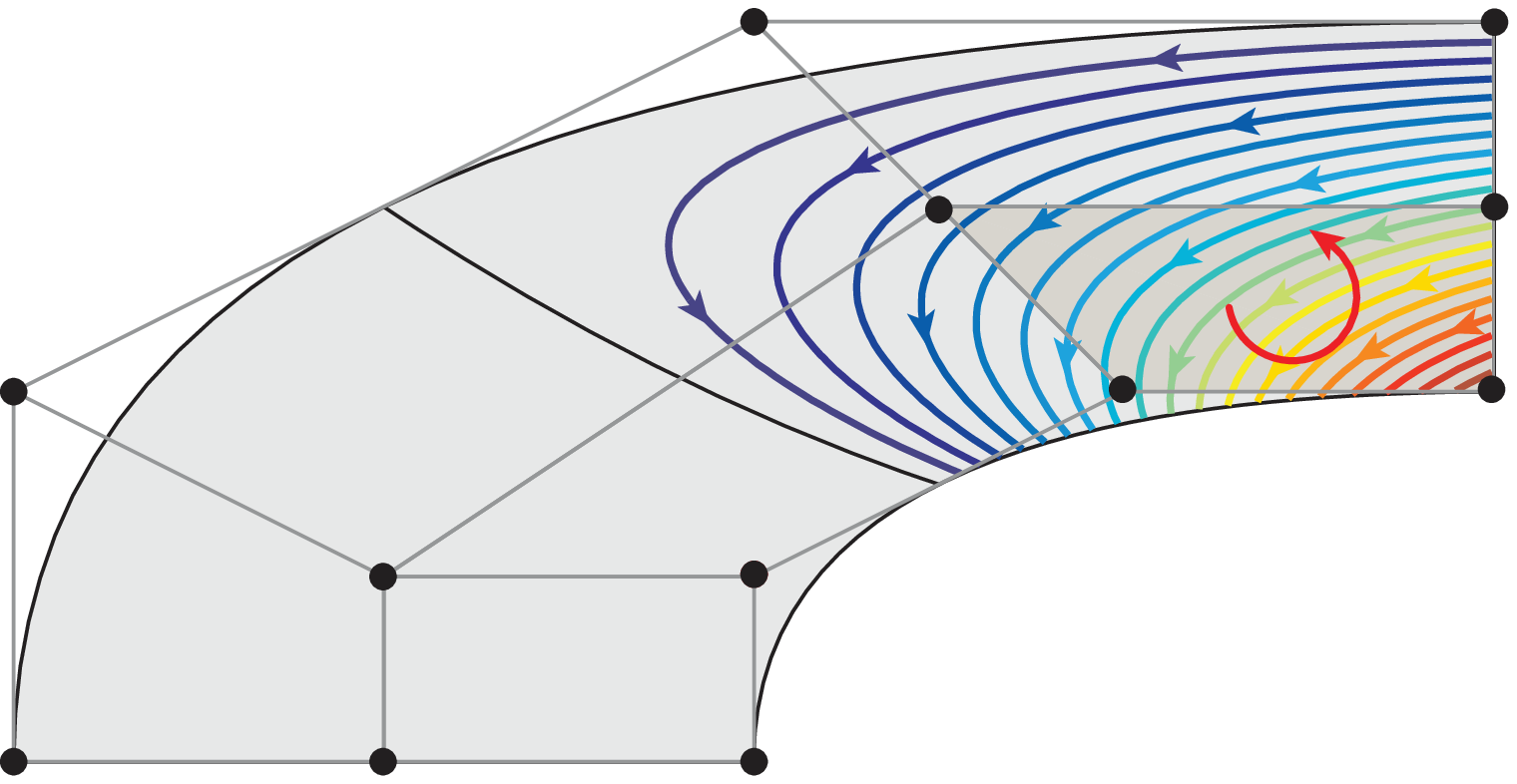}}
	\caption{2D examples of inner oriented 'node' (b and c), 'edge' (d and e) and 'face' (f and g) basis functions. The basis functions are derived from NURBS of bi-degree 2 and knot vectors $\KnotVector{U}{1}{0,0,0,1,1,1}$ and $\KnotVector{U}{2}{0,0,0,0.5,1,1,1}$ with weights equal to 1 (NURBS reduces to B-spline).}
	\label{fig:Reconstruction}
\end{figure}

\begin{figure}
	\centering
		\subfigure[Outer oriented cell complex]{\includegraphics[trim = 0cm 1cm -1cm 1.5cm,clip,width=0.60\textwidth]{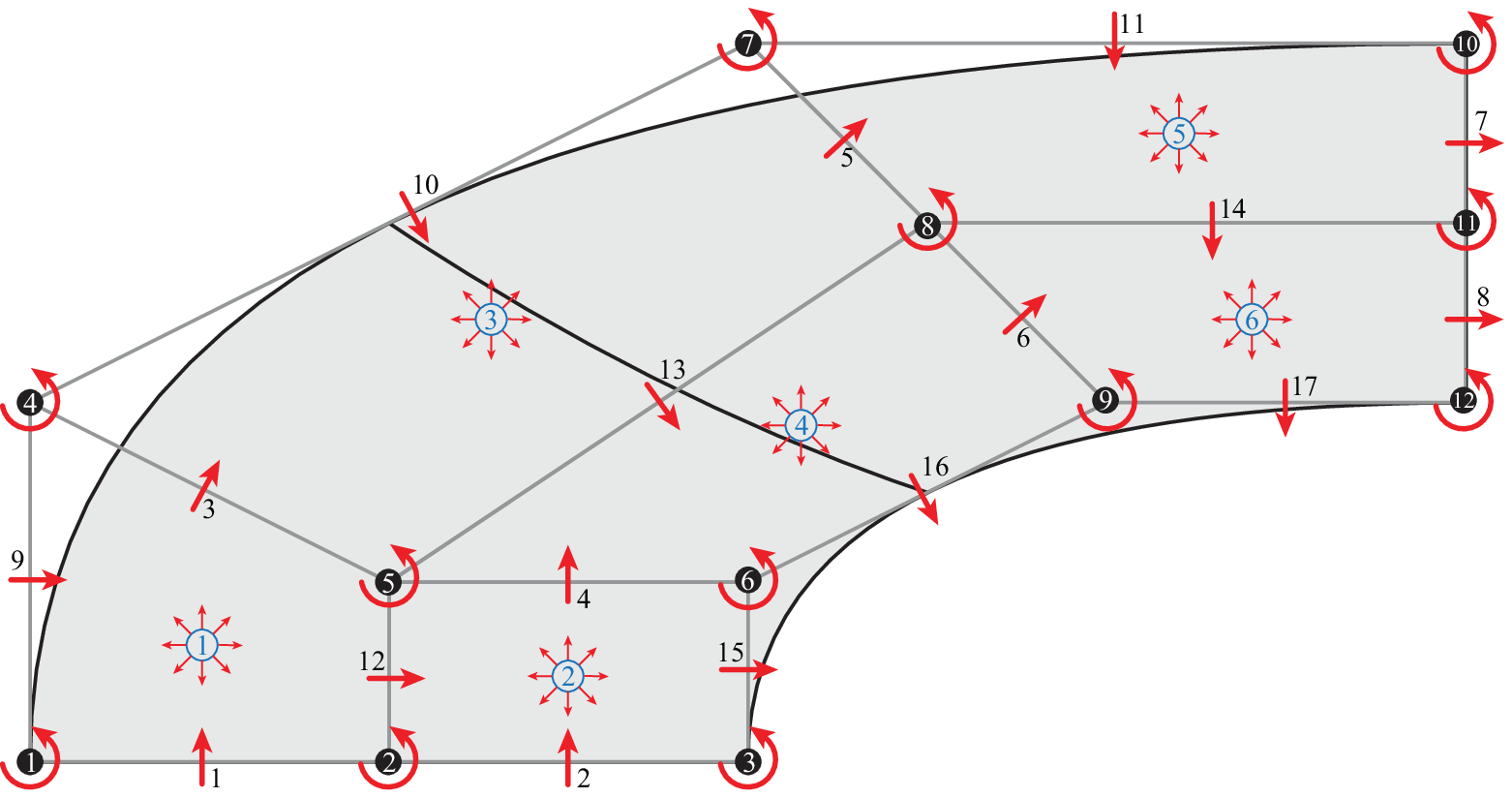}}
		\subfigure[Outer oriented 'node' function]{\includegraphics[trim = 0cm 1cm -2cm 1cm,clip,width=0.45\textwidth]{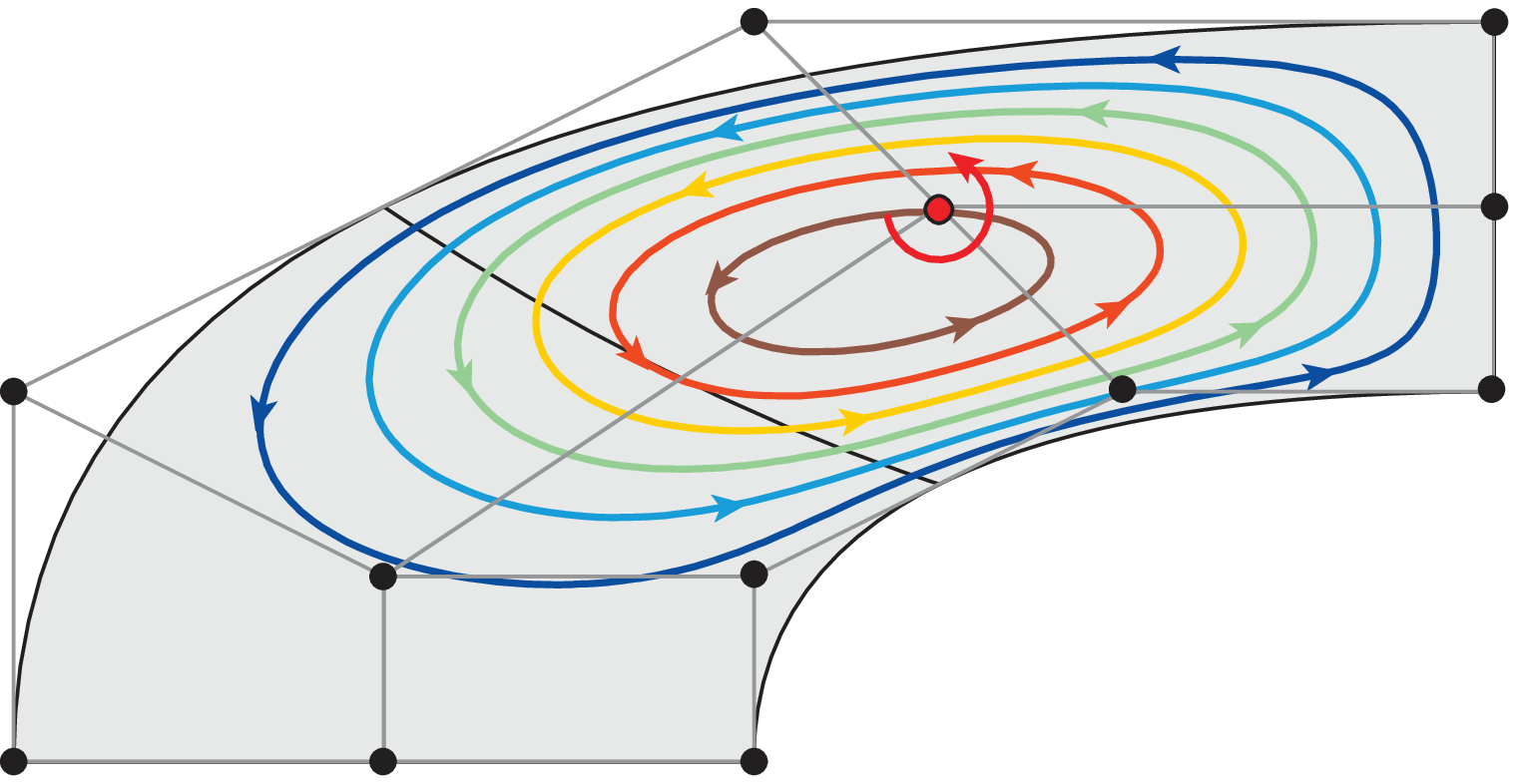}}
		\subfigure[Outer oriented 'node' function]{\includegraphics[trim = -2cm 1cm 0cm 1cm,clip,width=0.45\textwidth]{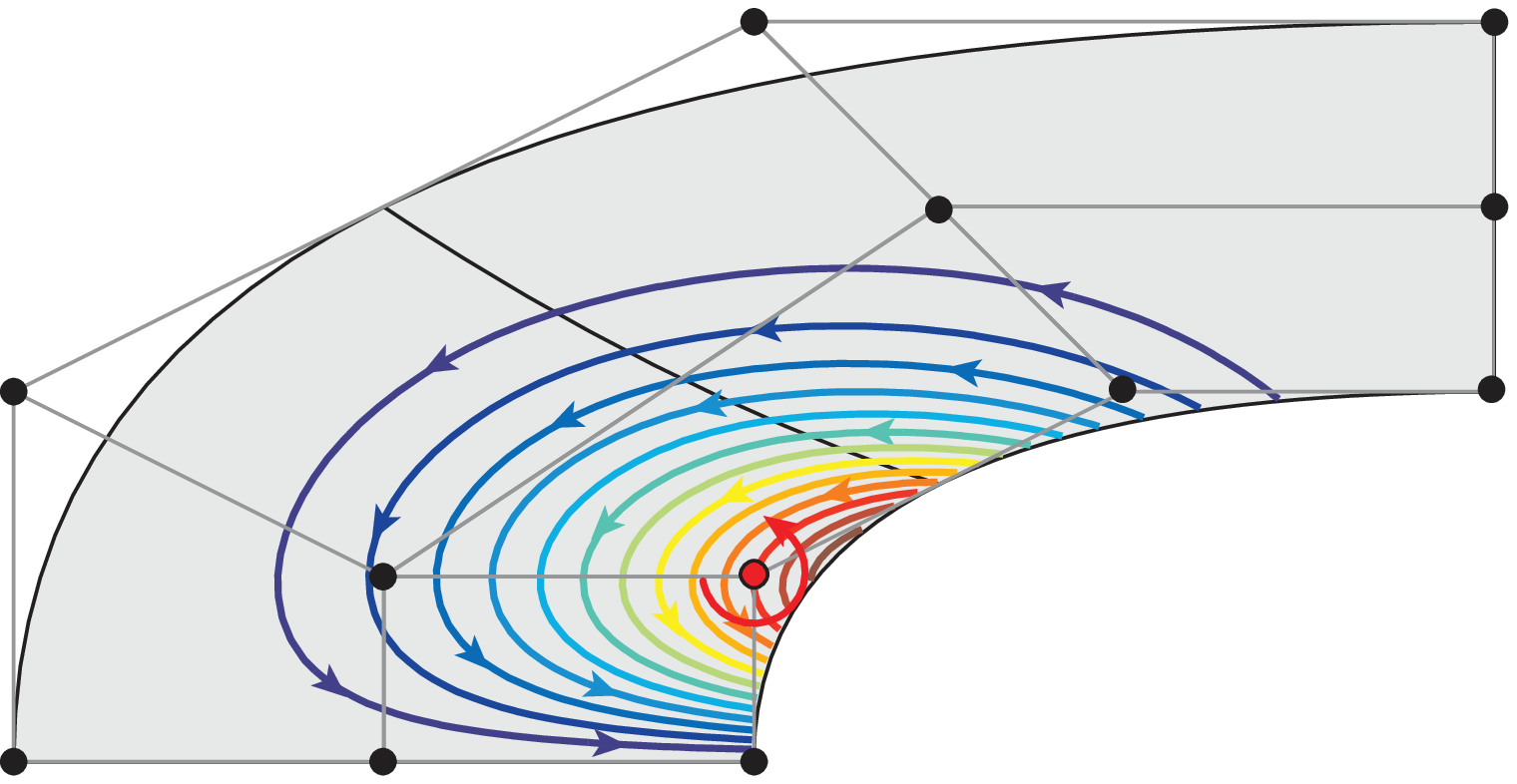}}
		\subfigure[Outer oriented 'edge' function]{\includegraphics[trim = 0cm 1cm -2cm 1cm,clip,width=0.45\textwidth]{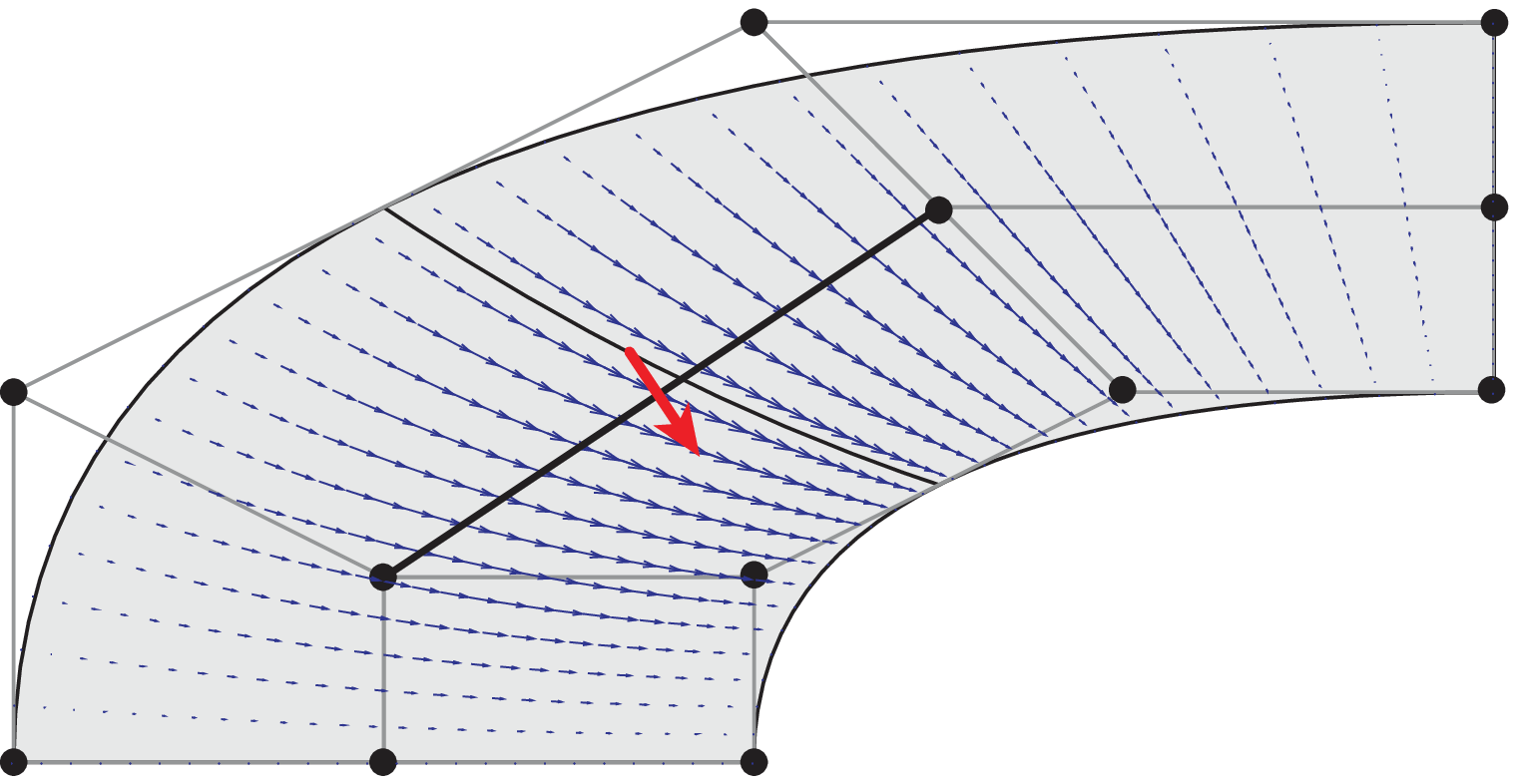}}
		\subfigure[Outer oriented 'edge' function]{\includegraphics[trim = -2cm 1cm 0cm 1cm,clip,width=0.45\textwidth]{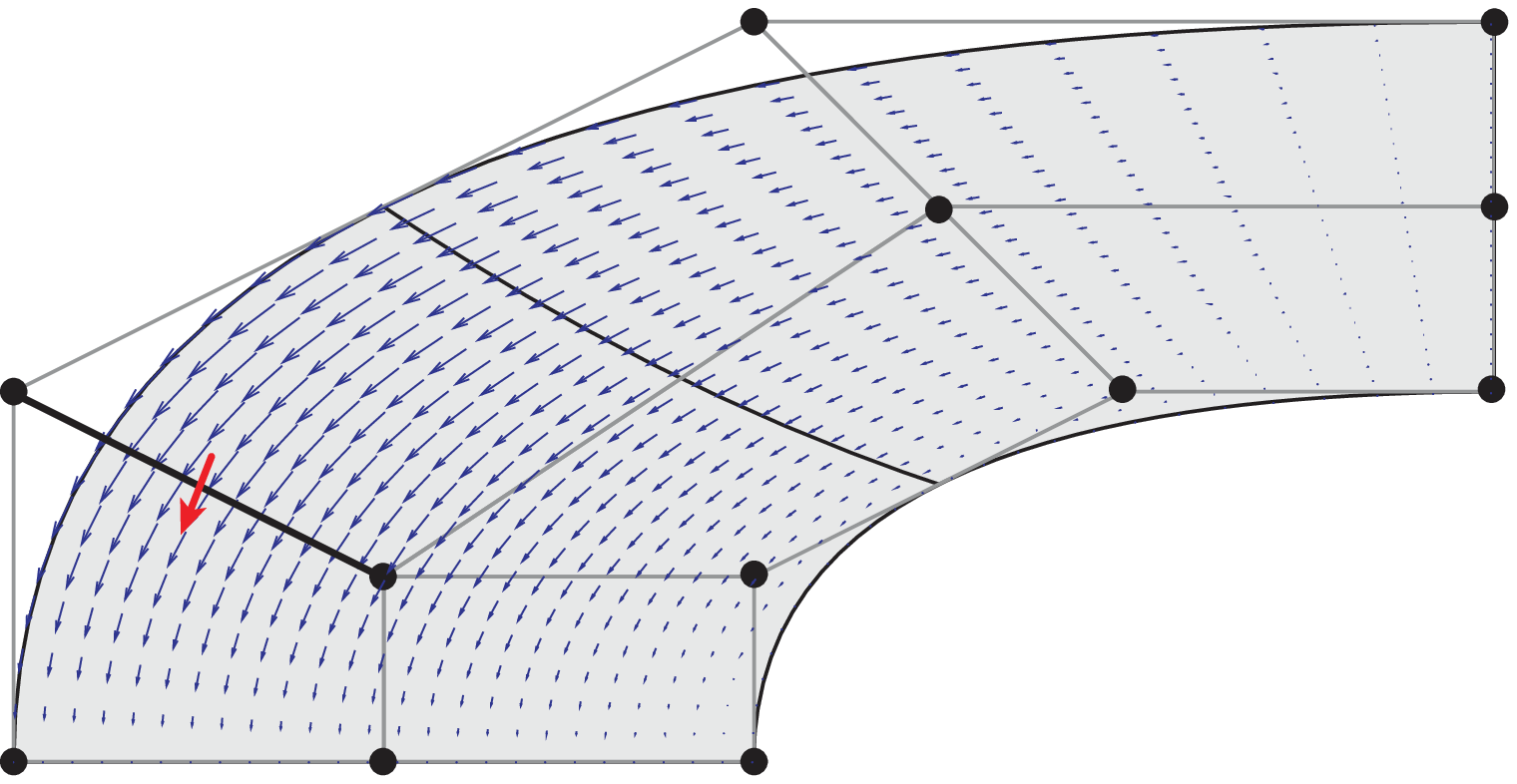}}
		\subfigure[Outer oriented 'face' function]{\includegraphics[trim = 0cm 1cm -2cm 1cm,clip,width=0.45\textwidth]{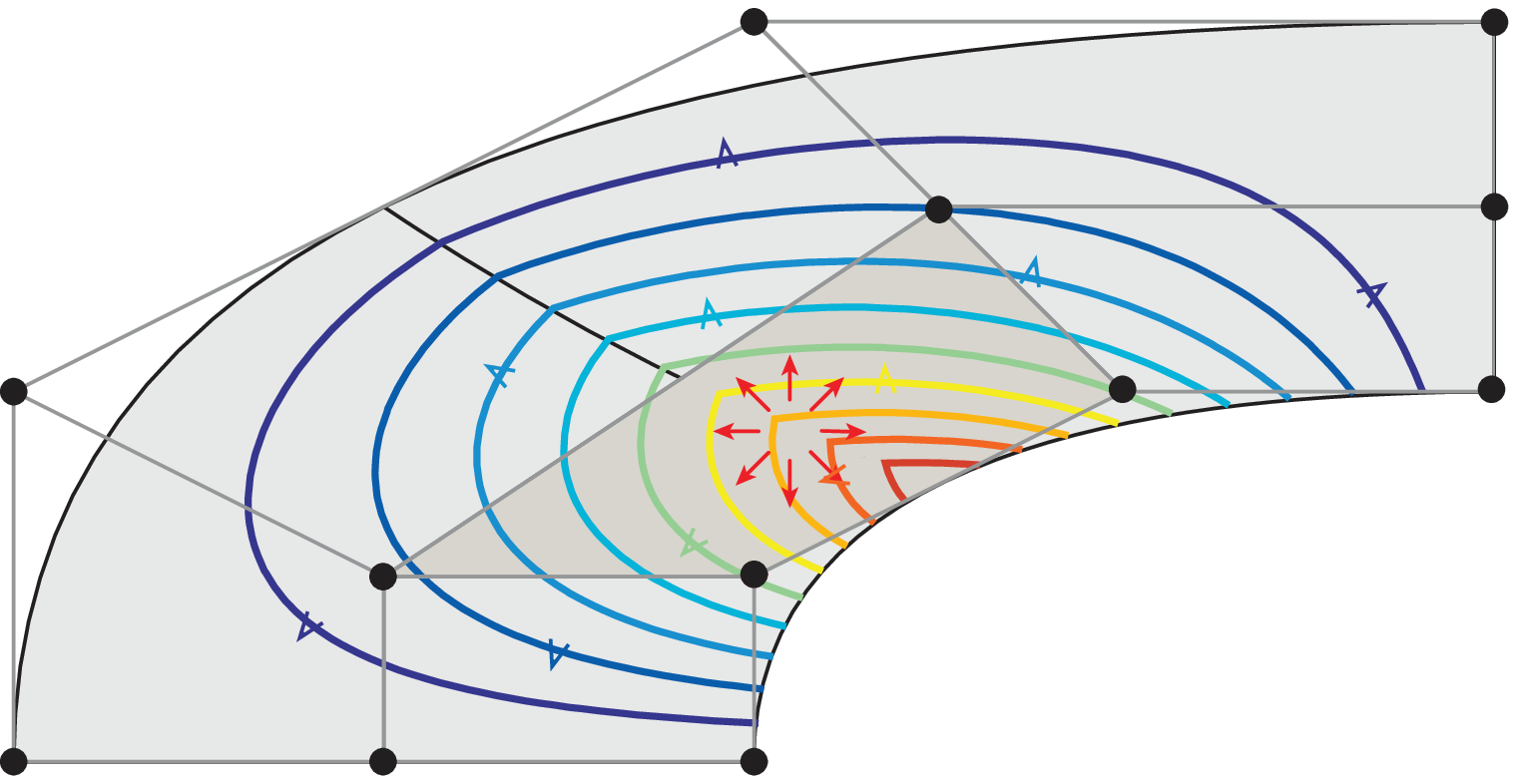}}
		\subfigure[Outer oriented 'face' function]{\includegraphics[trim = -2cm 1cm 0cm 1cm,clip,width=0.45\textwidth]{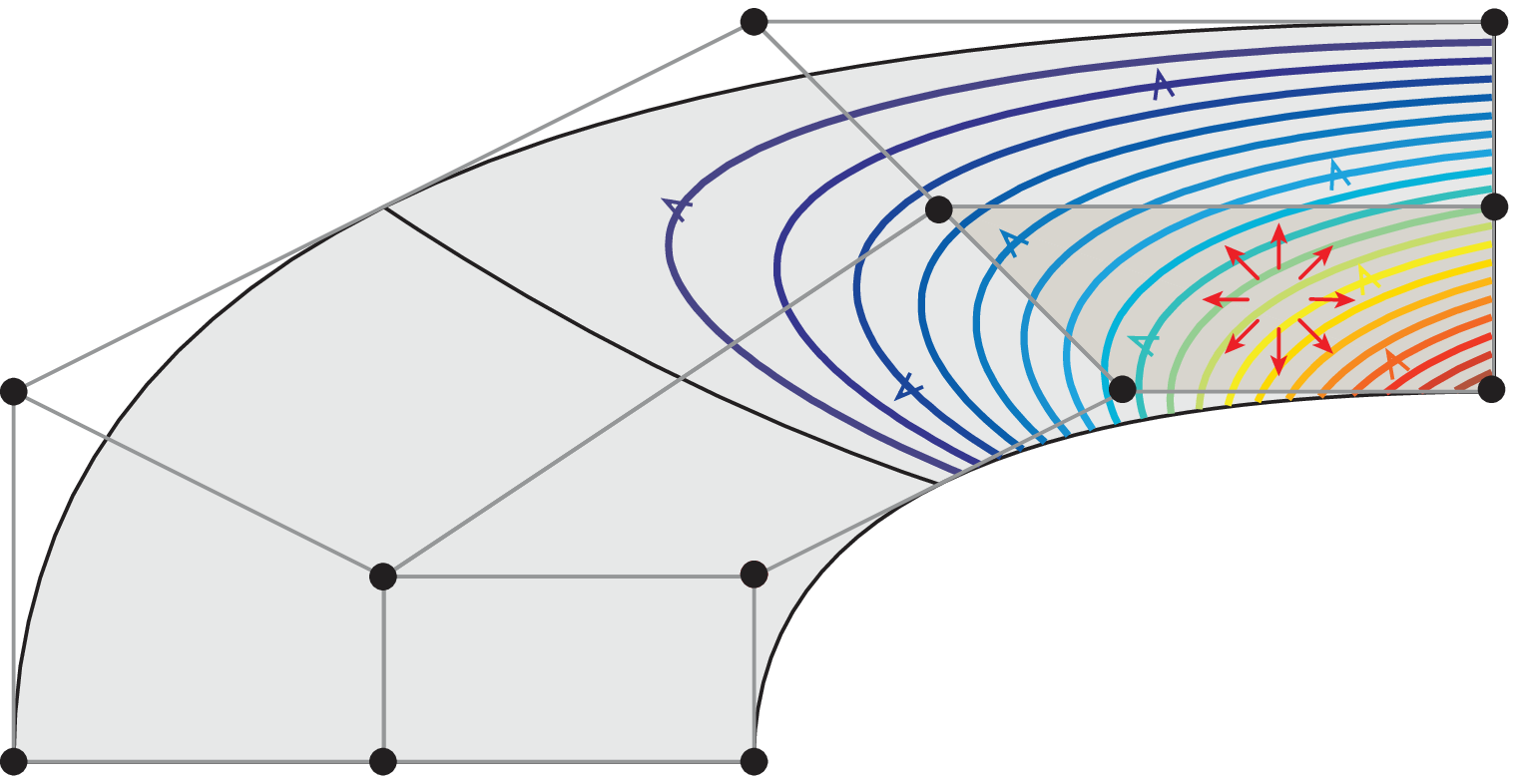}}
	\caption{2D examples of outer oriented 'node' (b and c), 'edge' (d and e) and 'face' (f and g) basis functions. The basis functions are derived from NURBS of bi-degree 2 and knot vectors $\KnotVector{U}{1}{0,0,0,1,1,1}$ $\KnotVector{U}{2}{0,0,0,0.5,1,1,1}$ with weights equal to 1 (NURBS reduces to B-spline).}
	\label{fig:Reconstruction2}
\end{figure}

Examples of 2-dimensional basis functions associated with inner and outer oriented points, edges, and faces are displayed in Figure \ref{fig:Reconstruction} and \ref{fig:Reconstruction2}. These spaces follow the following sequence, which is exact on contractible domains.
%\begin{multicols}{2}
\begin{align}
\begin{diagram}
	\Lambda^{0}(\Omega)	   &  \rTo^{\diff}  & \Lambda^{1}(\Omega)    & 	\rTo^{\diff} & \Lambda^{2}(\Omega)				\\
	\dTo^{\pi_h} 				   &                & \dTo^{\pi_h} 					 &               & \dTo^{\pi_h} 							\\ 
	\Lambda^{0}_h(\Omega)	 &  \rTo^{\diff}  & \Lambda^{1}_h(\Omega)  & 	\rTo^{\diff} & \Lambda^{2}_h(\Omega)
\end{diagram}
\end{align}

%\begin{align}
%\begin{diagram}
%	\Lambda^{0}(\Omega)	   &  \rTo^{\diff}_{\grad^{\perp}}  & \Lambda^{1}(\Omega)    & 	\rTo^{\diff}_{\div} & \Lambda^{2}(\Omega)				\\
%	\dTo^{\pi_h} 				   &                & \dTo^{\pi_h} 					   &               & \dTo^{\pi_h} 							\\ 
%	\Lambda^{0}_h(\Omega)	 &  \rTo^{\diff}_{\grad^{\perp}}  & \Lambda^{1}_h(\Omega)  & 	\rTo^{\diff}_{\div} & \Lambda^{2}_h(\Omega)
%\end{diagram}
%\end{align}
%%\end{multicols}
%Note that inner oriented variables feature the regular gradient and the curl operator which in 2D is usually denoted by $\rot$. The outer oriented variables have the perpendicular gradient $\grad^{\perp}$ and the regular 2D divergence operator.
%
%Some examples of finite dimensional differential forms in 2D are given in example \ref{ex:Dforms1}. Remember that in $\mathbb{R}^2$, 0- and 2-forms are scalar valued functions and 1-forms are vector valued functions.
%\begin{example}[Discrete differential forms]
%Examples of the finite dimensional reconstructions of 0-, 1-, 2-form fields.
%\begin{align*}
%	\text{0-form:} \quad & \Dform{T}{(0)}_h(\vect{x}{}) = 
%								 \sum_{i,j} \bar{T}_{i,j} \; N_{i}(\x{1})\cdot N_{j}(\x{2})  \\
%	\text{1-form:} \quad & \Dform{u}{(1)}_h(\vect{x}{}) = 
%								 \sum_{i,j} \bar{u}^{(1)}_{i,j} \; M_{i}(\x{1})\cdot N_{j}(\x{2}) + 
%								 \sum_{i,j} \bar{u}^{(2)}_{i,j} \; N_{i}(\x{1})\cdot M_{j}(\x{2}) \\
%	\text{2-form:} \quad & \Dform{f}{(2)}_h(\vect{x}{}) = 
%								 \sum_{i,j} \bar{f}_{i,j} \; M_{i}(\x{1})\cdot M_{j}(\x{2})  \\
%\end{align*}
%\label{ex:Dforms1}
%\end{example}

We will more deeply discuss the discrete spaces of differential forms in $\mathbb{R}^3$ to stay close to our discussion in section \ref{sec:Introduction}. We introduce the following notation for 3-dimensional tensor product basis functions,
\begin{align*}
	& B_{i,j,k}(\vect{x}{})					 = N_i(\x{1})  N_j(\x{2})  N_k(\x{3})  &  &&  && 	\\
	& B^{(1)}_{i,j,k} (\vect{x}{})	 = M_i(\x{1})  N_j(\x{2})  N_k(\x{3}), &  
	& B^{(2)}_{i,j,k} (\vect{x}{}) 	 = N_i(\x{1})  M_j(\x{2})  N_k(\x{3}), & 
	& B^{(3)}_{i,j,k} (\vect{x}{}) 	 = N_i(\x{1})  N_j(\x{2})  M_k(\x{3})  & 	\\
	& B^{(2,3)}_{i,j,k} (\vect{x}{}) = N_i(\x{1})  M_j(\x{2})  M_k(\x{3}), &
	& B^{(3,1)}_{i,j,k} (\vect{x}{}) = M_i(\x{1})  N_j(\x{2})  M_k(\x{3}), &
	& B^{(1,2)}_{i,j,k} (\vect{x}{})	= M_i(\x{1})  M_j(\x{2})  N_k(\x{3}) &	\\
	& B^{(1,2,3)}_{i,j,k} (\vect{x}{}) = M_i(\x{1}) M_j(\x{2})  M_k(\x{3}) &  &&  &&
\end{align*}
Note that the superscripts indicate in which direction an edge function is used. We can then define the following finite dimensional approximation spaces for 0-, 1-, 2- and 3-forms,
\begin{align*}
	\Lambda^0_h \left(\Omega \right) &= \Span{B_{i,j,k}(\vect{x}{})}_{i=0,j=0,k=0}^{n_1,n_2,n_3} \\
	\Lambda^1_h \left(\Omega \right) &= \Span{B^{(1)}_{i,j,k}(\vect{x}{})}_{i=1,j=0,k=0}^{n_1,n_2,n_3} \times 
																		 	\Span{B^{(2)}_{i,j,k}(\vect{x}{})}_{i=0,j=1,k=0}^{n_1,n_2,n_3} \times 
																		 	\Span{B^{(3)}_{i,j,k}(\vect{x}{})}_{i=0,j=0,k=1}^{n_1,n_2,n_3} \\
	\Lambda^2_h \left(\Omega \right) &= \Span{B^{(2,3)}_{i,j,k}(\vect{x}{})}_{i=0,j=1,k=1}^{n_1,n_2,n_3} \times 
																		 	\Span{B^{(3,1)}_{i,j,k}(\vect{x}{})}_{i=1,j=0,k=1}^{n_1,n_2,n_3} \times 
																		 	\Span{B^{(1,2)}_{i,j,k}(\vect{x}{})}_{i=1,j=1,k=0}^{n_1,n_2,n_3} \\
	\Lambda^3_h \left(\Omega \right) &= \Span{B^{(1,2,3)}_{i,j,k}(\vect{x}{})}_{i=1,j=1,k=1}^{n_1,n_2,n_3} 
\end{align*}

These spaces are constructed such that they satisfy the following sequence, which is exact on contractible domains,
\begin{align}
\begin{diagram}
	\Lambda^{0}(\Omega)	   &  \rTo^{\diff}_{\grad}  & \Lambda^{1}(\Omega)    & 	\rTo^{\diff}_{\curl} & \Lambda^{2}(\Omega)    &  	\rTo^{\diff}_{\div}  & \Lambda^{3}(\Omega) 					\\
	\dTo^{\pi_h} 				   &                        & \dTo^{\pi_h} 					   &                       & \dTo^{\pi_h} 					  &                        & \dTo^{\pi_h} 			 						\\ 
	\Lambda^{0}_h(\Omega)	 &  \rTo^{\diff}_{\grad}  & \Lambda^{1}_h(\Omega)  & 	\rTo^{\diff}_{\curl} & \Lambda^{2}_h(\Omega)  &  	\rTo^{\diff}_{\div}  & \Lambda^{3}_h(\Omega)
\end{diagram}
\end{align}

We emphasize that 0- and 3-forms are scalar valued functions and 1- and 2-forms are vector valued functions. Some examples of quantities such as temperature, velocity, vorticity and density are displayed in the following example.
\begin{example}[Reconstruction in $\mathrm{R}^3$]
Examples of the finite dimensional reconstructions of 0-, 1-, 2- and 3-form fields.
\begin{align*}
	\text{0-form:} \quad & \Dform{T}{(0)}_h(\vect{x}{}) = \sum_{i,j,k} \bar{T}_{i,j,k} \; B_{i,j,k}(\vect{x}{}) \\
	\text{1-form:} \quad & \Dform{u}{(1)}_h(\vect{x}{}) = 
																\sum_{i,j,k} \bar{u}^{(1)}_{i,j,k} \; B^{(1)}_{i,j,k}(\vect{x}{}) +
														    \sum_{i,j,k} \bar{u}^{(2)}_{i,j,k} \; B^{(2)}_{i,j,k}(\vect{x}{}) +  
														    \sum_{i,j,k} \bar{u}^{(3)}_{i,j,k} \; B^{(3)}_{i,j,k}(\vect{x}{}) \\
	\text{2-form:} \quad & \Dform{\omega}{(2)}_h(\vect{x}{}) = 
																\sum_{i,j,k} \bar{\omega}^{(1)}_{i,j,k} \; B^{(2,3)}_{i,j,k}(\vect{x}{}) +
														    \sum_{i,j,k} \bar{\omega}^{(2)}_{i,j,k} \; B^{(3,1)}_{i,j,k}(\vect{x}{}) +  
														    \sum_{i,j,k} \bar{\omega}^{(3)}_{i,j,k} \; B^{(1,2)}_{i,j,k}(\vect{x}{}) \\
	\text{3-form:} \quad & \Dform{f}{(3)}_h(\vect{x}{}) = \sum_{i,j,k} \bar{f}_{i,j,k} \; B^{(1,2,3)}_{i,j,k}(\vect{x}{})
\end{align*}
\label{ex:Dforms}
\end{example}

As in univariate space, the balance laws are completely independent of the basis functions and we can perform differentiation discretely. We thus have a discrete matrix representation of the gradient, curl and divergence operator which is exact, even on curved and coarse meshes, see Section \ref{sec:DiscreteModeling}. These relations hold for Lagrange, Bezier, B-spline, NURBS or any other basis $\left\{N_i(x)\right\}_{i=0}^n$ that forms a partition of unity. This means that these relations hold independent of a specific numerical method.

\begin{example}[Discrete gradient operator.] 
Consider $ \Dform{u}{1}_h(\vect{x}{}) = \grad \; \Dform{T}{0}_h(\vect{x}{})$ where $\Dform{u}{1}_h(\vect{x}{})$ and $\Dform{T}{0}_h(\vect{x}{})$ are of the form as in example \ref{ex:Dforms}. Then the degrees of freedom $\bar{u}^{(d)}_{i,j,k}$, where $d=1,2,3$ represents component direction, are associated with edges and are given by,
\begin{align*}
	& \bar{u}^{(1)}_{i,j,k} = \bar{T}_{i,j,k} - \bar{T}_{i-1,j,k}, &   
	& \bar{u}^{(2)}_{i,j,k} = \bar{T}_{i,j,k} - \bar{T}_{i,j-1,k}  & &\text{and}&
	& \bar{u}^{(3)}_{i,j,k} = \bar{T}_{i,j,k} - \bar{T}_{i,j,k-1}  &
\end{align*}
This discrete representation allows a matrix vector product $\vect{u}{} =  \Mat{D}{1,0} \; \vect{T}{}$ where $ \Mat{D}{1,0} = \Mat{E}{0,1}^T$ is an incidence matrix, see Example \ref{ex:IncidenceMatrix}. This discrete gradient is exact and does not depend on geometric transformations, as proven in (\ref{CD:PullBack}).
\end{example}

\begin{example}[Discrete curl operator.] 
Consider $ \Dform{\omega}{2}_h(\vect{x}{}) = \curl \; \Dform{u}{1}_h(\vect{x}{})$ where $\Dform{\omega}{2}_h(\vect{x}{})$ and $\Dform{u}{1}_h(\vect{x}{})$ are of the form as in example \ref{ex:Dforms}. Then the degrees of freedom $\bar{\omega}^{(d)}_{i,j,k}$ are associated with faces and are given by,
\begin{align*}
	& \bar{\omega}^{(1)}_{i,j,k} = \bar{u}^{(3)}_{i,j,k} - \bar{u}^{(3)}_{i,j-1,k} + \bar{u}^{(2)}_{i,j,k-1} - \bar{u}^{(2)}_{i,j,k} & 					 \\
	& \bar{\omega}^{(2)}_{i,j,k} = \bar{u}^{(1)}_{i,j,k} - \bar{u}^{(1)}_{i,j,k-1} + \bar{u}^{(3)}_{i-1,j,k} - \bar{u}^{(3)}_{i,j,k} & \nonumber \\
	& \bar{\omega}^{(3)}_{i,j,k} = \bar{u}^{(2)}_{i,j,k} - \bar{u}^{(2)}_{i-1,j,k} + \bar{u}^{(1)}_{i,j-1,k} - \bar{u}^{(1)}_{i,j,k} & \nonumber
\end{align*}
which allows a matrix vector product $\vect{\omega}{} = \Mat{D}{2,1} \; \vect{u}{}$ where $\Mat{D}{2,1} = \Mat{E}{1,2}^T$ is an incidence matrix, see Example \ref{ex:IncidenceMatrix}. This discrete curl is exact and does not depend on geometric transformations, as proven in (\ref{CD:PullBack}).
\end{example}

\begin{example}[Discrete divergence operator.] 
Consider $ \Dform{f}{3}_h(\vect{x}{}) = \div \; \Dform{\omega}{2}_h(\vect{x}{})$ where $\Dform{\omega}{2}_h(\vect{x}{})$ and $\Dform{f}{3}_h(\vect{x}{})$ are of the form as in example \ref{ex:Dforms}. Then the degrees of freedom $\bar{f}_{i,j,k}$ are associated with volumes and follow the discrete relation,
\begin{align*}
	\bar{f}_{i,j,k} = \bar{\omega}^{(1)}_{i,j,k} - \bar{\omega}^{(1)}_{i-1,j,k} + \bar{\omega}^{(2)}_{i,j,k} - \bar{\omega}^{(2)}_{i,j-1,k}  + \bar{\omega}^{(3)}_{i,j,k} - \bar{\omega}^{(3)}_{i,j,k-1}
\end{align*}
which allows a matrix vector product $\vect{f}{} =  \Mat{D}{3,2} \; \vect{\omega}{}$ where $ \Mat{D}{3,2} = \Mat{E}{2,3}^T$ is an incidence matrix, see Example \ref{ex:IncidenceMatrix}. This discrete divergence is exact and does not depend on geometric transformations, as proven in (\ref{CD:PullBack}).
\end{example}

To summarize, with the appropriate basis functions in which to represent continuous differentiation, $ \Dform{\beta}{k}_h(\vect{x}{}) = \diff \; \Dform{\alpha}{k-1}_h(\vect{x}{})$ ($k=1,2,3$), the gradient, curl and divergence reduce to a relation between the expansion coefficients.
\section{Mixed formulation for Stokes flow in terms of differential forms \label{sec:MixedGalerkin}}
In this section we will use the finite dimensional spaces of differential forms, developed in the previous section, in a mixed continuous Galerkin setting. We will apply the method to several 2D test problems concerning Stokes flow (see section \ref{sec:Introduction}). Which mixed formulation we choose to attack the problem, will solely dependent on physical considerations. Following the reasoning in section \ref{sec:Introduction}, we choose all variables as outer oriented differential forms, such that we can apply the topological divergence for mass conservation, see Figure \ref{diag:OuterOrientation}. This ultimately results in a point-wise divergence free velocity field in the discrete setting. Furthermore, choosing all variables as outer oriented differential forms, has the advantage that normal velocity is prescribed in a strong sense, while tangential velocity is prescribed in a weak sense.

Consider an $n$-dimensional domain $\Omega$, filled with an incompressible fluid with constant viscosity $\nu$. Under these assumptions, Stokes flow can be described in terms of differential forms as, 
\begin{align}
	& \diff^{\star} \Dform{u}{n-1} = \Dform{\omega}{n-2},& &\diff \; \Dform{\omega}{n-2} + \diff^{\star} \; \Dform{p}{n} = \Dform{f}{n-1},& &\diff \; \Dform{u}{n-1} = 0&
	\label{eq:Stokes1}
\end{align}
Here velocity is an outer oriented $(n-1)$-form, $\Dform{u}{n-1}$, pressure is an outer oriented $n$-form, $\Dform{p}{n}$, the vorticity is an outer oriented $(n-2)$-form and the right-hand-side is an outer oriented $(n-1)$-form, $\Dform{f}{n-1}$. 

We follow the same reasoning as in \cite{Kreeft:2012a} and derive the mixed formulation of the Stokes problem in 2 steps,
\begin{enumerate}
	\item Multiply the equations in (\ref{eq:Stokes1}) by test functions $\Dform{\alpha}{n-2}$, $\Dform{\beta}{n-1}$ and  $\Dform{\gamma}{n}$ to obtain $L^2$ inner products terms.
	\item Use integration by parts to replace the coderivative $\diff^{\star}$ by the exterior derivative $\diff$ and an additional boundary integral using (\ref{eq:coderivative}).
\end{enumerate}

The Stokes problem can now be posed as follows: find $\left\{\Dform{\omega}{n-2} \in \Lambda^{n-2},\Dform{u}{n-1} \in \Lambda^{n-1},\Dform{p}{n} \in \Lambda^{n} \right\}$ for all $\left\{\Dform{\alpha}{n-2} \in \Lambda^{n-2}\right.$, $\Dform{\beta}{n-1} \in \Lambda^{n-1}$, $\left.\Dform{\gamma}{n} \in \Lambda^{n} \right\}$, such that 
\begin{align}
	\left(\diff \; \Dform{\alpha}{n-2}, \Dform{u}{n-1} \right)_{\Omega} - \left(\Dform{\alpha}{n-2},\Dform{\omega}{n-2}\right)_{\Omega} &= 
			\int_{\partial \Omega} \Dform{\alpha}{n-2} \wedge \star \Dform{u}{n-1}  \\
	\left(\Dform{\beta}{n-1}, \diff \; \Dform{\omega}{n-2} \right)_{\Omega} + \left(\diff \; \Dform{\beta}{n-1}, \Dform{p}{n} \right)_{\Omega} 
			- \int_{\partial \Omega} \Dform{\beta}{n-1} \wedge \star \Dform{p}{n} &= \left(\Dform{\beta}{n-1}, \Dform{f}{n-1} \right)_{\Omega} \\
	\left(\Dform{\gamma}{n}, \diff \; \Dform{u}{n-1} \right)_{\Omega} &= 0.
\end{align}
This is the vorticity-velocity-pressure (VVP) formulation also used in \cite{Bochev:2010,Kreeft:2012a,Kreeft:2012b,Dubois:2002} which is well posed \cite{Kreeft:2012b,Girault:1986}. In the discrete problem we merely replace the infinite dimensional spaces of differential forms by $\left\{\Dform{\omega_h}{n-2} \in \Lambda_h^{n-2} \right.$, $\Dform{u_h}{n-1} \in \Lambda_h^{n-1}$, $\left.\Dform{p_h}{n} \in \Lambda_h^{n} \right\}$. Because, $\Lambda_h(\Omega) \subset \Lambda(\Omega)$, the discrete problem is automatically well posed as well. The problem is closed by imposing normal and tangential velocity components at the boundary. See \cite{Kreeft:2012a,Kreeft:2012b} for an overview of all admissible boundary conditions. The resulting system is symmetric and given by,
\begin{align}
	\begin{pmatrix}  
		-\Mat{M}{n-2} 						& \Mat{D}{n-1,n-2}^T \Mat{M}{n-1} &  \oslash 								\\
		 \Mat{M}{n-1} \Mat{D}{n-1,n-2} 	& \oslash 												&  \Mat{D}{n,n-1}^T \Mat{M}{n} - \Mat{B}{2}(\vect{p}{})\\
		 \oslash                  & \Mat{M}{n} \Mat{D}{n,n-1}      	&  \oslash
	\end{pmatrix}
	\begin{pmatrix}  
		\vect{\omega}{} \\ \vect{u}{} \\ \vect{p}{} 
	\end{pmatrix}
	=
	\begin{pmatrix}  
		\Mat{B}{1}(\Dform{u}{n-1})  \\ \Mat{M}{n-1}(\Dform{f}{n-1}) \\ \oslash 
	\end{pmatrix}
\end{align}
The discrete curl $\Mat{D}{n-1,n-2} = \Mat{E}{n-2,n-1}^T$ and divergence $ \Mat{D}{n,n-1} = \Mat{E}{n-1,n}^T$ are incidence matrices, see Section \ref{sec:DiscreteModeling}, depending solely on mesh topology and are exact irrespective of the geometry mapping $\Phi$ and coarseness of the mesh. The inner product matrices $\Mat{M}{n-2} $, $\Mat{M}{n-1}$ and $\Mat{M}{n}$ are calculated by pulling the finite dimensional spaces of differential forms back to the reference domain $\Omega'$ (section \ref{sec:DifferentialModeling}) and performing Gauss numerical integration there. $\Mat{B}{1}(\Dform{u}{n-1})$ and $\Mat{B}{2}(\vect{p}{})$ represent the natural boundary terms, tangential velocity and tangential pressure, that arise due to integration by parts using (\ref{eq:coderivative}). Since the tangential velocity is weakly enforced, it appears in the right hand side, while the pressure term $\Mat{B}{2}$ is not enforced and thus remains in the left hand side.

\section{Numerical results \label{sec:NumericalResults}}
We will test the method on several common numerical test cases: 1) We compare with an analytical solution (same test problem as in \cite{Kreeft:2012a}) and show results for $h$-refinement on a Cartesian and skewed, non-linear mesh; 2) Taylor Couette flow and compare with the analytical solution; and 3) lid driven cavity flow and compare with benchmark results of \cite{Sahin:2003}. 

\subsection{Manufactured solution}
\begin{figure}
		\centering
		\subfigure[grid 1]{\includegraphics[width=0.3\textwidth]{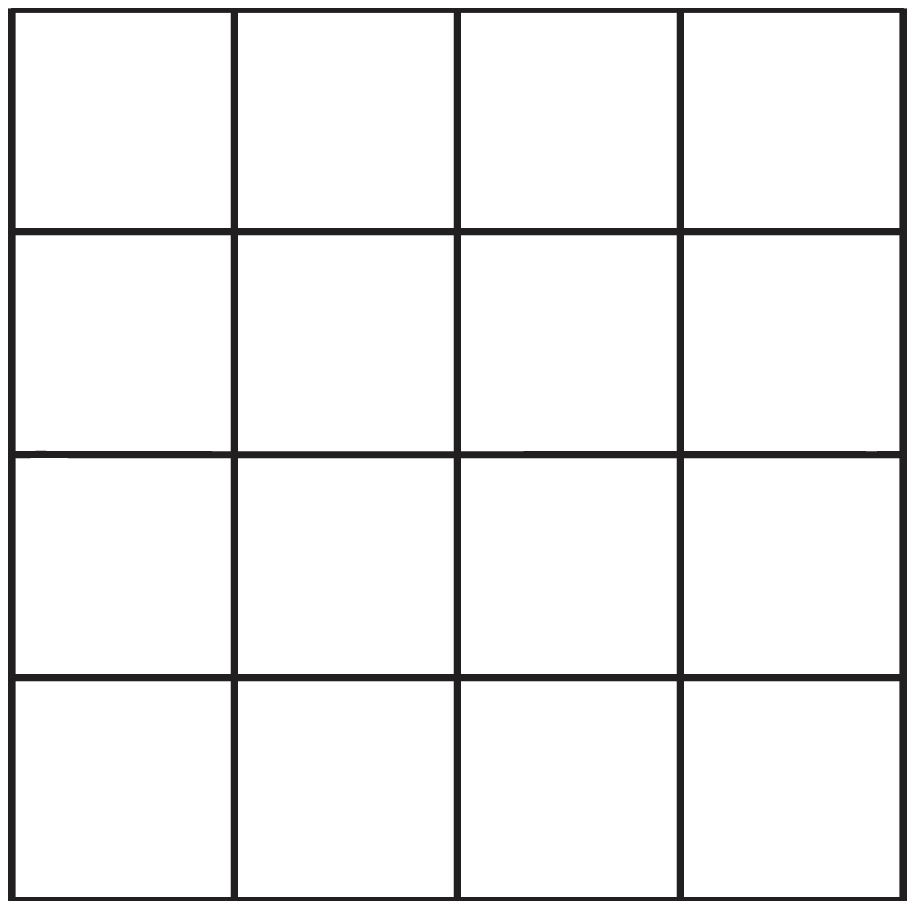}}
		\subfigure[grid 2]{\includegraphics[width=0.3\textwidth]{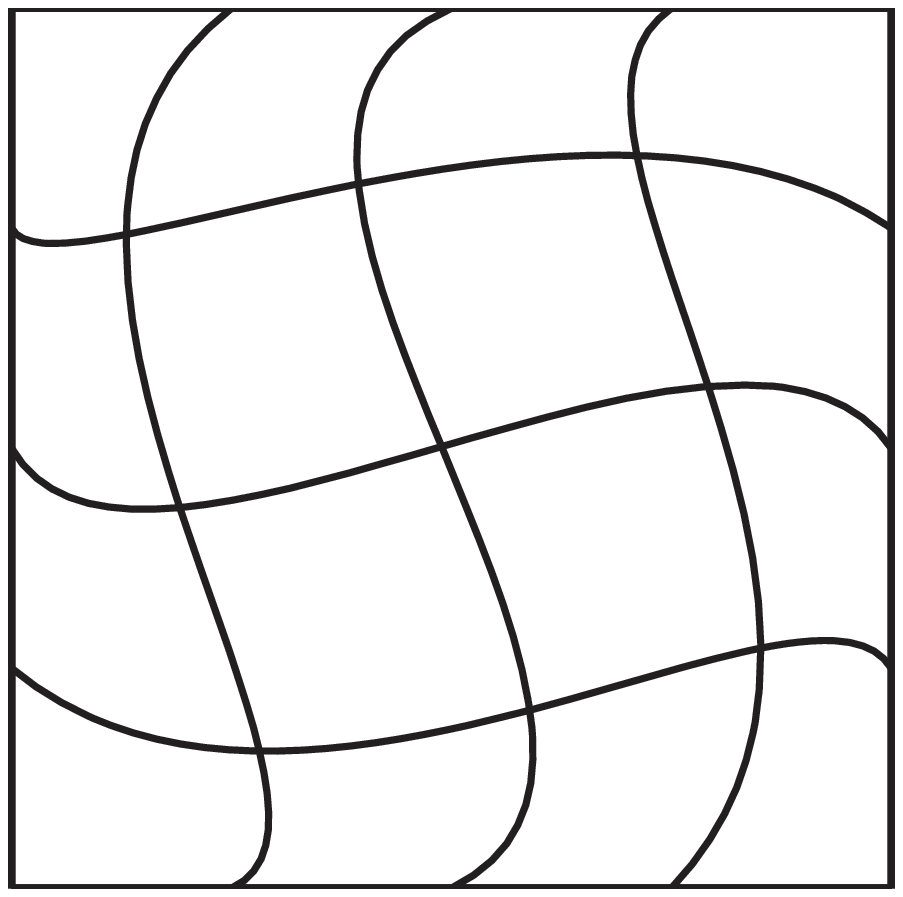}}
		\caption{The Cartesian and curved mesh used in the convergence study.}
		\label{MS1}
\end{figure}

Consider as computational domain the unit square, $\Omega = [0,1]^2$, filled with an incompressible fluid of viscosity, $\nu = 1$. As an analytical solution we choose,
\begin{align*}
 	\Dform{\omega}{0}(x,y) = & -4 \pi \sin{(2\pi x)}\sin{(2\pi y)}.			\\
 	\Dform{u}{1}(x,y) = &-\sin{(2\pi x)} \cos{(2\pi y)} \cdot dy   -\cos{(2\pi x)} \sin{(2\pi y)} \cdot dx \\
	\Dform{p}{2}(x,y) = &\sin{(\pi x)}\sin{(\pi y)} \cdot dx \; dy. 		\\
		\Dform{f}{1}(x,y) = 
		& \left(-8\pi^2 \sin{(2\pi x)} \cos{(2\pi y)} + \pi \cos{(\pi x)} \sin{(\pi y)} \right) \cdot dy + \\ 
		& \left(-8\pi^2 \cos{(2\pi x)} \sin{(2\pi y)} - \pi \sin{(\pi x)} \cos{(\pi y)} \right) \cdot dx. 
\end{align*}
We perform calculations on a Cartesian mesh as well as on a bicubic degree curved mesh, see Figure \ref{MS1}.

\begin{figure}
		\centering
%		\subfigure[grid 1]{\includegraphics[width=0.3\textwidth]{figures/Grid1}}
%		\subfigure[grid 2]{\includegraphics[width=0.3\textwidth]{figures/Grid2}}
%		\caption{\label{MS1} The Cartesian and curved mesh used in the convergence study.}
		\subfigure[Vorticity - grid 1]{\includegraphics[width=0.33\textwidth]{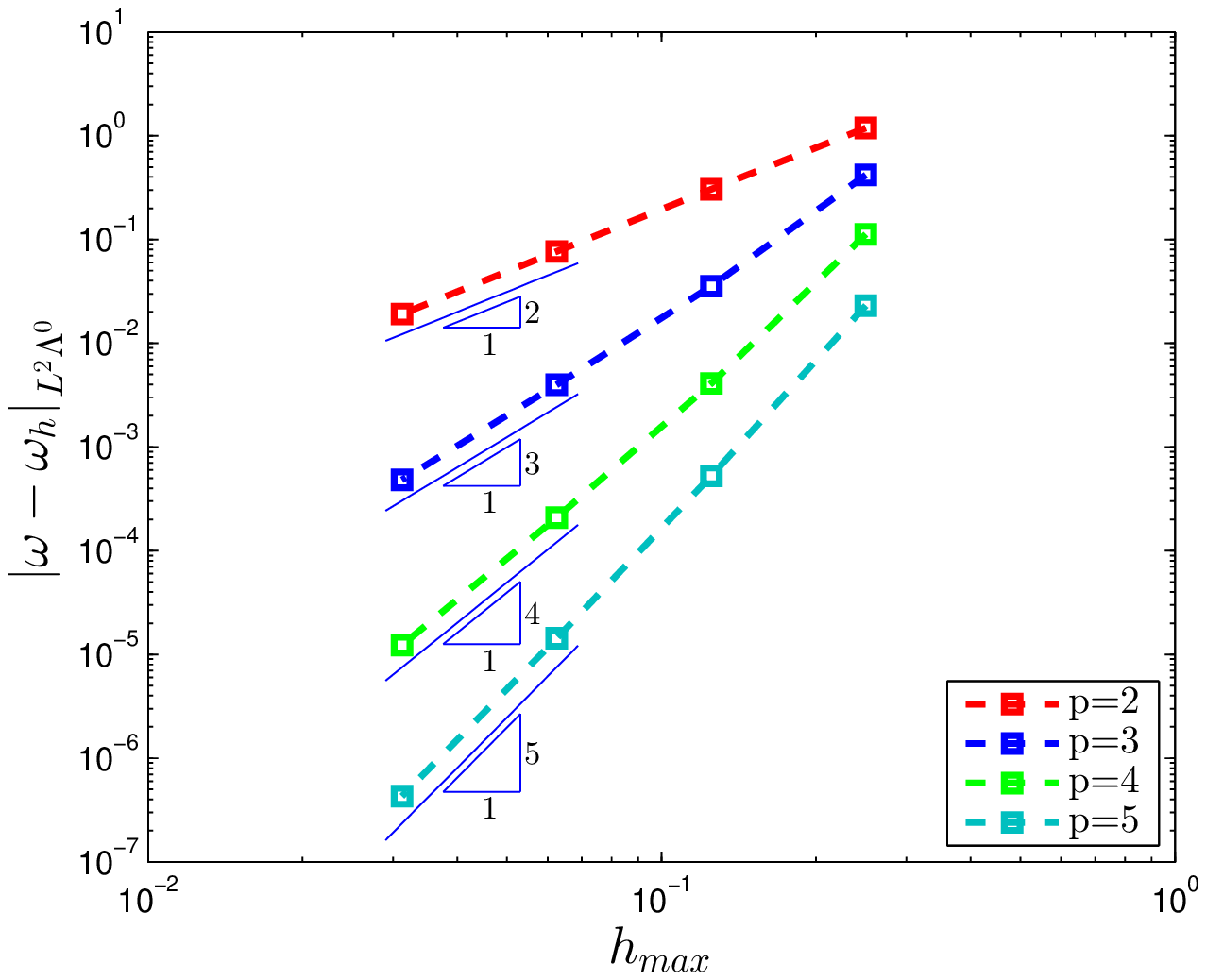}}
		\subfigure[Velocity - grid 1]{\includegraphics[width=0.33\textwidth]{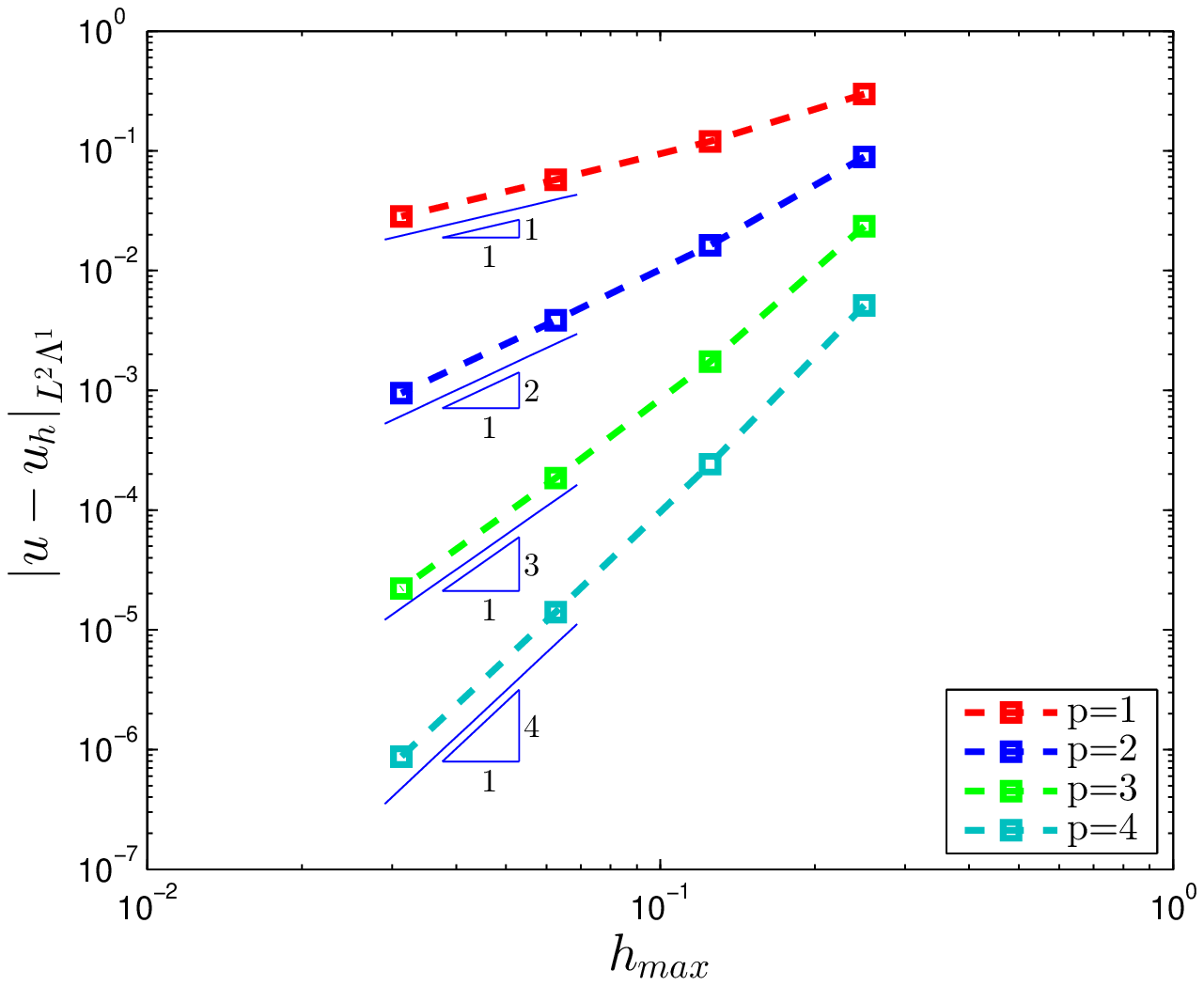}}
		\subfigure[Pressure - grid 1]{\includegraphics[width=0.33\textwidth]{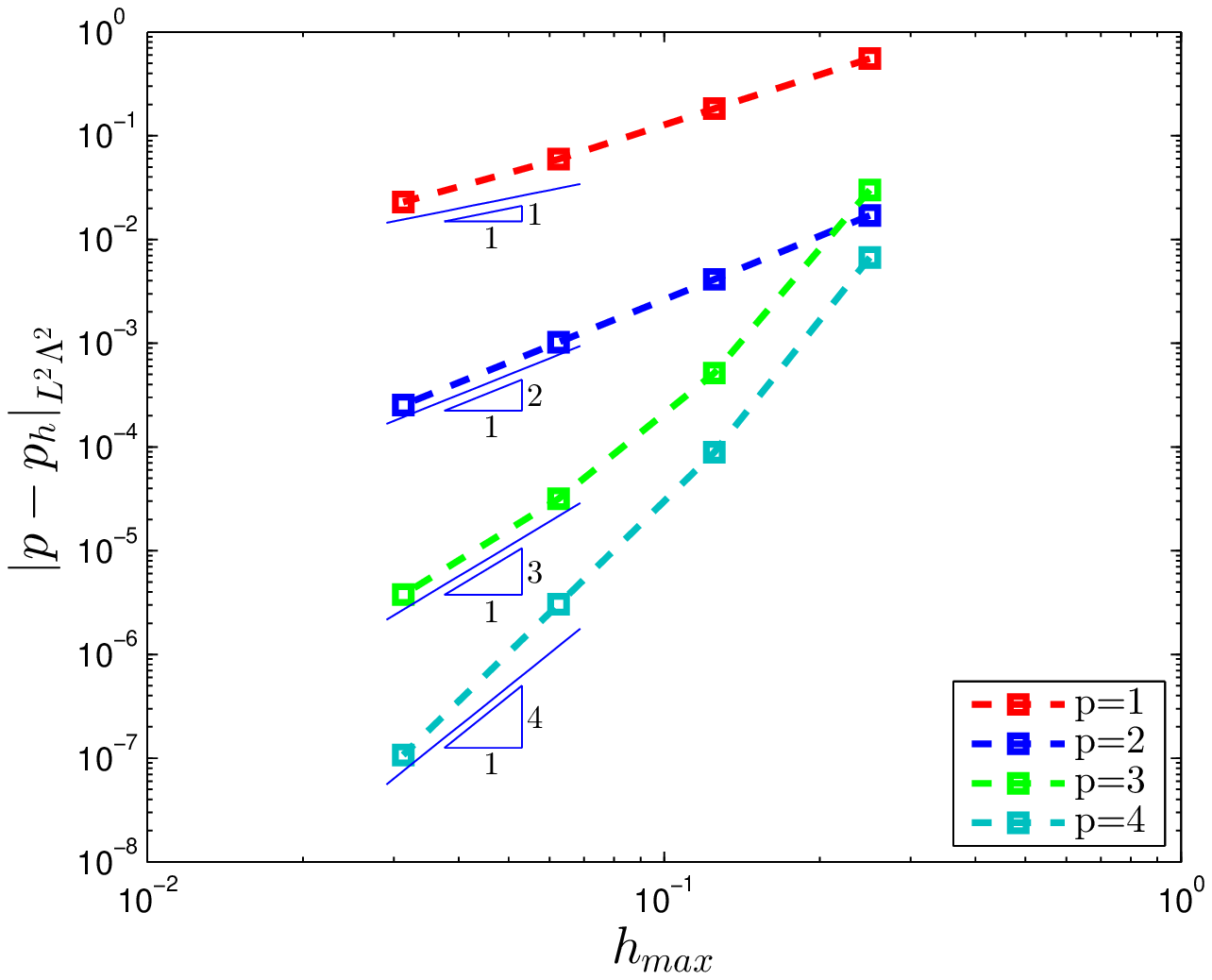}} \\
		\subfigure[Vorticity - grid 2]{\includegraphics[width=0.33\textwidth]{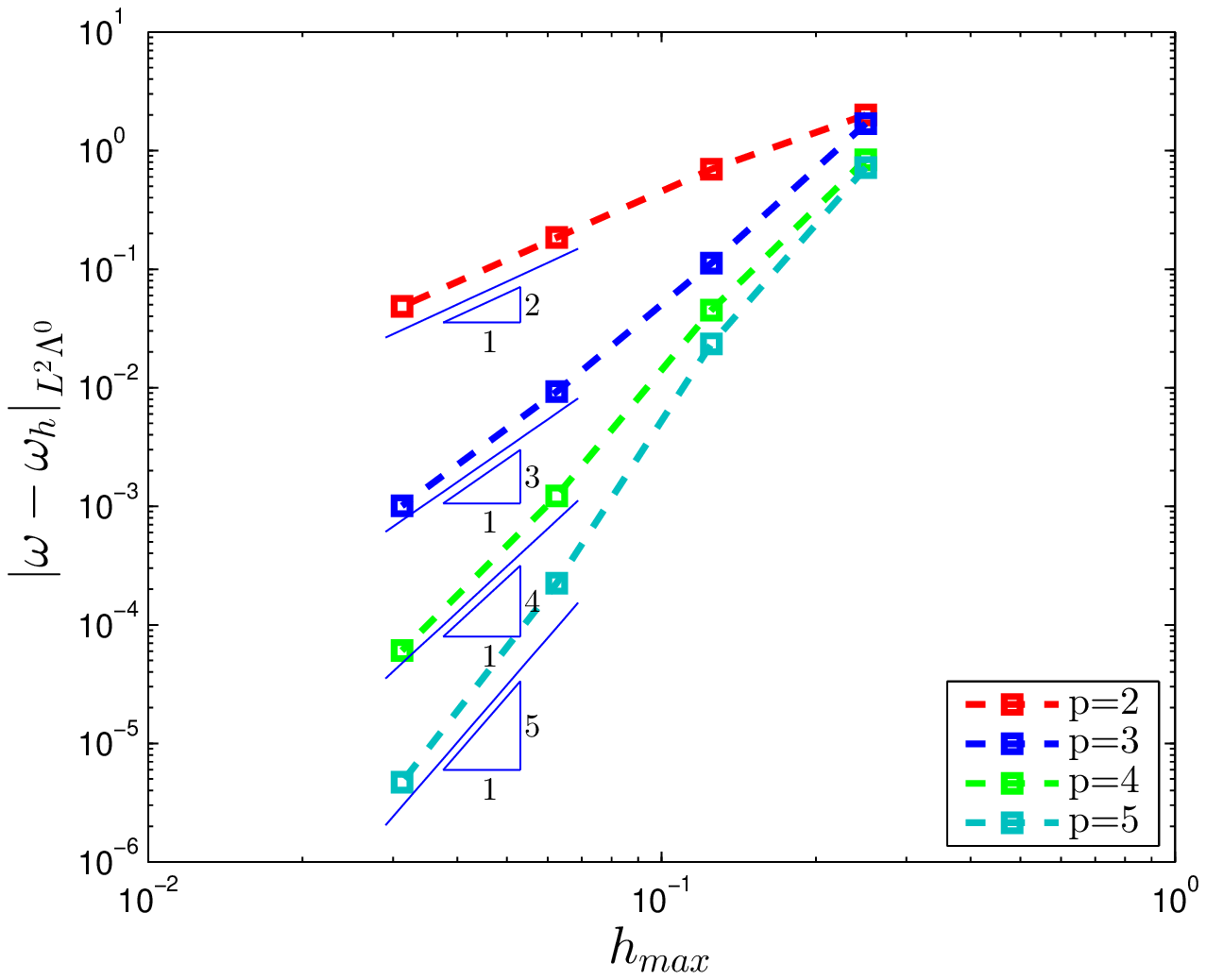}}
		\subfigure[Velocity - grid 2]{\includegraphics[width=0.33\textwidth]{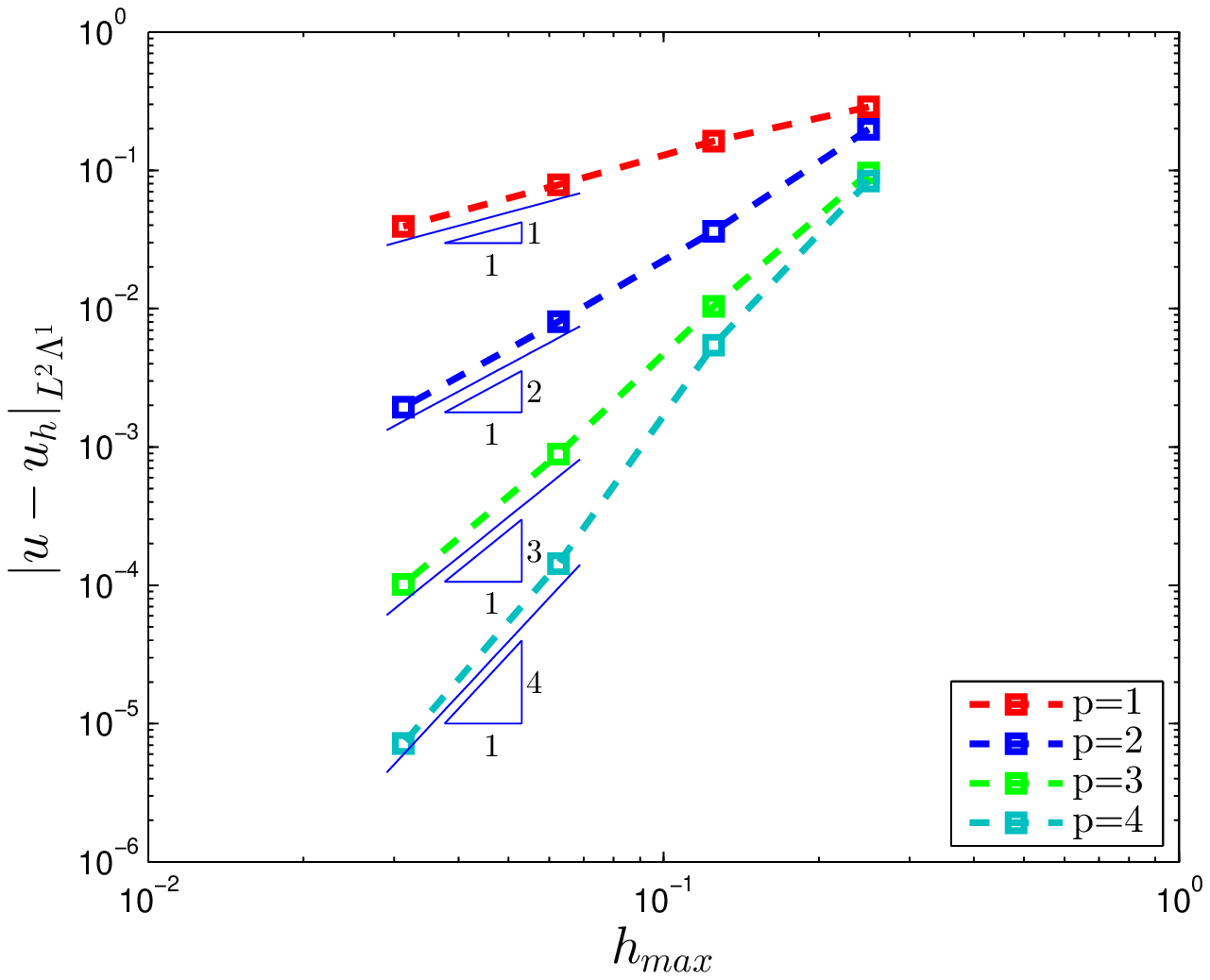}}
		\subfigure[Pressure - grid 2]{\includegraphics[width=0.33\textwidth]{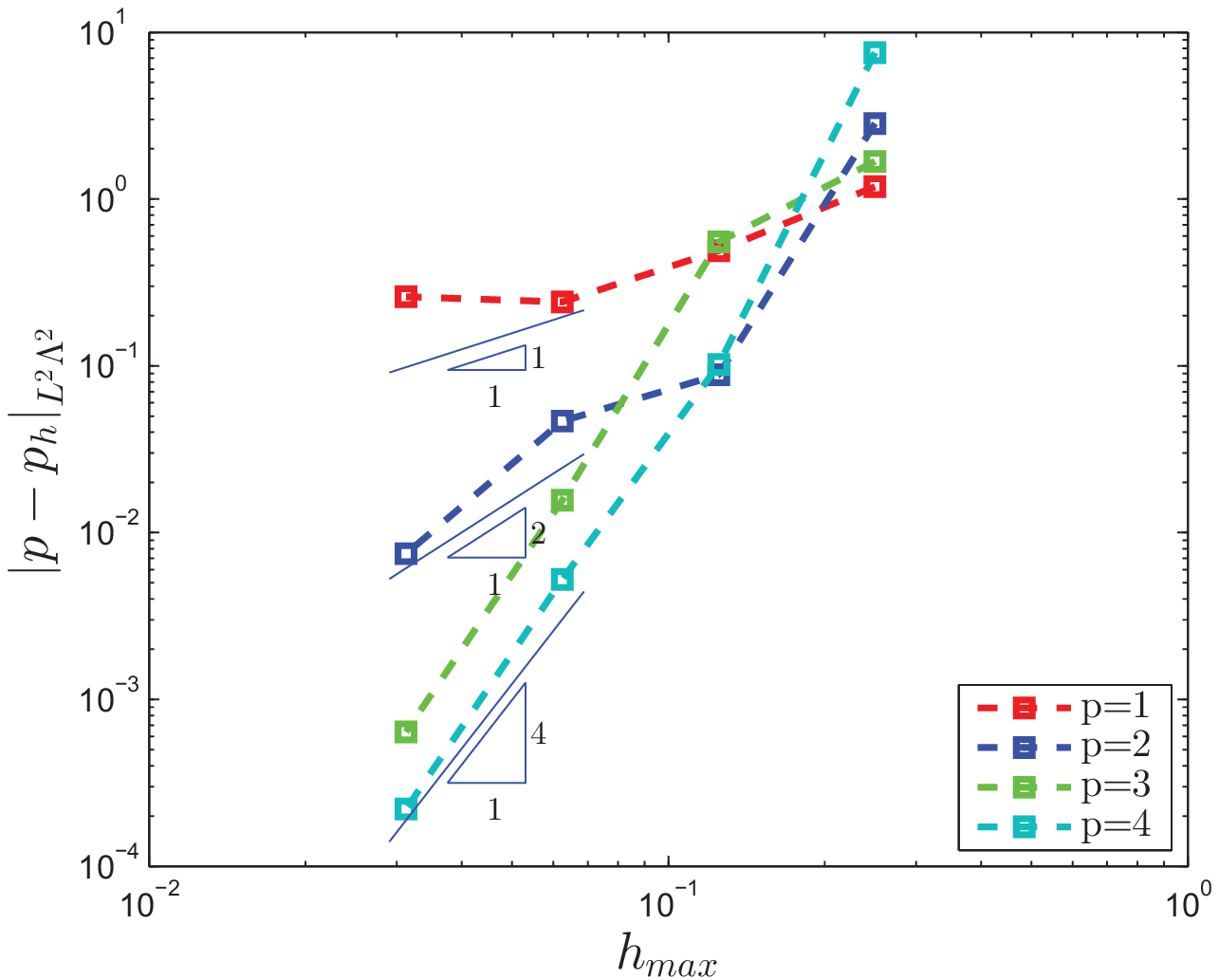}}
		\caption{\label{MS2} H-convergence for vorticity, velocity and pressure on grid 1 (a,b,c) and 2 (d,e,f).}
\end{figure}

The divergence of $\Dform{u}{1}_h$ is point wise zero at all stages of refinement. Figure \ref{MS2} shows the convergence behavior under $h$-refinement of the vorticity, velocity and pressure for polynomial orders ranging from $p=2$ to $p=5$ for the vorticity and $p=1$ to $p=4$ for the velocity and pressure. The maximum mesh-size $h_{max}$ is calculated as the maximum diagonal length of all elements divided by the square root of 2. We can observe optimal order of convergence for all variables on both grids. Only the pressure on grid 2 of polynomial order $p=1$ shows some odd behavior. This is probably due to the low order of approximation of the variables in conjunction with a higher order ($p=4$) non linear mapping of the geometry (super-parametric). 

\subsection{Taylor Couette flow}
We consider a somewhat more realistic problem with an incompressible fluid with viscosity, $\nu = 1$, in a domain bounded by two contra rotating cylinders. The outer cylinder with radius $r_{out} = 2$ is fixed, while the inner cylinder with radius $r_{in} = 1$ rotates in counter clockwise direction with an angular velocity equal to 1.

\begin{figure}
		\centering
		\subfigure[velocity field and NURBS geometry]{\includegraphics[trim = 10cm 4cm 10cm 4cm,clip,height=0.25\textheight]{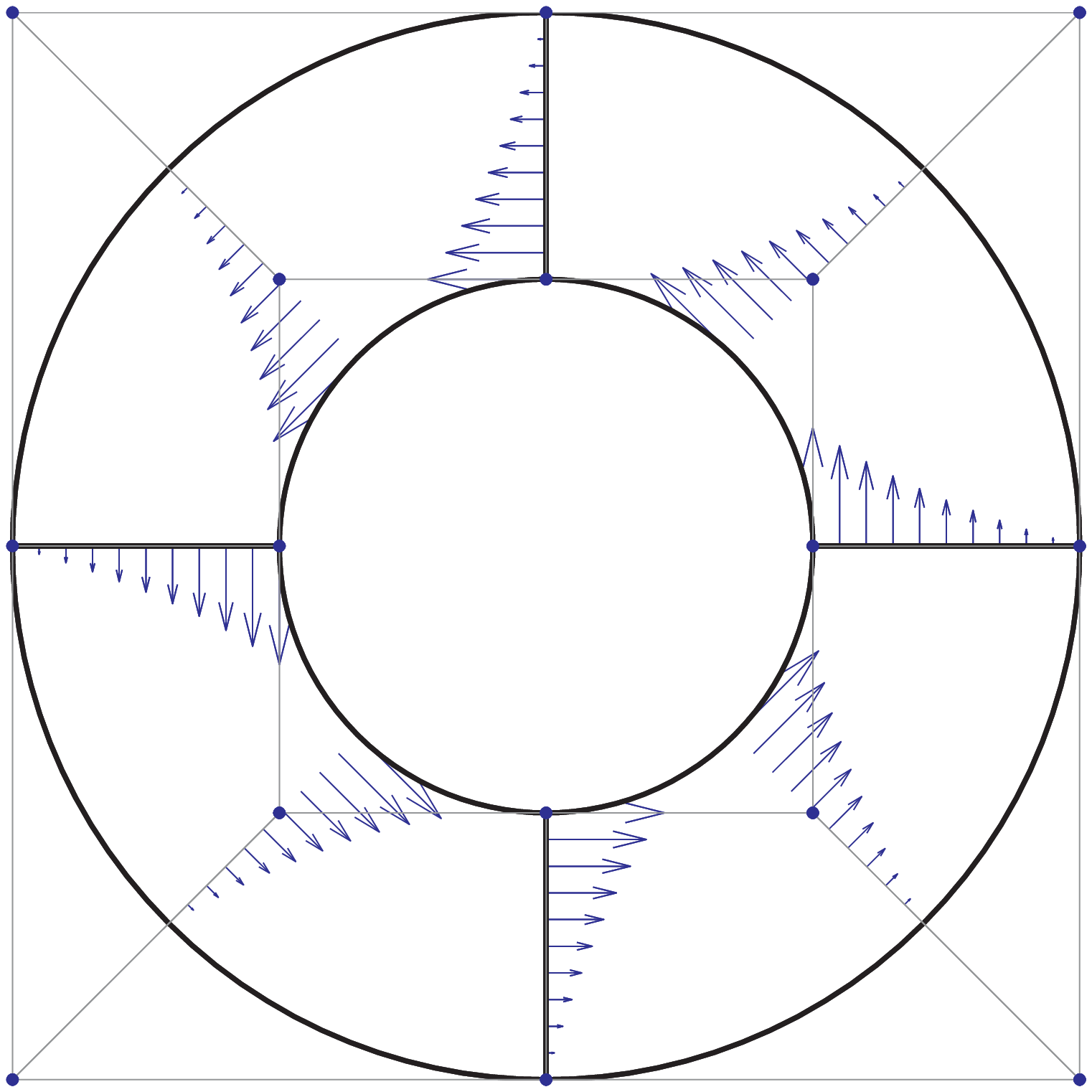}}
		\subfigure[$h$-convergence velocity]{\includegraphics[height=0.25\textheight]{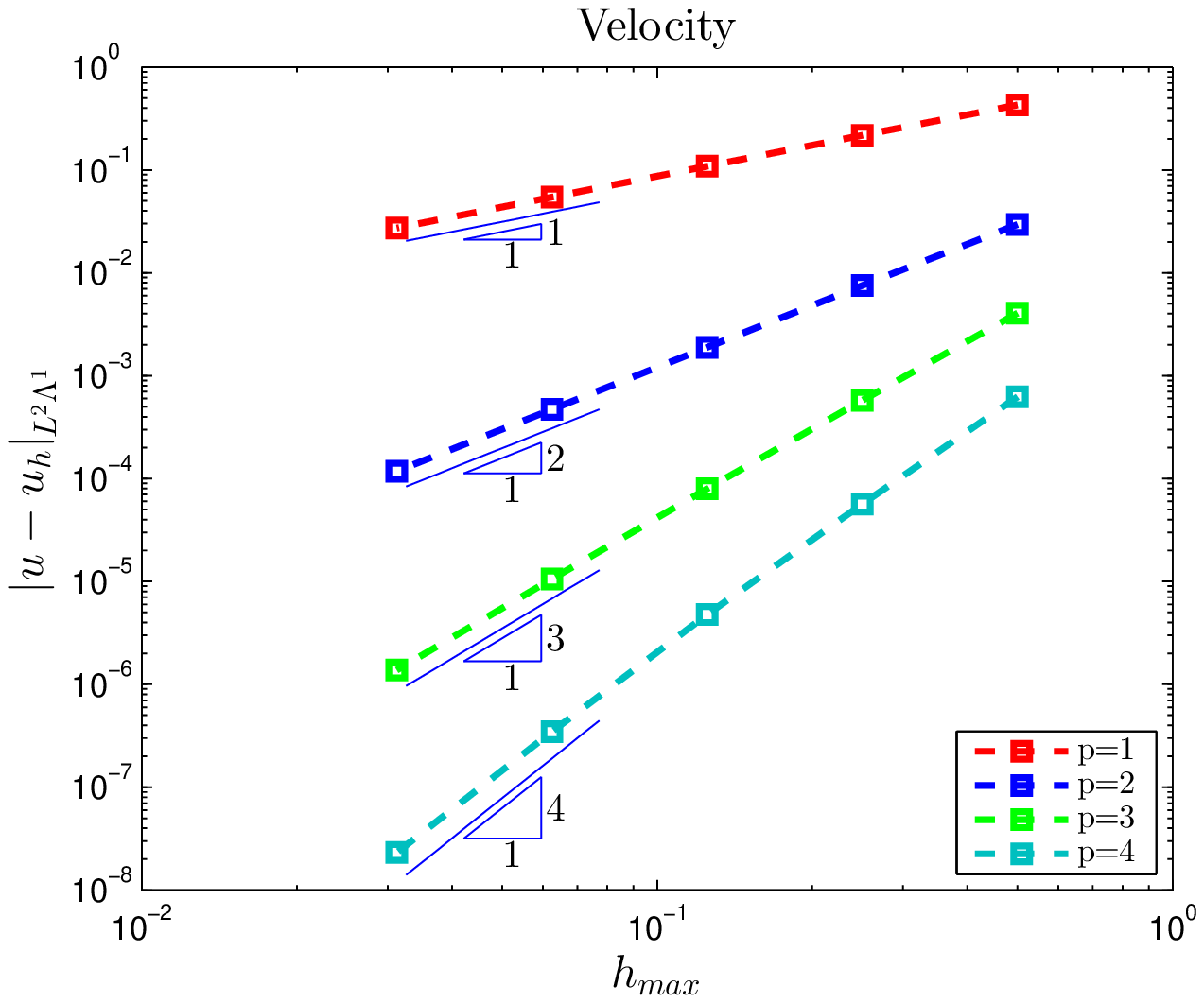}}
		\caption{Flow field and NURBS geometry representation (a) using 4 $C^0$ patches; and $h$-convergence for velocity (b).}
		\label{TCF}
\end{figure}

The flow field in Cartesian coordinates (x,y) as a function of radius $r$ is then described by the following velocity field,
\begin{align}
	\Dform{u}{1}(r) = \left(-\frac{1}{3} r + \frac{4}{3r}  \right) \left( -\sin{\theta} dy + \cos{\theta} dx \right) 
\end{align}
and is depicted in Figure \ref{TCF}a. 

The geometry is exactly represented by 4 $C^0$ NURBS patches of degree 2 by 1. Note that this test-case is more interesting since not all weights are equal to one. The weights associated with the four inner and four outer corners of the domain have weights equal to $0.5 \cdot \sqrt{2}$. Figure \ref{TCF}b shows optimal convergence of the velocity under $h$-refinement for polynomial orders of $p=1$ to $p=4$ in the radial direction. The pressure is constant throughout the domain and is up to machine precision at all stages of refinement. Here $h_{max}$  represents the mesh-size in radial direction. 

\subsection{Lid driven cavity flow}
The lid driven cavity flow is one of the classical benchmark cases to asses numerical methods and to verify Navier Stokes codes. Consider the unit square domain $\Omega = [0,1]^2$ filled with an incompressible fluid of viscosity $\nu = 1$. On the top of the domain we apply a tangential velocity $u(x,1) = 1$ while on the other sides of the domain the tangential velocity is set to zero. The result is a clockwise rotating flow with small counter rotating eddies in the two lower corners. Because of the discontinuity of the velocity in the two upper corners, both the vorticity and pressure are infinite at these places. These singularities make the lid driven cavity flow a challenging test case.

\begin{figure}
		\centering
		\subfigure[streamfunction]{\includegraphics[width=0.32\textwidth]{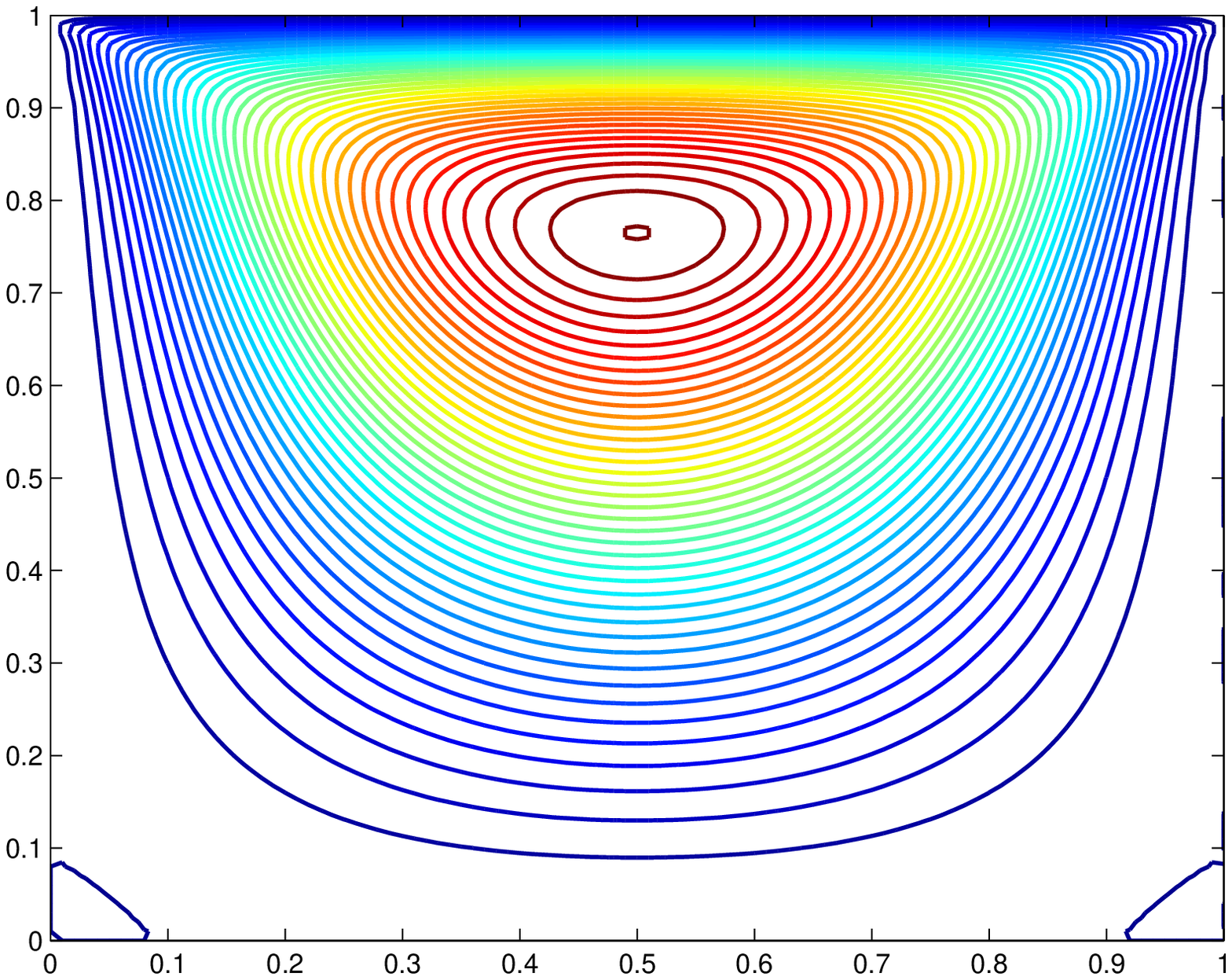}}
		\subfigure[vorticity]{\includegraphics[width=0.32\textwidth]{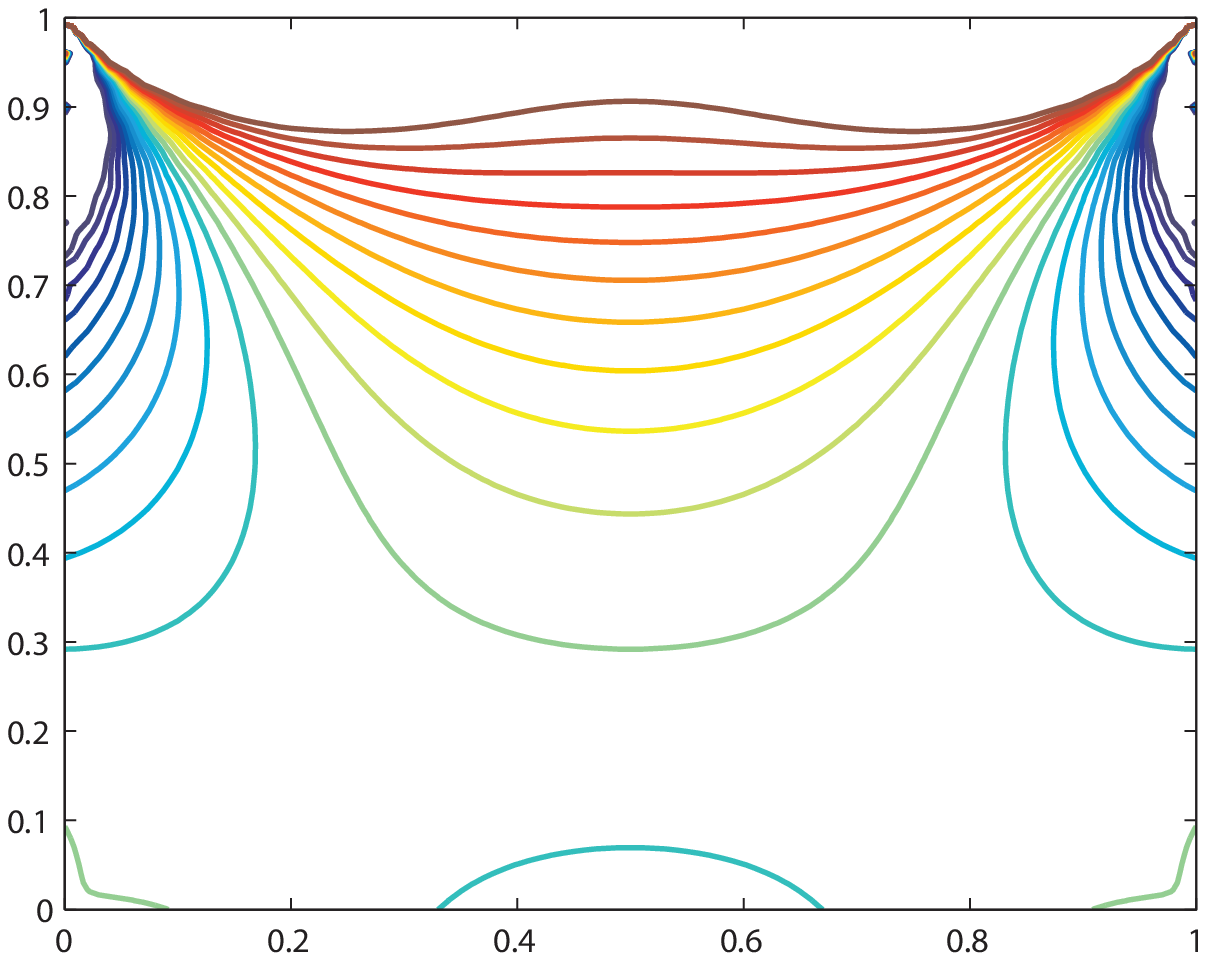}}
		\subfigure[pressure]{\includegraphics[width=0.32\textwidth]{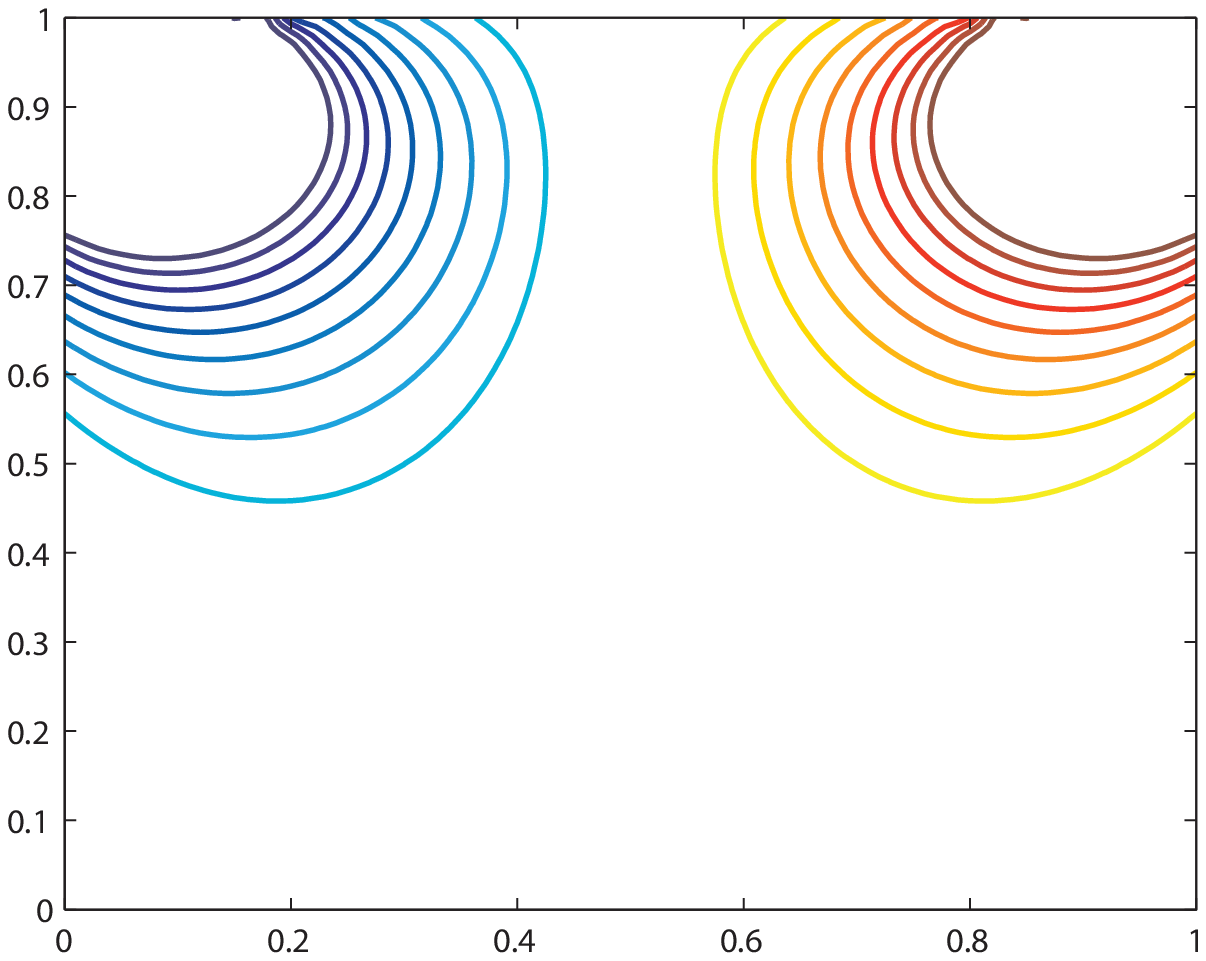}}
		\caption{Results lid driven cavity flow on a 60x60 uniform bi-cubic NURBS mesh (weights are all equal to one) of highest regularity.}
		\label{LDC1}
\end{figure}

Figure \ref{LDC1} shows results for the stream function, vorticity and pressure on a bi-cubic uniform grid of maximum regularity and $60x60$ degrees of freedom. We note that the divergence of velocity is point-wise zero in whole $\Omega$. Furthermore no special treatment has been given to the corner singularities.

\begin{figure}
		\centering
		\subfigure[]{\includegraphics[width=0.32\textwidth]{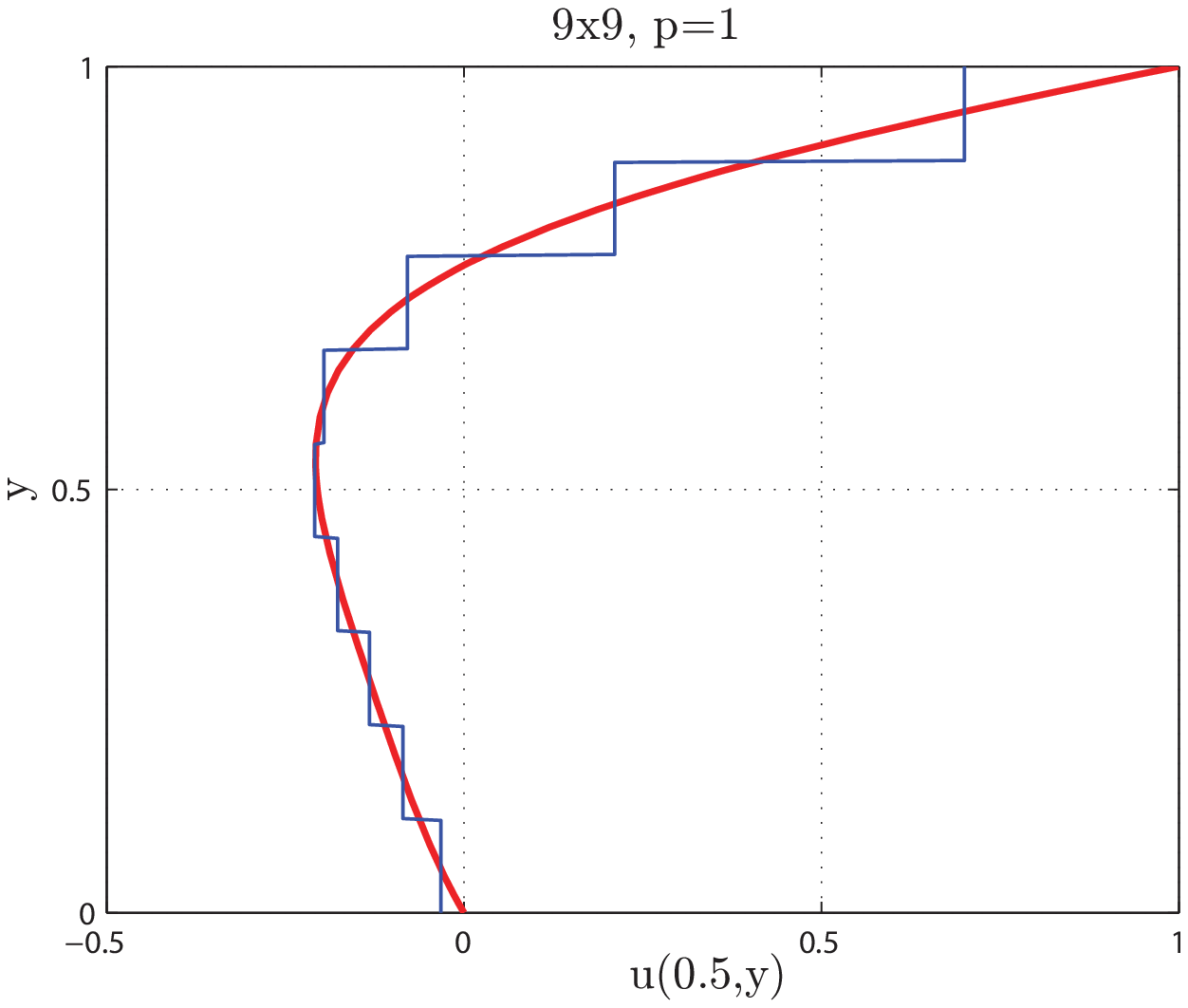}}
		\subfigure[]{\includegraphics[width=0.32\textwidth]{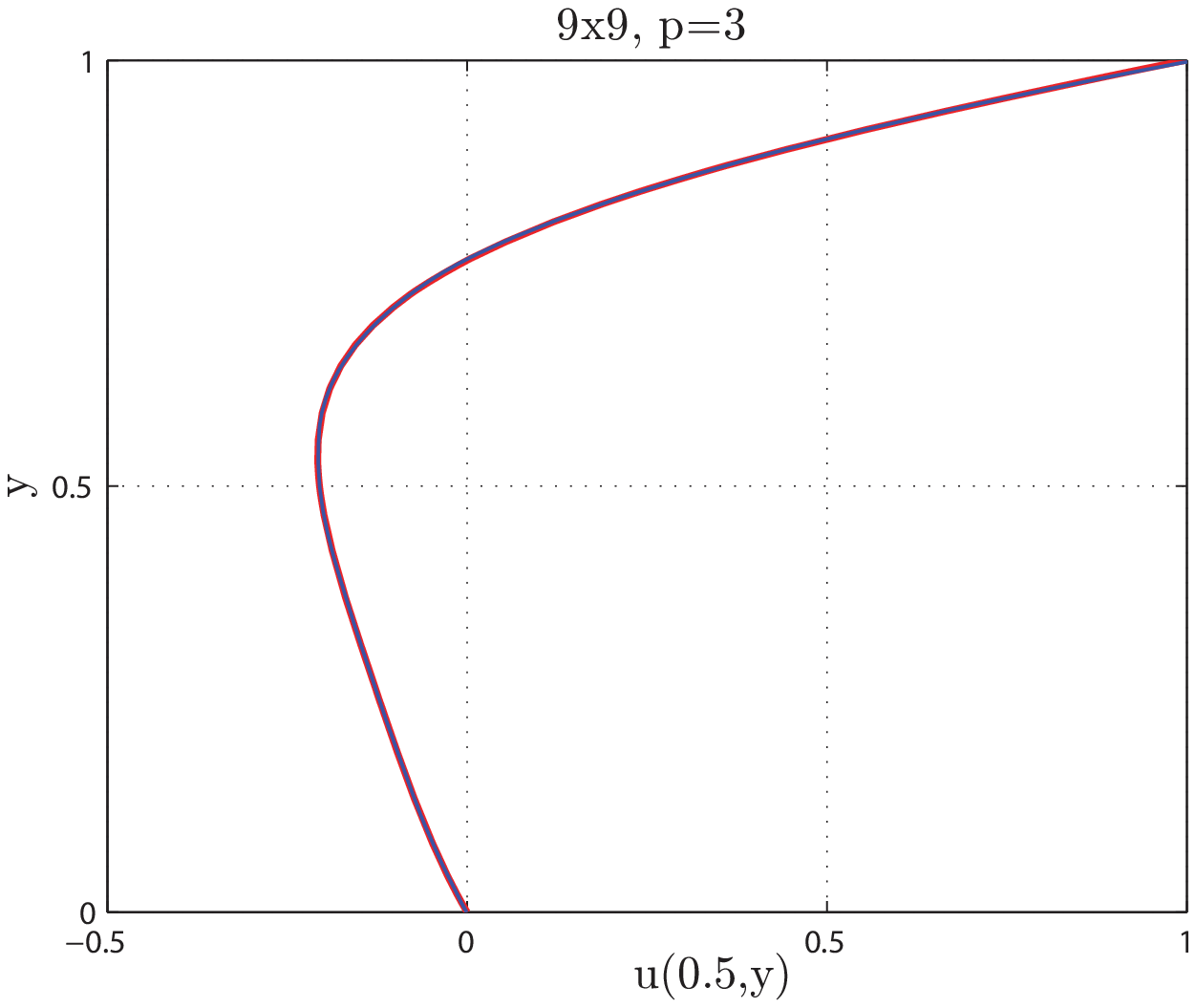}}
		\subfigure[]{\includegraphics[width=0.32\textwidth]{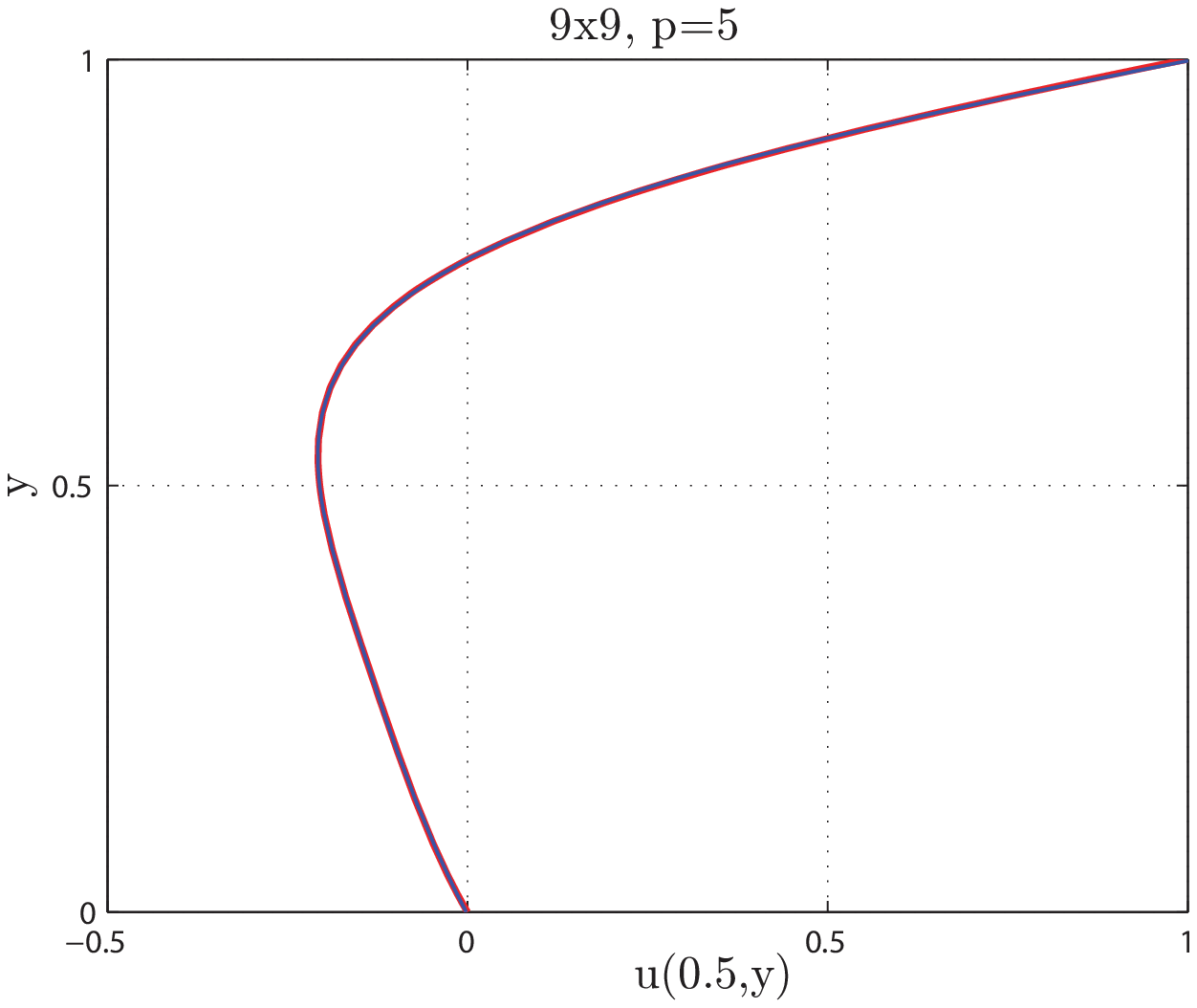}} \\
		\subfigure[]{\includegraphics[width=0.32\textwidth]{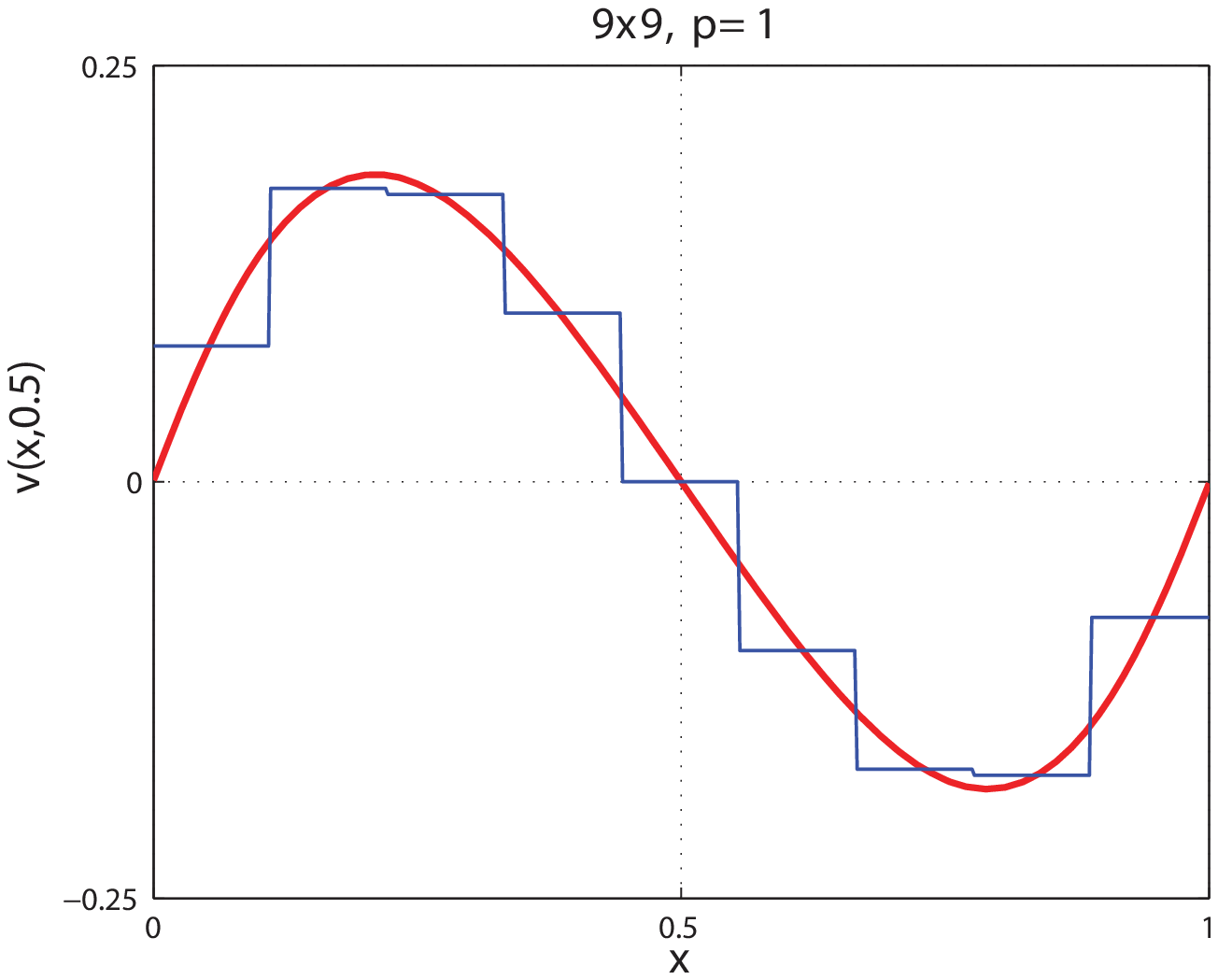}}
		\subfigure[]{\includegraphics[width=0.32\textwidth]{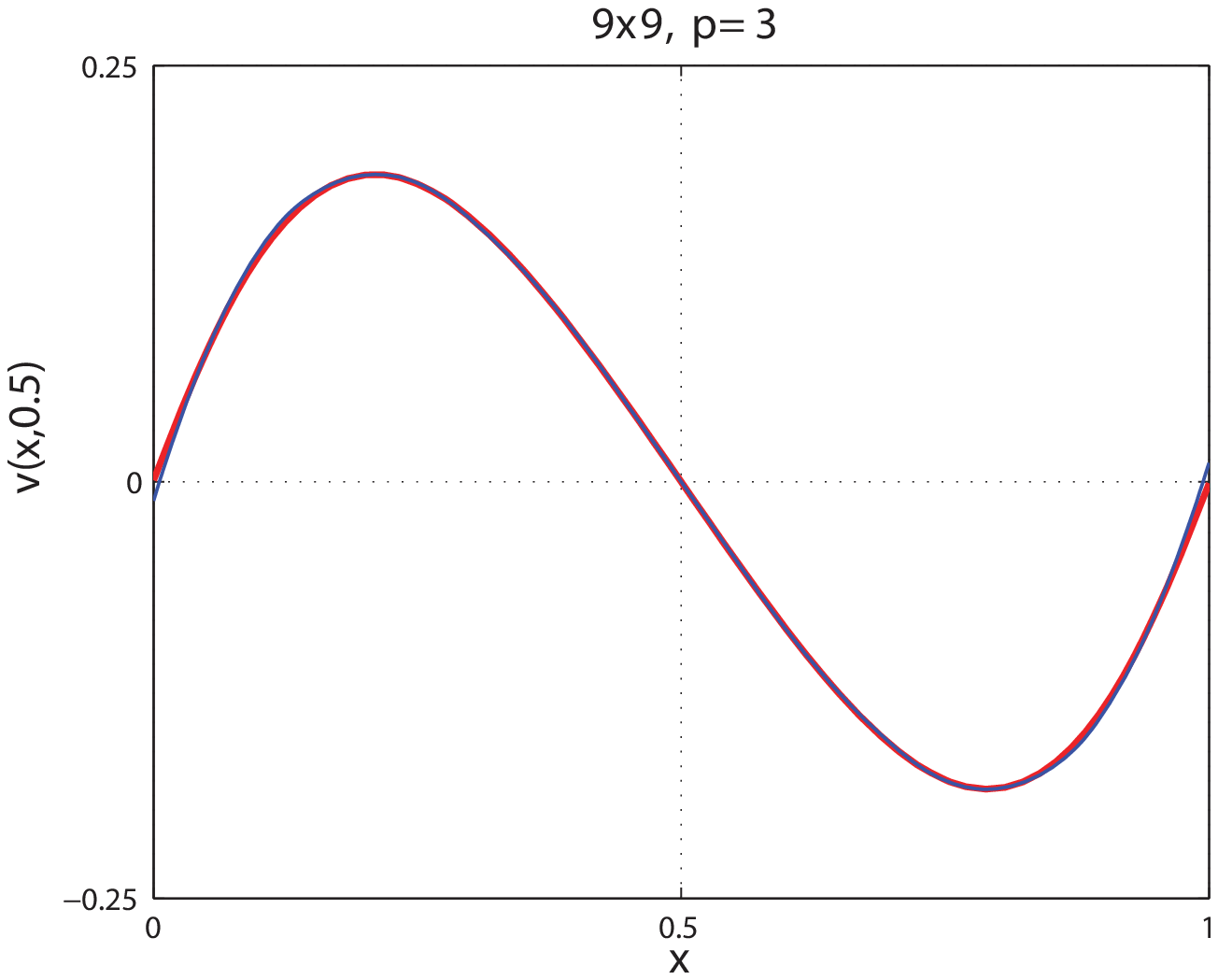}}
		\subfigure[]{\includegraphics[width=0.32\textwidth]{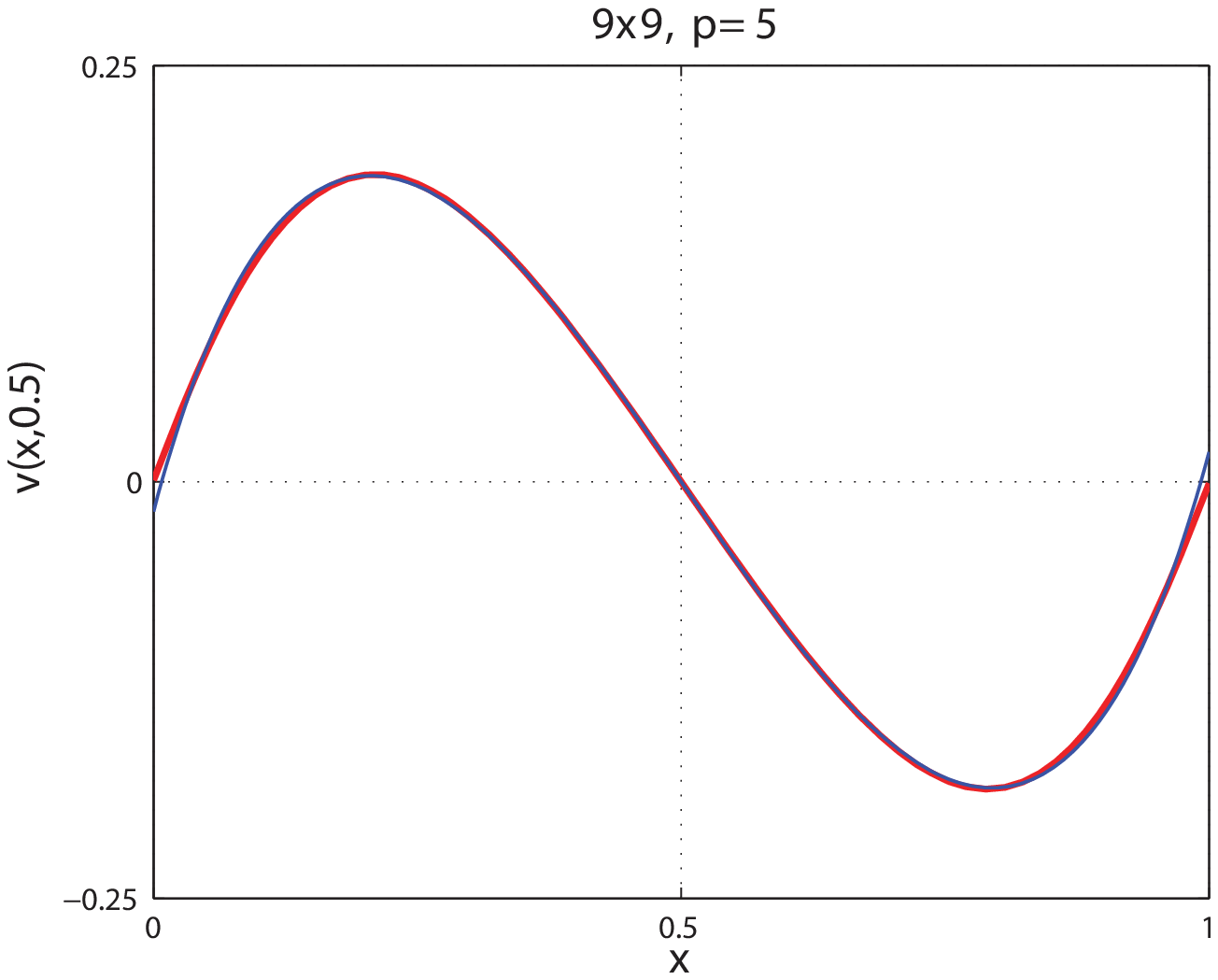}}
		\caption{Comparison of numerical approximation (blue) with benchmark results of \cite{Sahin:2003} (red). Horizontal (a,b,c) and vertical (d,e,f) velocity profile at centerlines of cavity for a 9x9 uniform grid of order 1, 3 and 5.}
		\label{LDC2}
\end{figure}

In Figure \ref{LDC2} the horizontal component of velocity has been plotted along the vertical centerline $(0.5,y)$ and the vertical component of velocity has been plotted along the horizontal centerline $(y,0.5)$. The results confirm very well with the benchmark results of \cite{Sahin:2003} even though the mesh is very coarse (9x9 uniform grid of polynomial order 1, 3 and 5). The most striking result is, however, that of the low order approximation. Note that although there is a clear point-wise difference, the integral values seem to match very well. This is a direct consequence of the conservation properties that we have build into the basis.

\section{Conclusion \label{sec:Conclusion}}
We have developed arbitrary order interpolants, from any basis that is a partition of unity, that satisfies a discrete Stokes theorem. The resulting gradient, curl and divergence conforming spaces have the property that the conservation laws become completely independent of the basis functions. As an example, we have derived conforming spaces form NURBS, and thereby generalized the discrete spaces of differential forms introduced in \cite{Buffa:2011b}. We have applied these new spaces in a mixed Galerkin setting which amongst others resulted in an exactly divergence free discretization of Stokes flow.

\bibliographystyle{elsarticle/elsarticle-num-names}
\bibliography{mybiblio} 

\end{document}